**Raul Oliveira Ribeiro**

# Higher-Order Spectral Element Methods for Electromagnetic Modeling of Complex Anisotropic Waveguides

**Tese de Doutorado**

Thesis presented to the Programa de Pós–graduação em Engenharia Elétrica, do Departamento de Engenharia Elétrica da PUC-Rio in partial fulfillment of the requirements for the degree of Doutor em Engenharia Elétrica.

Advisor     :   Prof. José Ricardo Bergmann
Co-advisor: Prof. Guilherme Simon da Rosa

Rio de Janeiro
April 2024

**Raul Oliveira Ribeiro**

# Higher-Order Spectral Element Methods for Electromagnetic Modeling of Complex Anisotropic Waveguides

Thesis presented to the Programa de Pós–graduação em Engenharia Elétrica da PUC-Rio in partial fulfillment of the requirements for the degree of Doutor em Engenharia Elétrica. Examination Committee:

**Prof. José Ricardo Bergmann**
Advisor
Departamento de Engenharia Elétrica – PUC-Rio

**Prof. Guilherme Simon da Rosa**
Departamento de Engenharia Elétrica – Unesp

**Prof. Fernando Lisboa Teixeira**
ElectroScience Laboratory – The Ohio State University

**Prof. Fernando José da Silva Moreira**
Departamento de Engenharia Elétrica – UFMG

**Prof. Guilherme Penello Temporão**
Departamento de Engenharia Elétrica – PUC-Rio

**Prof. Rafael Abrantes Penchel**
Departamento de Engenharia Elétrica – Unesp

Rio de Janeiro, April the 19th, 2024

# Acknowledgments

First and foremost, I would like to express my gratitude to my advisors, Professor José Ricardo Bergmann and Professor Guilherme Simon da Rosa, for their unwavering support and guidance throughout this research. I would also like to extend my heartfelt appreciation to my overseas advisor, Professor Fernando Lisboa Teixeira, for his invaluable support throughout the research process.

I am also grateful to all the professors who taught me the fascinating aspects of Electromagnetic Theory. I want to thank my colleagues at PUC-Rio and OSU for their friendship and valuable technical discussions.

I want to thank PUC-Rio for exempting school fees offered by the Programa de Pós-Graduação em Engenharia Elétrica.

Lastly, I would like to sincerely thank my family for their endless support.

I must also acknowledge the Brazilian agency CNPq for supporting this work under Grants 140056/2020-3 and 200228/2022-6.

# Abstract

Ribeiro, Raul O.; Bergmann, José R. (Advisor); Rosa, Guilherme S. (Co-Advisor). **Higher-Order Spectral Element Methods for Electromagnetic Modeling of Complex Anisotropic Waveguides**. Rio de Janeiro, 2024. 101p. Tese de Doutorado – Departamento de Engenharia Elétrica, Pontifícia Universidade Católica do Rio de Janeiro.

This research thesis presents a novel higher-order spectral element method (SEM) formulated in cylindrical coordinates for analyzing electromagnetic fields in waveguides filled with complex anisotropic media. In this study, we consider a large class of cylindrical waveguides: radially-bounded and radially-unbounded domains; homogeneous and inhomogeneous waveguides; concentric and non-concentric geometries; Hermitian and non-Hermitian anisotropic media tensors. This work explores different wave equation formulations for one-layer eccentric and multilayer cylindrical waveguides. For the first case, we can define a new normalized scalar Helmholtz equation for decoupling TM and TE modes, and for the second, a vectorial Helmholtz equation for hybrid modes in multilayered anisotropic structures. Additionally, we formulate a transformation optics (TO) framework to include non-symmetric and non-Hermitian media tensors for non-concentric multilayer waveguides. Lastly, we model excitation sources for logging sensors applied in geophysical problems using the fields obtained by SEM. We validate the proposed approach against analytical solutions, perturbation-based and mode-matching-based methods, finite-elements, and finite-integration numerical methods. Our technique obtains accurate results with fewer elements and degrees of freedom (DoF) than Cartesian-based SEM and ordinary finite-element approaches. To this end, we use higher-order two-dimensional basis functions associated with the zeros of the completed Lobatto polynomial to model the fields in each reference element. The convergence analysis demonstrates the absence of the Runge effect as the expansion order increases. Numerical results show that our formulation is efficient and accurate for modeling cylindrical waveguided geometries filled with complex media.

## Resumo

Ribeiro, Raul O.; Bergmann, José R.; Rosa, Guilherme S.. **Método de Elementos Espectrais de Ordem Superior para a Modelagem Eletromagnética de Guias de Ondas Anisotrópicos**. Rio de Janeiro, 2024. 101p. Tese de Doutorado – Departamento de Engenharia Elétrica, Pontifícia Universidade Católica do Rio de Janeiro.

Esta tese apresenta um novo método de elementos espectrais (SEM) de ordem superior formulada em coordenadas cilíndricas para análise de campos eletromagnéticos em guias de onda preenchidos com meios anisotrópicos complexos. Neste estudo, consideramos uma grande classe de guias de onda cilíndricos: domínios radialmente limitados e radialmente ilimitados; guias de onda homogêneos e não homogêneos; geometrias concêntricas e não concêntricas; tensores anisotrópicos hermitianos e não hermitianos. Este trabalho explora diferentes formulações de equações de onda para guias de onda excêntricos de uma camada e cilíndricos multicamadas. Para o primeiro caso, podemos definir uma nova equação escalar normalizada de Helmholtz para desacoplar os modos TM e TE, e para o segundo, uma equação vetorial de Helmholtz para modos híbridos. Além disso, estabelecemos uma transformada óptica (TO) para incluir tensores dos meios não simétricos e não hermitianos para guias de onda multicamadas não concêntricos. Por fim, modelamos fontes de excitação para sensores de perfilagem aplicados em problemas geofísicos utilizando os campos obtidos pela SEM. Validamos nossa abordagem utilizando como referência soluções analíticas, métodos baseados em perturbações e em casamento de modos, e métodos numéricos de elementos finitos e integração finita. Nossa técnica obtém resultados precisos com menos elementos e graus de liberdade (DoF) comparados com a SEM formulada em coordenadas cartesianas e o método dos elementos finitos. Para tanto, utilizamos funções de base bidimensionais de ordem superior associadas aos zeros do polinômio de Lobatto completo para modelar os campos em cada elemento de referência. A análise de convergência demonstra a ausência do efeito Runge à medida que a ordem de expansão aumenta. Os resultados numéricos mostram que nossa formulação é eficiente e precisa para modelar geometrias cilíndricas guiadas preenchidas com meios complexos.

## Palavras-chave

# Table of contents

# List of figures

# List of tables

# List of Abbreviations

AMM – Axial mode-matching

BCS – Bipolar coordinate system

CEM – Computational electromagnetics

CMPM – Cavity-material perturbation method

DoF – Degrees of freedom

FEM – Finite element method

FIT – Finite-integration technique

GAT – Graf's addition theorem

LWD – Logging-while-drilling

MMT – Mode-matching technique

NMM – Numerical mode-matching

PEC – Perfect electric conductor

PML – Perfectly matched layer

RMM – Radial mode-matching

RPM – Regular perturbation method

SEM – Spectral element method

TE – Transverse electric mode

TM – Transverse magnetic mode

TO – Transformation optics

# 1
# Introduction

## 1.1
## General Introduction

The analysis of electromagnetic field propagation in complex media (lossy, inhomogeneous, and anisotropic) is a task of fundamental importance in many engineering applications, such as designing metamaterial-based microwave devices, high-frequency optical circuits, and modeling low-frequency electromagnetic sensors for oil prospecting in geophysical anisotropic formations. The complex shapes of such structures and the non-trivial characteristics of the media require studying novel numerical and analytical approaches for solving Maxwell's equations in a robust and numerically stable manner.

In recent decades, several computational electromagnetics (CEM) methods have been developed or improved to solve new technological problems associated with inhomogeneous and anisotropic waveguides. Predominantly, numerical techniques such as finite-difference, finite-element, and integral-equation methods have gained prominence due to their straightforward numerical implementations. Consequently, these methods have enjoyed widespread adoption, with many commercial CEM software packages incorporating them into their simulation frameworks [1–3]. However, the computational cost regarding memory utilization and CPU processing time can present significant barriers, particularly in large-scale problems relative to the electrical wavelength or abrupt variations in the medium properties.

Metamaterials constructed through artificial media have received significant attention due to their ability to exhibit electromagnetic properties that are not present in conventional materials. Artificial magnetism and negative refractive index are two specific types of behavior demonstrated over the past few years, illustrating the new physics and applications possible when we expand our view of what constitutes a material [4]. Some practical metamaterial realizations involve anisotropic cylindrically-layered structures employing dielectric nanoporous or thin film-layered. These geometries can support propagation at frequencies well below the cutoff, backward-, forward-, and slow-wave propagation at RF and optical frequencies where the metamaterial exhibits exotic constitutive parameters [5–8].

Another approach to these cylindrical constructions involves displacing the center conductor, facilitating various engineering applications, including the creation of filters, impedance-matching tools, and sensing probes [9–12]. Additionally, another significant application involves radial-unbounded cylindrical waveguides, particularly in modeling tunnels and boreholes for geophysical exploration. In such cases, the constitutive parameters may exhibit fully anisotropic characteristics and extreme conductivity variation [13–16].

In this context, the spectral element method (SEM) is an advanced implementation of the finite element method (FEM), wherein the internal nodes of each element are strategically positioned at the zeros of specific families of orthogonal polynomials. SEM offers a notable advantage in achieving high accuracy with a comparatively modest number of elements compared to a standard FEM approach. In recent years, this method has demonstrated considerable success across a diverse array of electromagnetic problems, encompassing applications in photonic crystals, open waveguides, metasurfaces, lithography models, and various truncation boundary conditions [17–21].

## 1.2
## Major Research Contributions

This Ph.D. thesis explores a novel higher-order spectral element method formulated in cylindrical coordinates to model a large class of cylindrical waveguides: radially-bounded and radially-unbounded domains; homogeneous and inhomogeneous filled waveguides; concentric and non-concentric geometries; Hermitian and non-Hermitian anisotropic media tensors. This research thesis explores different wave equation formulations for one-layer eccentric and multilayer cylindrical waveguides. For the first case, we can define a new normalized scalar Helmholtz equation for decoupling TM and TE modes, and for the second, a vectorial Helmholtz equation for hybrid modes in multilayered anisotropic geometries. Additionally, we extend the prior transformation optics (TO) [10, 12] to include non-symmetric and non-Hermitian media tensors for non-concentric multilayer waveguides. Lastly, we model excitation sources for logging sensors applied in hydrocarbon prospection using the fields obtained by SEM. Accordingly, the following topics of study correspond to our original scientific contributions:

– A novel SEM formulation in cylindrical coordinates using generic higher-order two-dimensional basis functions. Our approach requires fewer elements and degrees of freedom (DoF) when compared with the conventional Cartesian-coordinate-based SEM and FEM.

– A normalized scalar variational formulation to model eccentric circular waveguides via conformal transformation optics. The proposed formulation obtains the solution for an entire set of configurations independently of the medium's parameters and, unlike perturbation-based techniques, does not restrict the eccentricity offset and presents excellent accuracy and time results.

– An improved SEM formulation via transformation optics to analyze eccentric two-layer waveguides filled with non-reciprocal media. Our approach required far fewer DoF in comparison to the conventional FEM.

– A new approach for modeling well-logging sensors via the fields obtained by higher-order SEM. Our technique can be applied to model inhomogeneous oil invasions, a task beyond the capabilities of standard mode-matching methods.

## 1.3
## Thesis Organization

All research chapters in this Ph.D. thesis can be read independently and are organized as follows. Chapter 2 presents a novel spectral element method in cylindrical coordinates for analyzing electromagnetic fields in waveguides filled with complex (non-reciprocal and non-Hermitian) anisotropic media. Our method obtains accurate results with a small number of elements and DoFs when compared with the Cartesian-coordinate-based SEM and FEM. For this, we use higher-order expansion functions in which the internal nodes of each element are defined by the zeros of the completed Lobatto polynomial. Numerical results show that our technique is efficient and accurate for modeling cylindrical waveguide geometries filled with complex media. The results of this chapter can also be found in [22, 23].

In Chapter 3, we present a higher-order SEM in cylindrical coordinates to rapidly analyze electromagnetic fields in eccentric coaxial waveguides filled with lossy uniaxially anisotropic media. The formulation utilizes transformation optics to map the eccentric waveguide into an equivalent concentric problem. The convergence analysis demonstrates the absence of Runge's phenomenon as the expansion order increases. We validate the proposed approach against an analytical solution using the Bessel-Fourier harmonics supplemented by Graf's addition theorem, as well as perturbation-based, finite-element-based, finite-integration-based numerical solutions. This work can be found in [24] and is an extension of the prior work [25].

In Chapter 4, we present an improved SEM to model eccentric two-layer waveguides filled with anisotropic non-reciprocal media. The formulation employs transformation optics theory to map the original eccentric waveguide into an equivalent concentric problem. The TO is extended to include non-symmetric and non-Hermitian media tensors. We validated our approach against finite element solutions, and the results show that our technique requires far fewer degrees of freedom than FEM. This work can also be found in [26].

In Chapter 5, we model excitation sources for logging sensors employed in geophysical exploration using the fields obtained by SEM. We utilize horizontal coil antennas for both the sensor's transmitter and receiver. We examine several validating scenarios, covering the absence and presence of oil invasion zones. Results are consistent with the mode-matching technique for the amplitude ratio and phase difference of the sensor voltage.

Additionally, Appendix A shows an analytical solution for the eccentric coaxial waveguide in terms of a Bessel–Fourier series, where the eigenvalues are obtained through the zeros of a characteristic determinant equation derived from the application of Graf's addition theorem (GAT) for cylindrical functions. This approach can also be found in [24]. In Appendix B, we model uniaxially anisotropic circular waveguides loaded with eccentric rods via SEM. This work can also be found in [27].

Finally, Chapter 6 summarizes the essential results of this Ph.D. thesis and proposed activities for future works.

# 2

# A Novel Spectral Element Method for the Analysis of Cylindrical Waveguides Filled with Complex Media

## 2.1
## Introduction

The study of electromagnetic propagation in cylindrically-layered anisotropic media is essential in several applications, such as analyzing and designing metamaterial-based microwave devices and optical structures. Some practical metamaterial realizations involve biaxially anisotropic media through dielectric nanoporous or thin film-layered structures. These metamaterials support modal fields associated with complex-order Bessel functions [5–8]. Analysis of such structures using analytical approaches is very challenging because of the need to calculate special functions where a strong coupling occurs between its order and argument (where both are complex-valued quantities). Another relevant application involving cylindrically-layered anisotropic media is modeling antennas and sensors for propagation in tunnels and boreholes for geophysical prospecting. In this case, the constitutive parameters may exhibit fully anisotropic characteristics and extreme conductivity values [13–16].

The spectral element method (SEM) is an advanced implementation of the FEM, which also involves partitioning the computational domain into discrete elements. However, in SEM, the internal nodes of each element are located at the zeros of specific orthogonal polynomial families. Typically, Lobatto, Chebyshev, and Laguerre polynomials are used in the SEM. This feature enables the avoidance of the Runge effect, which refers to large oscillations that occur in the polynomial approximation of regions near the edges of a surface in higher-order interpolations. As a result, SEM can achieve higher accuracy with fewer elements compared to a corresponding FEM implementation [28, 29].

SEM has been applied successfully to various electromagnetic problems, such as a dielectric-loaded cavity, photonic crystals, metasurfaces, lithography models, and for different truncation boundary conditions [17–21]. This chapter will investigate a novel 2D SEM in cylindrical coordinates for modeling cylindrically-layered waveguides composed of radial lossy and anisotropic layers in bounded and unbounded domains. The presented formulation is based on cylindrical coordinates because many simplifications (in terms of both the development effort and the resulting computational load) are possible when solving cylindrical geometries. The results indicate that the proposed SEM requires fewer elements and degrees of freedom (DoF) than the conventional Cartesian SEM.

## 2.2
## Variational Waveguide Problem

### 2.2.1
### Electromagnetic Fields in Cylindrical Coordinates

Assuming and omitting the time-harmonic dependence $e^{-i\omega t}$, Maxwell's equations in an anisotropic and source-free region are given by

$$\boldsymbol{\nabla} \times \mathbf{E} = i\omega\mu_0\bar{\bar{\mu}}_r \cdot \mathbf{H} \tag{2-1a}$$

$$\boldsymbol{\nabla} \times \mathbf{H} = -i\omega\epsilon_0\bar{\bar{\epsilon}}_c \cdot \mathbf{E} \tag{2-1b}$$

$$\boldsymbol{\nabla} \cdot (\bar{\bar{\epsilon}}_c \cdot \mathbf{E}) = 0 \tag{2-1c}$$

$$\boldsymbol{\nabla} \cdot (\bar{\bar{\mu}}_r \cdot \mathbf{H}) = 0, \tag{2-1d}$$

where $\mathbf{E}$ and $\mathbf{H}$ are the electric and the magnetic fields, respectively. As usual, $\mu_0$ and $\epsilon_0$ denote the vacuum permeability and permittivity, respectively. The medium is characterized by real-valued relative permeability tensor

$$\bar{\bar{\mu}}_r = \begin{bmatrix} \bar{\bar{\mu}}_s & \bar{0} \\ \bar{0}^T & \mu_{zz} \end{bmatrix} \text{ with } \bar{\bar{\mu}}_s = \begin{bmatrix} \mu_{\rho\rho} & \mu_{\rho\phi} \\ \mu_{\phi\rho} & \mu_{\phi\phi} \end{bmatrix} \tag{2-2}$$

and complex-valued relative permittivity tensor

$$\bar{\bar{\epsilon}}_c = \bar{\bar{\epsilon}}_r + \frac{i}{\omega\epsilon_0}\bar{\bar{\sigma}}, \tag{2-3}$$

where

$$\bar{\bar{\epsilon}}_r = \begin{bmatrix} \bar{\bar{\epsilon}}_{rs} & \bar{0} \\ \bar{0}^T & \epsilon_{rz} \end{bmatrix} \text{ with } \bar{\bar{\epsilon}}_{rs} = \begin{bmatrix} \epsilon_{\rho\rho} & \epsilon_{\rho\phi} \\ \epsilon_{\phi\rho} & \epsilon_{\phi\phi} \end{bmatrix}, \tag{2-4}$$

and

$$\bar{\bar{\sigma}} = \begin{bmatrix} \bar{\bar{\sigma}}_s & \bar{0} \\ \bar{0}^T & \sigma_{zz} \end{bmatrix} \text{ with } \bar{\bar{\sigma}}_s = \begin{bmatrix} \sigma_{\rho\rho} & \sigma_{\rho\phi} \\ \sigma_{\phi\rho} & \sigma_{\phi\phi} \end{bmatrix} \tag{2-5}$$

are the real-valued relative permittivity and conductive tensors, respectively. Such tensors are decomposed into axial (or longitudinal, i.e., along with $z$) and transverse components, with subscripts $z$ and $s$, respectively. The transverse components are expressed in polar coordinates $(\rho, \phi)$.

We are interested in solving Maxwell's equations in guided structure invariant regarding the axial direction $z$. Consider a closed waveguide defined by a right cylinder with cross-section $\Omega \subset \mathbb{R}^2$. Assuming the propagation along the positive $z-$direction, the electric and magnetic fields can be expressed as

$$\mathbf{E}(\rho, \phi, z) = \mathbf{e}(\rho, \phi)e^{ik_z z} \tag{2-6a}$$

$$\mathbf{H}(\rho, \phi, z) = \mathbf{h}(\rho, \phi)e^{ik_z z}, \tag{2-6b}$$

where the vector-fields $\mathbf{e}$ and $\mathbf{h}$ are defined in the $\Omega$ domain and $k_z \in \mathbb{C}$ is the axial wavenumber. Our goal is to rewrite Maxwell's equations (2-1) in a system of equations for the electric field components $(\mathbf{e}_s, e_z)$ defined in the $\Omega$, i.e., to reduce a 3D problem to a 2D waveguide problem.

The fields and del-operator can be decomposed into axial and transversal components as

$$\mathbf{e}(\rho, \phi) = \mathbf{e}_s(\rho, \phi) + \hat{z}e_z(\rho, \phi), \tag{2-7a}$$

$$\mathbf{h}(\rho, \phi) = \mathbf{h}_s(\rho, \phi) + \hat{z}h_z(\rho, \phi), \tag{2-7b}$$

$$\boldsymbol{\nabla} = \boldsymbol{\nabla}_s + \hat{z}ik_z. \tag{2-7c}$$

Returning to Maxwell's equations, using (2-2)-(2-7) into (2-1a), we obtain

$$\boldsymbol{\nabla}_s \times \mathbf{e}_s = i\omega\mu_0\mu_{zz}\hat{z}h_z, \tag{2-8a}$$

$$\boldsymbol{\nabla}_s \times \hat{z}e_z + ik_z\hat{z} \times \mathbf{e}_s = i\omega\mu_0(\bar{\bar{\mu}}_s \cdot \mathbf{h}_s). \tag{2-8b}$$

Using (2-2)-(2-7) into (2-1b), we get

$$\boldsymbol{\nabla}_s \times \mathbf{h}_s = -i\omega\epsilon_0\epsilon_{cz}\hat{z}e_z, \tag{2-8c}$$

$$\boldsymbol{\nabla}_s \times \hat{z}h_z + ik_z\hat{z} \times \mathbf{h}_s = -i\omega\epsilon_0(\bar{\bar{\epsilon}}_{cs} \cdot \mathbf{e}_s). \tag{2-8d}$$

Also, using (2-2)-(2-7) into (2-1c), we get

$$\boldsymbol{\nabla}_s \cdot (\bar{\bar{\epsilon}}_{cs} \cdot \mathbf{e}_s) + ik_z\epsilon_{cz}e_z = 0. \tag{2-8e}$$

Now, applying $\boldsymbol{\nabla}_s \times \mu_{zz}^{-1}$ to (2-8a), and using (2-8d) and (2-8b), we obtain

$$\boldsymbol{\nabla}_s \times \mu_{zz}^{-1}\boldsymbol{\nabla}_s \times \mathbf{e}_s + ik_z\bar{\bar{\mu}}_s^{-1} \cdot \boldsymbol{\nabla}_s e_z - k_0^2\bar{\bar{\epsilon}}_{cs} \cdot \mathbf{e}_s + k_z^2\bar{\bar{\mu}}_s^{-1} \cdot \mathbf{e}_s = 0. \tag{2-9a}$$

Similarly, applying $\boldsymbol{\nabla}_s \times \bar{\bar{\mu}}_s^{-1}$ to (2-8b) and using (2-8c), we get

$$\boldsymbol{\nabla}_s \cdot (\bar{\bar{\mu}}_s^{-1} \cdot \boldsymbol{\nabla}_s e_z) - ik_z\boldsymbol{\nabla}_s \cdot (\bar{\bar{\mu}}_s^{-1} \cdot \mathbf{e}_s) + k_0^2\epsilon_{cz}e_z = 0, \tag{2-9b}$$

where $k_0^2 = \omega^2\epsilon_0\mu_0$ is the free-space wavenumber.

In possession of the equations (2-9a), (2-9b) and (2-8e), we can define a system of vector equations in terms of the electric field components, such as

$$\begin{cases} \boldsymbol{\nabla}_s \times \mu_{zz}^{-1} \boldsymbol{\nabla}_s \times \mathbf{e}_s + i k_z \bar{\bar{\mu}}_s^{-1} \cdot \boldsymbol{\nabla}_s e_z - k_0^2 \bar{\bar{\epsilon}}_{cs} \cdot \mathbf{e}_s + k_z^2 \bar{\bar{\mu}}_s^{-1} \cdot \mathbf{e}_s &= 0 \\ \boldsymbol{\nabla}_s \cdot (\bar{\bar{\mu}}_s^{-1} \cdot \boldsymbol{\nabla}_s e_z) - i k_z \boldsymbol{\nabla}_s \cdot (\bar{\bar{\mu}}_s^{-1} \cdot \mathbf{e}_s) + k_0^2 \epsilon_{cz} e_z &= 0 \\ \boldsymbol{\nabla}_s \cdot (\bar{\bar{\epsilon}}_{cs} \cdot \mathbf{e}_s) + i k_z \epsilon_{cz} e_z &= 0 \end{cases} \quad \text{in } \Omega. \tag{2-10}$$

For the physical problem of interest, perfect electric conductor (PEC) boundary conditions are imposed according to

$$\begin{aligned} \mathbf{e}_s \times \hat{n} &= \mathbf{0} \\ e_z &= 0 \end{aligned} \quad \text{on } \partial\Omega, \tag{2-11}$$

where $\hat{n}$ is the unit outward normal on $\partial\Omega$. Since no sources are given, (2-10) represents a 2D waveguide problem. The goal is to find all possible pairs of eigenvalues and eigenfunctions $(k_z, (\mathbf{e}_s, e_z))$ that solve (2-10) and satisfy the boundary conditions in (2-11).

## 2.2.2
**Variational Formulation**

On the grounds of the theory in [30–35], before we construct our variational formulation, let us develop a functional setting suitable for the problem at hand. The following Hilbert spaces are introduced:

$$H_0^1(\Omega) = \{w \in L_2(\Omega) \mid \boldsymbol{\nabla}_s w \in L_2(\Omega)^2, \; w|_{\partial\Omega} = 0\} \tag{2-12a}$$

$$\mathbf{H}(\text{div}, \Omega) = \{\mathbf{u} \in L_2(\Omega)^2 \mid \boldsymbol{\nabla}_s \cdot \mathbf{u} \in L_2(\Omega)\} \tag{2-12b}$$

$$\mathbf{H}_0(\text{curl}, \Omega) = \{\mathbf{u} \in L_2(\Omega)^2 \mid \boldsymbol{\nabla}_s \times \mathbf{u} \in L_2(\Omega), \; \hat{n} \times \mathbf{u}|_{\partial\Omega} = 0\}, \tag{2-12c}$$

where $L_2(\Omega)$ is a square-integrable function space. The $L_2(\Omega)$ inner product extends trivially to vector functions. Supposing that $\mathbf{u} = (u_1, u_2)^T \in L_2(\Omega)^2$ and $\mathbf{v} = (v_1, v_2)^T \in L_2(\Omega)^2$, we can then write the $L_2(\Omega)^2$ inner product as

$$\langle \mathbf{u}, \mathbf{v} \rangle = \int_\Omega \sum_{j=1}^2 u_j v_j^* \; dS. \tag{2-13}$$

In the above, notice we choose to conjugate the second term as a convention.

Given $\mathbf{u} \in \mathbf{H}_0(\text{curl}, \Omega)$, using results of Divergence Theorem (Corollary 3.20 in [35]), we get

$$\begin{aligned} \langle \boldsymbol{\nabla}_s \times \mathbf{u}, \mathbf{w} \rangle &= \int_\Omega \boldsymbol{\nabla}_s \times \mathbf{u} \cdot \mathbf{w} \; dS \\ &= \int_\Omega \mathbf{u} \cdot \boldsymbol{\nabla}_s \times \mathbf{w} \; dS + \int_{\partial\Omega} \hat{n} \times \mathbf{u} \cdot \mathbf{w} \; dl \\ &= \int_\Omega \mathbf{u} \cdot \boldsymbol{\nabla}_s \times \mathbf{w} \; dS = \langle \mathbf{u}, \boldsymbol{\nabla}_s \times \mathbf{w} \rangle, \end{aligned} \tag{2-14a}$$

for all $\mathbf{w} \in \mathbf{H}_0(\text{curl}, \Omega)$, where $\hat{n}$ is the unit outward normal on $\partial\Omega$. Similarly,

given $\mathbf{u} \in \mathbf{H}(\mathrm{div}, \Omega)$, we readily obtain

$$
\begin{aligned}
\langle \boldsymbol{\nabla}_s \cdot \mathbf{u}, w \rangle &= \int_\Omega \boldsymbol{\nabla}_s \cdot \mathbf{u} \; w \, dS \\
&= -\int_\Omega \mathbf{u} \cdot \boldsymbol{\nabla}_s w \; dS + \int_{\partial\Omega} \hat{n} \cdot \mathbf{u} \; w \, dl \\
&= -\int_\Omega \mathbf{u} \cdot \boldsymbol{\nabla}_s w \; dS = -\langle \mathbf{u}, \boldsymbol{\nabla}_s w \rangle,
\end{aligned} \tag{2-14b}
$$

for all $w \in H_0^1(\Omega)$.

To construct a weak formulation of waveguide problem (2-10), we first take the inner product of the vectorial Helmholtz's equation with appropriate weight functions $(\mathbf{w}, w) \in \mathbf{H}_0(\mathrm{curl}, \Omega) \times H_0^1(\Omega)$. Using (2-14a), the inner product between (2-9a) and the vector weight function $\mathbf{w} \in \mathbf{H}_0(\mathrm{curl}, \Omega)$ is given by

$$
\langle \mu_{zz}^{-1} \boldsymbol{\nabla}_s \times \mathbf{e}_s, \boldsymbol{\nabla}_s \times \mathbf{w} \rangle + ik_z \langle \bar{\bar{\mu}}_s^{-1} \boldsymbol{\nabla}_s e_z, \mathbf{w} \rangle - k_0^2 \langle \bar{\bar{\epsilon}}_{cs} \mathbf{e}_s, \mathbf{w} \rangle + k_z^2 \langle \bar{\bar{\mu}}_s^{-1} \mathbf{e}_s, \mathbf{w} \rangle = 0. \tag{2-15a}
$$

Using (2-14b), the inner product between (2-9b) and (2-8e) with scalar weight function $w \in H_0^1(\Omega)$, respectively, results in

$$
\langle \bar{\bar{\mu}}_s^{-1} \boldsymbol{\nabla}_s e_z, \boldsymbol{\nabla}_s w \rangle - ik_z \langle \bar{\bar{\mu}}_s^{-1} \mathbf{e}_s, \boldsymbol{\nabla}_s w \rangle - k_0^2 \langle \epsilon_{cz} e_z, w \rangle = 0 \tag{2-15b}
$$

and

$$
\langle \bar{\bar{\epsilon}}_{cs} \mathbf{e}_s, \boldsymbol{\nabla}_s w \rangle - ik_z \langle \epsilon_{cz} e_z, w \rangle = 0. \tag{2-15c}
$$

The set of equations (2-15) forms an overdetermined system for the electric field components. As already demonstrated in the works in [18, 32], this system is linearly dependent. For $k_z \neq 0$, consider $\mathbf{w} = \boldsymbol{\nabla}_s w$, we can rewrite (2-15a) as

$$
ik_z \langle \bar{\bar{\mu}}_s^{-1} \boldsymbol{\nabla}_s e_z, \boldsymbol{\nabla}_s w \rangle - k_0^2 \langle \bar{\bar{\epsilon}}_{cs} \mathbf{e}_s, \boldsymbol{\nabla}_s w \rangle + k_z^2 \langle \bar{\bar{\mu}}_s^{-1} \mathbf{e}_s, \boldsymbol{\nabla}_s w \rangle = 0, \tag{2-16}
$$

because $\boldsymbol{\nabla}_s \times \boldsymbol{\nabla}_s w = \mathbf{0}$. Multiplying (2-15b) by $ik_z$, we get

$$
ik_z \langle \bar{\bar{\mu}}_s^{-1} \boldsymbol{\nabla}_s e_z, \boldsymbol{\nabla}_s w \rangle + k_z^2 \langle \bar{\bar{\mu}}_s^{-1} \mathbf{e}_s, \boldsymbol{\nabla}_s w \rangle - ik_z k_0^2 \langle \epsilon_{cz} e_z, w \rangle = 0. \tag{2-17}
$$

Subtracting (2-16) from (2-17) and dividing the result by $k_0^2$, we return (2-15c). This completes the proof of linear dependence of the equations in (2-15).

Based on the above result, we do not need to consider all three equations of (2-15); just two are enough. Previous works show that considering equations (2-15a) and (2-15c), i.e., including Gauss' law in the system of equations, the variational eigenvalue problem does not present spurious modes [18, 19, 32, 33].

Now, we can propose a new eigenvalue problem. For a given $k_0 \neq 0$, we want to find all pairs $(k_z, (\mathbf{e}_s, e_z)) \in (\mathbb{C}, \mathbf{H}_0(\mathrm{curl}, \Omega) \times H_0^1(\Omega))$ that satisfies

$$\begin{cases} \langle \mu_{zz}^{-1} \boldsymbol{\nabla}_s \times \mathbf{e}_s, \boldsymbol{\nabla}_s \times \mathbf{w} \rangle + ik_z \langle \bar{\bar{\mu}}_s^{-1} \boldsymbol{\nabla}_s e_z, \mathbf{w} \rangle - k_0^2 \langle \bar{\bar{\epsilon}}_{cs} \mathbf{e}_s, \mathbf{w} \rangle + k_z^2 \langle \bar{\bar{\mu}}_s^{-1} \mathbf{e}_s, \mathbf{w} \rangle = 0 \\ \langle \bar{\bar{\epsilon}}_{cs} \mathbf{e}_s, \boldsymbol{\nabla}_s w \rangle - ik_z \langle \epsilon_{cz} e_z, w \rangle = 0 \end{cases}$$

(2-18)

for all $(\mathbf{w}, w) \in \mathbf{H}_0(\text{curl}, \Omega) \times H_0^1(\Omega)$. The problem above is known as the *weak formulation* for the waveguide problem in (2-10).

## 2.3
## Spectral Element Method

### 2.3.1
### Discretized Domain

Our computational domain was first discretized into a series of elements covering the waveguide cross-section. Supposing $\Omega$ is a circular or coaxial domain, we can discretize $\Omega$ by a finite set of curvilinear elements $\Omega^e$ for $e = 1, 2, ..., M$ such that $\Omega = \bigcup_{e=1}^{M} \Omega^e$, where $\text{area}(\Omega^e \cap \Omega^{e'}) = 0$ for all $e \neq e'$. Each element is described by the polar coordinates $(\rho, \phi)$, which is then mapped into the reference coordinates $(\xi, \eta)$ via

$$\rho^e(\xi) = 0.5(\rho_0^e + \rho_1^e) + 0.5(\rho_1^e - \rho_0^e)\xi, \quad \forall \xi \in [-1, 1], \tag{2-19a}$$

$$\phi^e(\eta) = 0.5(\phi_0^e + \phi_1^e) + 0.5(\phi_1^e - \phi_0^e)\eta, \quad \forall \eta \in [-1, 1]. \tag{2-19b}$$

The mapping above allows us to relate the transverse points of the reference element (covering $-1 \leq \xi \leq 1$ and $-1 \leq \eta \leq 1$) with the physical element. The mapping between reference and physical elements is illustrated in Fig. 2.1, and an example of discretization of a coaxial domain considering $M = 12$ is presented in Fig. 2.2. The associated Jacobian matrix of this coordinate transformation is given by

$$\bar{\bar{J}}_s^e = \begin{bmatrix} \frac{\partial \rho^e}{\partial \xi} & \frac{\partial \rho^e}{\partial \eta} \\ \rho \frac{\partial \phi^e}{\partial \xi} & \rho \frac{\partial \phi^e}{\partial \eta} \end{bmatrix} = \begin{bmatrix} 0.5(\rho_1^e - \rho_0^e) & 0 \\ 0 & 0.5\rho(\phi_1^e - \phi_0^e) \end{bmatrix}. \tag{2-20}$$

The associated Jacobian shows that $\rho$- and $\eta$-directions are decoupled, and consequently, directions $\phi$ and $\xi$ are also decoupled.

### 2.3.2
### Lagrange Interpolation

To represent the field components in reference coordinates, we propose using Lagrange interpolation basis functions associated with Gauss-Lobatto-Legendre (GLL) sampling points. For example, with the sampling points $\xi_k \in [-1, 1]$, for all $j = 1, 2, ..., N + 1$, along the coordinate $\xi$, the $N$th-order

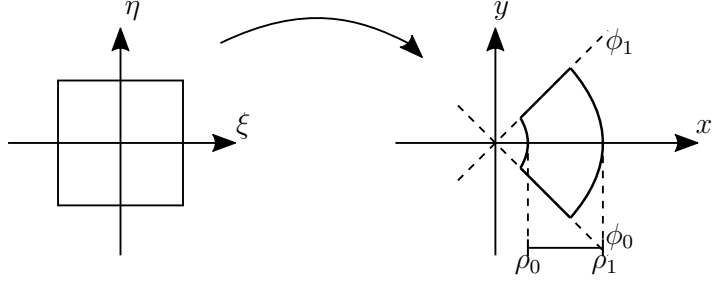

Figure 2.1: Mapping between transverse points of the reference element and the physical element.

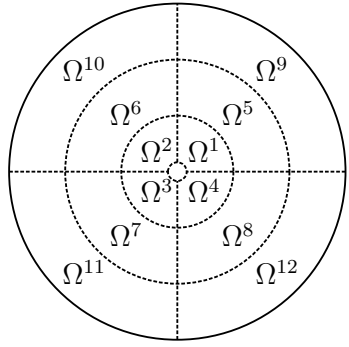

Figure 2.2: Discretization of the computational domain considering twelve curvilinear elements.

one-dimensional (1D) GLL basis function is given by

$$\phi_j^N(\xi) = \prod_{k \neq j}^{N+1} \frac{\xi - \xi_k}{\xi_j - \xi_k} \quad \forall \xi \in [-1, 1], \tag{2-21}$$

where the $N + 1$ interpolation nodes $\xi_k$ are the roots of the $(N + 1)$-degree completed Lobatto polynomial $(1 - \xi^2)Lo_{N-1}(\xi) = 0$. Cursory inspection reveals that these polynomials satisfy the cardinal interpolation property

$$\phi_j^N(\xi_k) = \begin{cases} 1 & \text{if } j = k \\ 0 & \text{if } j \neq k \end{cases}. \tag{2-22}$$

Stated differently, $\phi_j^N(\xi_k)$ is the identity matrix represented by Kronecker's delta, $\phi_j^N(\xi_k) = \delta_{jk}$ [28]. Accordingly, the two-dimensional (2D) GLL basis function is obtained via tensor-product expansions

$$\hat{\psi}_{z,ij}(\xi, \eta) = \phi_i^N(\xi)\phi_j^N(\eta) \quad \forall (\xi, \eta) \in [-1, 1]^2, \tag{2-23}$$

and the mixed-order curl-conforming vector edge-based basis functions are written as [18, 19]

$$\hat{\Phi}_{\xi,ij}(\xi,\eta) = \phi_i^{N-1}(\xi)\phi_j^N(\eta)$$

$$\hat{\Phi}_{\eta,ij}(\xi,\eta) = \phi_i^N(\xi)\phi_j^{N-1}(\eta)$$

(2-24)

for all $(\xi,\eta) \in [-1,1]^2$.

On the other hand, the corresponding basis functions $\psi_{z,ij}^e$ and $\Phi_{l,ij}^e$ for $l \in \{\xi,\eta\}$, and their gradient and curl in the physical element can be obtained by the covariant mappings [35]

$$\begin{cases} \psi_{z,ij}^e(\rho,\phi) = \hat{\psi}_{z,ij}(\xi,\eta) \\[2mm] \Phi_{l,ij}^e(\rho,\phi) = \bar{\bar{J}}_s^{e-1}\hat{\Phi}_{l,ij}(\xi,\eta) \\[2mm] \boldsymbol{\nabla}_s\psi_{z,ij}^e(\rho,\phi) = \bar{\bar{J}}_s^{e-1}\hat{\boldsymbol{\nabla}}_s\hat{\psi}_{z,ij}(\xi,\eta) \\[2mm] \boldsymbol{\nabla}_s \times \Phi_{l,ij}^e(\rho,\phi) = \frac{1}{|\bar{\bar{J}}_s^e|}\hat{\boldsymbol{\nabla}}_s \times \hat{\Phi}_{l,ij}(\xi,\eta) \end{cases},$$

(2-25)

where $\boldsymbol{\nabla}_s = \left(\frac{\partial}{\partial\rho}, \frac{1}{\rho}\frac{\partial}{\partial\phi}\right)^T$, $\hat{\boldsymbol{\nabla}}_s = \left(\frac{\partial}{\partial\xi}, \frac{\partial}{\partial\eta}\right)^T$, the Jacobian matrix $\bar{\bar{J}}_s^e$ is defined by (2-20), and $|\bar{\bar{J}}_s^e|$ is its respective determinant.

The choices for the interpolation basis functions result in $2N(N+1) + (N+1)^2$ degrees of freedom (DoF) on each element, where $(N+1)^2$ for $\hat{\psi}_{z,ij}$, $N(N+1)$ for $\hat{\Phi}_{\xi,ij}$, plus $N(N+1)$ for $\hat{\Phi}_{\eta,ij}$. Fig. 2.3 shows the distribution of these basis functions in the reference element considering $N=5$.

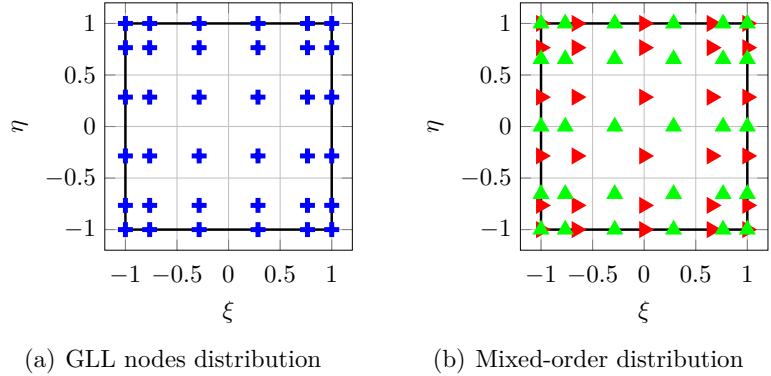

(a) GLL nodes distribution      (b) Mixed-order distribution

Figure 2.3: Lagrange interpolation basis functions in the reference element $[-1,1]^2$ with the order $N=5$ generated in MATLAB.

### 2.3.3
### Eigenvalue Problem

To construct the discrete form of the problem (2-18), we will approximate the electric field components using GLL basis functions. The transversal component in the reference domain, $\hat{\mathbf{e}}_s = (\hat{e}_\xi, \hat{e}_\eta)$, can be expanded as

$$\hat{e}_\xi(\xi, \eta) = \sum_{j=1}^{N(N+1)} \hat{e}_{\xi,j} \hat{\Phi}_{\xi,j}(\xi, \eta), \tag{2-26a}$$

$$\hat{e}_\eta(\xi, \eta) = \sum_{j=1}^{N(N+1)} \hat{e}_{\eta,j} \hat{\Phi}_{\eta,j}(\xi, \eta), \tag{2-26b}$$

and the axial field $\hat{e}_z$ is approximated as

$$\hat{e}_z(\xi, \eta) = \sum_{j=1}^{(N+1)^2} \hat{e}_{z,j} \hat{\psi}_{z,j}(\xi, \eta). \tag{2-26c}$$

For weight functions, we will use Galerkin's method, where the vector weight function $\mathbf{w} = (\Phi_{\xi,i}, \Phi_{\eta,i})^T$ for $i = 1, ...N(N+1)$, and scalar weight function $w = \psi_{z,i}$ for $i = 1, ..., (N+1)^2$. Using expansions (2-26) and covariant mappings in (2-25), we can define elemental matrices that represent the inner products in the problem (2-18). These matrices can be written in reference coordinates [18] as

$$(\mathbf{S}_{ss})_{ij}^e = \int_{-1}^{1} \int_{-1}^{1} \frac{1}{|\bar{\bar{J}}_s^e|} (\mu_{zz}^{-1} \boldsymbol{\nabla}_s \times \hat{\Phi}_{l,j})^T \cdot \boldsymbol{\nabla}_s \times \hat{\Phi}_{l,i} d\xi d\eta \tag{2-27a}$$

$$(\mathbf{M}_{ss1})_{ij}^e = \int_{-1}^{1} \int_{-1}^{1} (\bar{\bar{\epsilon}}_{cs} \bar{\bar{J}}_s^{e-1} \hat{\Phi}_{l,j})^T \cdot \bar{\bar{J}}_s^{e-1} \hat{\Phi}_{l,i} |\bar{\bar{J}}_s^e| d\xi d\eta \tag{2-27b}$$

$$(\mathbf{M}_{ss2})_{ij}^e = \int_{-1}^{1} \int_{-1}^{1} (\bar{\bar{\mu}}_s^{-1} \bar{\bar{J}}_s^{e-1} \hat{\Phi}_{l,j})^T \cdot \bar{\bar{J}}_s^{e-1} \hat{\Phi}_{l,i} |\bar{\bar{J}}_s^e| d\xi d\eta \tag{2-27c}$$

$$(\mathbf{M}_{zz})_{ij}^e = \int_{-1}^{1} \int_{-1}^{1} (\epsilon_{cz} \hat{\psi}_{z,j}) \cdot \hat{\psi}_{z,i} |\bar{\bar{J}}_s^e| d\xi d\eta \tag{2-27d}$$

$$(\mathbf{K}_{sz})_{ij}^e = \int_{-1}^{1} \int_{-1}^{1} (\bar{\bar{\mu}}_s^{-1} \bar{\bar{J}}_s^{e-1} \boldsymbol{\nabla}_s \hat{\psi}_{z,j})^T \cdot \hat{\Phi}_{l,i} |\bar{\bar{J}}_s^e| d\xi d\eta \tag{2-27e}$$

$$(\mathbf{K}_{zs})_{ij}^e = \int_{-1}^{1} \int_{-1}^{1} (\bar{\bar{\epsilon}}_{cs} \bar{\bar{J}}_s^{e-1} \hat{\Phi}_{l,j})^T \cdot \bar{\bar{J}}_s^{e-1} \boldsymbol{\nabla}_s \hat{\psi}_{z,i} |\bar{\bar{J}}_s^e| d\xi d\eta \tag{2-27f}$$

where the superscript $e$ is used to indicate the parameters of the reference element, and $l \in \{\xi, \eta\}$. The global matrices $\mathbf{S}_{ss}$, $\mathbf{M}_{ss1}$, $\mathbf{M}_{ss2}$, $\mathbf{M}_{zz}$, $\mathbf{K}_{sz}$, and $\mathbf{K}_{zs}$ are assembled by elemental matrices and are subject to the boundary conditions in (2-11). Now, we can define the system of matrix equations

$$\mathbf{S}_{ss} \hat{\mathbf{e}}_{s,mn} + i k_z \mathbf{K}_{sz} \hat{\mathbf{e}}_{z,p} - k_0^2 \mathbf{M}_{ss1} \hat{\mathbf{e}}_{s,mn} + k_z^2 \mathbf{M}_{ss2} \hat{\mathbf{e}}_{s,mn} = 0, \tag{2-28a}$$

$$\mathbf{K}_{zs} \hat{\mathbf{e}}_{s,mn} - i k_z \mathbf{M}_{zz} \hat{\mathbf{e}}_{z,p} = 0, \tag{2-28b}$$

where $\hat{\mathbf{e}}_{s,mn} = (\hat{e}_{\xi,1}, ..., \hat{e}_{\xi,m}, \hat{e}_{\eta,1}, ..., \hat{e}_{\eta,n})^T$ and $\hat{\mathbf{e}}_{z,p} = (\hat{e}_{z,1}, ..., \hat{e}_{z,p})^T$ are vectors that contain the unknowns of the electric field for transversal and axial components in reference coordinates, respectively. Finally, isolating $\hat{\mathbf{e}}_{zp}$ in (2-28b)

$$\hat{\mathbf{e}}_{z,p} = \frac{1}{i k_z} \mathbf{M}_{zz}^{-1} \mathbf{K}_{zs} \hat{\mathbf{e}}_{s,mn}, \tag{2-29}$$

and substituting the result in (2-28a), we obtain the discrete eigenvalue problem associated with (2-18) given by

$$\left( -\mathbf{S}_{ss} - \mathbf{K}_{sz}\mathbf{M}_{zz}^{-1}\mathbf{K}_{zs} + k_0^2\mathbf{M}_{ss1} \right) \hat{\mathbf{e}}_{s,mn} = k_z^2\mathbf{M}_{ss2}\hat{\mathbf{e}}_{s,mn}, \qquad (2\text{-}30)$$

where $k_z^2$ and $\hat{\mathbf{e}}_{s,mn}$ are are the eigenvalues and eigenvectors, respectively. The eigenvalues are associated with the propagation modes in the axial direction, while the eigenvectors are associated with the transverse components of the electric field in the reference coordinates.

The fields in cylindrical coordinates can be obtained via the inverse transformation of the Jacobian matrix. Therefore, the electric field components in each element can be written as

$$\begin{bmatrix} e_{\rho,i} \\ e_{\phi,i} \\ e_{z,i} \end{bmatrix} = \begin{bmatrix} \frac{\partial \rho^e}{\partial \xi} & 0 & 0 \\ 0 & \rho\frac{\partial \phi^e}{\partial \eta} & 0 \\ 0 & 0 & 1 \end{bmatrix}^{-1} \begin{bmatrix} \hat{e}_{\xi,i} \\ \hat{e}_{\eta,i} \\ \hat{e}_{z,i} \end{bmatrix}, \qquad (2\text{-}31)$$

where this transformation is performed for each mesh point, with the subscript $i$ representing any point of the range.

It is important to mention that the present SEM formulation has two important aspects that contrast with the former approaches. First, the mapping in (2-19) simplifies waveguides that present cylindrical-conforming boundaries. Second, the associated Jacobian in (2-31) shows that $\rho$- and $\eta$- field components are decoupled. Components $\phi$ and $\xi$ are also decoupled. Such characteristics that simplify the computational algorithm are absent in conventional Cartesian-coordinate-based SEM formulations [18, 19].

### 2.3.4
### Lobatto Integration Quadrature

The integrals on the right-hand side of (2-27) are composed of higher-order polynomials. To perform the integration, we adopt the Lobatto integration quadrature, also called the Gauss-Lobatto quadrature. The Lobatto quadrature with $(N+1)$ points on each axis is defined such that, for a smooth function $f(\xi, \eta)$, the following equality holds

$$\int_{-1}^{1}\int_{-1}^{1} f(\xi,\eta)d\xi d\eta = \sum_{p=1}^{N+1}\sum_{q=1}^{N+1} f(z_p, z_q)w_p w_q, \qquad (2\text{-}32)$$

where $N \geq 1$ denotes the chosen order of quadrature, $(z_p, z_q)$ are the GLL nodes, and $(w_p, w_q)$ are the weights corresponding to the GLL nodes. The

weights are computed as

$$w_1 = w_{N+1} = \frac{2}{N(N+1)}, \tag{2-33a}$$

$$w_r = \frac{2}{N(N+1)} \frac{1}{L_N^2(z_r)}, \tag{2-33b}$$

for $r = 2, ..., N$, where $r \in \{p, q\}$ and $L_N$ denotes a Legendre polynomial [28, 29].

## 2.4
## Coordinate Transformations of Anisotropic Tensors

In our formulation, we assume the relative permittivity ($\bar{\bar{\epsilon}}_c$) and permeability ($\bar{\bar{\mu}}_r$) tensors are expressed in cylindrical coordinates. In many application problems, however, writing these parameters in Cartesian coordinates is more expedient. In what follows, distinctions between tensors in Cartesian and cylindrical formats will be made via the subscripts rec and cyl, respectively. Accordingly, the complex-valued tensors

$$\bar{\bar{p}}_{\text{rec}} = \begin{bmatrix} p_{xx} & p_{xy} & 0 \\ p_{yx} & p_{yy} & 0 \\ 0 & 0 & p_{zz} \end{bmatrix}, \ \bar{\bar{p}}_{\text{cyl}} = \begin{bmatrix} p_{\rho\rho} & p_{\rho\phi} & 0 \\ p_{\phi\rho} & p_{\phi\phi} & 0 \\ 0 & 0 & p_{zz} \end{bmatrix}, \tag{2-34}$$

with $\bar{\bar{p}} \in \{\bar{\bar{\epsilon}}_c, \bar{\bar{\mu}}_r\}$, can be transformed via [10, 36]

$$\bar{\bar{p}}_{\text{cyl}} = \bar{\bar{R}}(\phi) \cdot \bar{\bar{p}}_{\text{rec}} \cdot \bar{\bar{R}}^T(\phi), \tag{2-35}$$

where

$$\bar{\bar{R}}(\phi) = \begin{bmatrix} \cos\phi & \sin\phi & 0 \\ -\sin\phi & \cos\phi & 0 \\ 0 & 0 & 1 \end{bmatrix}. \tag{2-36}$$

The above formulas apply for transforming the relative permittivity and permeability, i.e., $\bar{\bar{p}} \in \{\bar{\bar{\epsilon}}_c, \bar{\bar{\mu}}_r\}$.

## 2.5
## Anisotropic Perfectly Matched Layer

The field propagation modeling in radially-unbounded waveguides such as rods and boreholes should employ appropriate absorbing boundary conditions to mimic an open space. We have explored the integration of a perfect matching layer (PML) in cylindrical coordinates [37–39] into the above presented SEM formulation. The PML layer extends over $r_{PML} \leq \rho \leq r_\infty$, as illustrated in Fig. 2.4. We first define a one-to-one mapping of the radial coordinate $\rho$ to a stretched $\tilde{\rho}$ as

$$\rho \longrightarrow \tilde{\rho} = \int_0^\rho s_\rho(\rho')d\rho'. \tag{2-37}$$

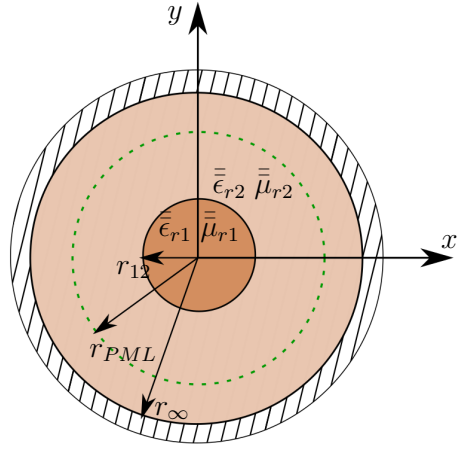

Figure 2.4: Cross-section of a circular waveguide bounded by a PML layer.

Using an integral in the above ensures that $\tilde{\rho}$ varies smoothly even if $s_\rho$ has discontinuities. Typically, the complex stretching variable $s_\rho$ employs a polynomial profile defined by

$$s_\rho(\rho) = \begin{cases} 1, & \text{for } \rho \in [0, r_{PML}) \\ 1 + i\alpha_{PML}\left(\frac{\rho - r_{PML}}{r_\infty - r_{PML}}\right)^q, & \text{for } \rho \in [r_{PML}, r_\infty] \end{cases} \tag{2-38}$$

which results in the mapping

$$\tilde{\rho} = \begin{cases} \rho, & \text{for } \rho \in [0, r_{PML}) \\ \rho + i\alpha_{PML}\frac{1}{1+q}\frac{(\rho - r_{PML})^{q+1}}{(r_\infty - r_{PML})^q}, & \text{for } \rho \in [r_{PML}, r_\infty] \end{cases} \tag{2-39}$$

where $\alpha_{PML} > 0$ [40].

This work will use the anisotropic PML formulation known as Maxwellian PML. In this formulation, the added degrees of freedom are entirely incorporated into the constitutive parameters, and the familiar form of Maxwell's equations is retained. The resultant electromagnetic fields inside the PML can be associated with anisotropic and lossy constitutive parameters. This allows for the interpretation of the fields inside the PML as physical fields and the straightforward application of the PML to different methods like FDTD, FEM, or SEM [41]. The PML constitutive tensors to arbitrary anisotropic media inside the absorbing layer can be written [38] as

$$\bar{\bar{\epsilon}}_{PML} = (\det \bar{\bar{S}})^{-1}[\bar{\bar{S}} \cdot \bar{\bar{\epsilon}}_r \cdot \bar{\bar{S}}], \tag{2-40a}$$

$$\bar{\bar{\mu}}_{PML} = (\det \bar{\bar{S}})^{-1}[\bar{\bar{S}} \cdot \bar{\bar{\mu}}_r \cdot \bar{\bar{S}}], \tag{2-40b}$$

where the stretching tensor $\bar{\bar{S}}$ is given by

$$\bar{\bar{S}} = \hat{\rho}\hat{\rho}\left(\frac{1}{s_\rho}\right) + \hat{\phi}\hat{\phi}\left(\frac{1}{s_\phi}\right) + \hat{z}\hat{z}\left(\frac{1}{s_z}\right), \tag{2-41}$$

where $s_\phi = \tilde{\rho}/\rho$, $s_z = 1$, and $s_\rho$ given by (2-38).

## 2.6
## Numerical Results

This section shows numerical examples for solving cylindrical waveguide problems using the above-presented SEM formulation. A variety of examples are considered here, including homogeneous and inhomogeneous filled waveguides, radially-bounded and radially-unbounded domains, lossless and lossy problems, and isotropic and anisotropic media. Our algorithms were written in the MATLAB environment [42], and our results are compared against analytical results from [43–45], the Cartesian-coordinate-based SEM from [18], and finite-element method (FEM) results from COMSOL Multiphysics [2].

### 2.6.1
### Homogeneous Coaxial Waveguide

As our first validation, we consider a homogeneous, isotropic, and lossless coaxial waveguide bounded by a PEC at the internal radius $r_0 = 0.1$ m and outer radius $r_1 = 1.0$ m. As a reference, we use well-known analytical eigenvalue solutions [43, Ch. 5] for this canonical problem. Considering the domain discretized into four elements and varying the expansion functions order $N$, Fig. 2.5 shows the relative error in the computed radial wavenumber $k_\rho$ for the first three modes used in our SEM algorithm. The relative error $\varepsilon$ is defined in terms of the numerical $f_n$ and analytical $f_a$ solutions by $\varepsilon = |f_n - f_a|/|f_a|$. We observe that the error decays exponentially with increasing polynomial order. The absence of numerical errors due to the Runge effect in higher-order interpolations is also evident, as predicted by theory [28, 29].

To validate the fields obtained via (2-31), we compare the numerical results of the present cylindrical-coordinate-based SEM (domain discretized into 4 elements and expansion functions of order 14) with analytical solutions for the first two TE modes in Fig. 2.6. The results show excellent agreement.

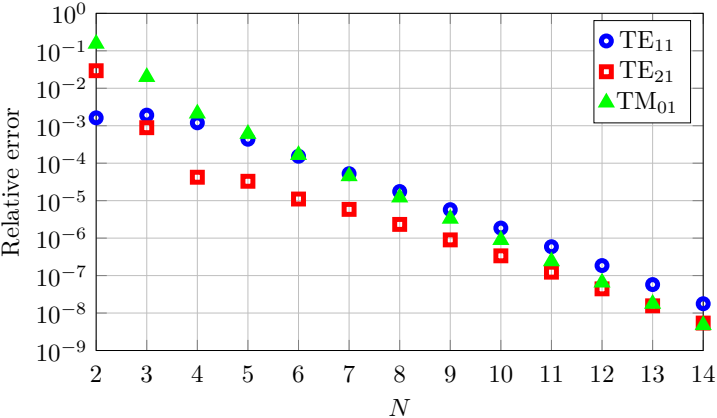

Figure 2.5: Relative error of the first three modes as a function of the basis function order $N$ used in our SEM algorithm.

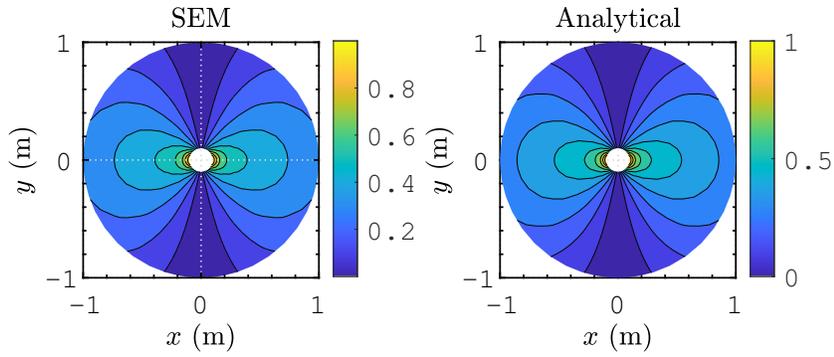

(a) Normalized $|e_\rho|$, TE$_{11}$

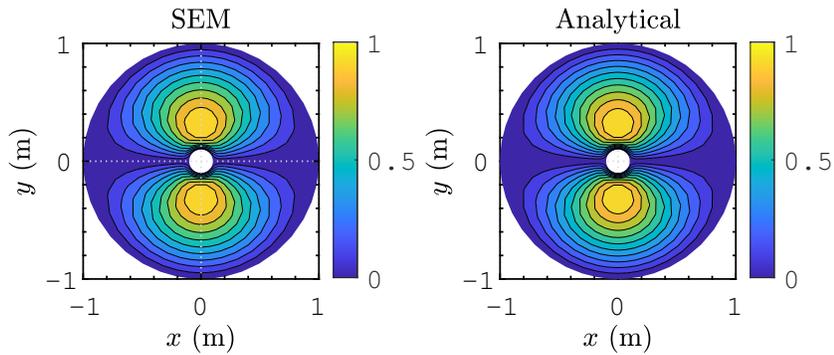

(b) Normalized $|e_\phi|$, TE$_{11}$

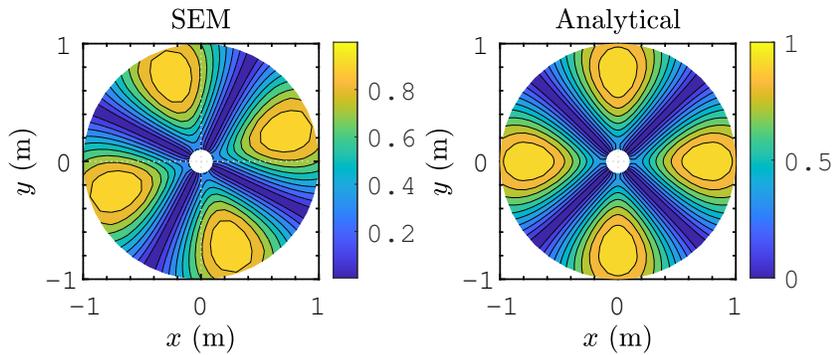

(c) Normalized $|e_\rho|$, TE$_{21}$

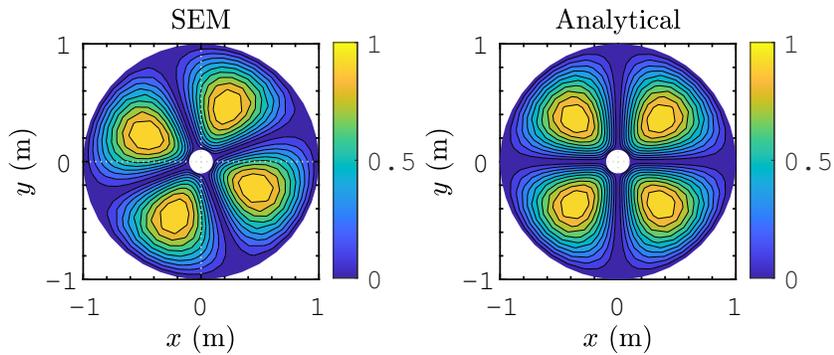

(d) Normalized $|e_\phi|$, TE$_{21}$

Figure 2.6: Electric field components in a homogeneous coaxial waveguide for TE modes were obtained via the present SEM (left) and via the analytical solution (right).

## 2.6.2
### Radially-Unbounded Circular Waveguides

Now, we investigate a specialization of the SEM for solving radially-unbounded problems. We use an anisotropic perfectly matched layer (PML) in cylindrical coordinates [37,38] to mimic the open boundary condition. As a first example, we consider a lossless isotropic dielectric rod defined by $\bar{\bar{\epsilon}}_c = 20\bar{\bar{I}}$ (where $\bar{\bar{I}}$ is the unity dyadic), with radius varying from $a = 0.4$ m to $a = 2.5$ m. The operating frequency is such that the vacuum wavenumber is given by $k_0 = 1$ m$^{-1}$. Other constitutive parameters are equal to the vacuum values. The computational domain is discretized into 12 elements (where four elements are for the PML) using expansion functions of order 5. This implies 484 DoF. Fig. 2.7 shows the results of the radial wavenumber $k_\rho$ in the innermost layer for the first eight modes obtained by the present SEM formulation. Excellent agreement is observed when comparing our solution with reference values from [45].

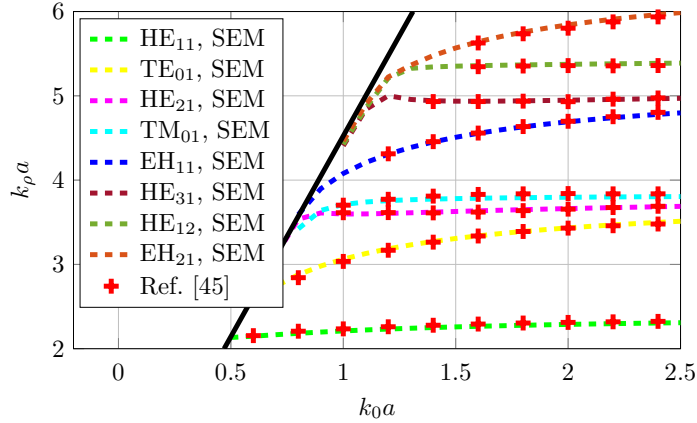

Figure 2.7: The first eight propagation modes of the dielectric rod waveguide are defined by $\bar{\bar{\epsilon}}_c = 20\bar{\bar{I}}$. The solid black line represents the cutoff for propagating modes.

As a second example, we consider a geometry similar to the latter problem. In the simulations, we consider a rod with radius $a$ filled with (a) the anisotropic permittivity defined in Cartesian rectangular coordinates by

$$\bar{\bar{\epsilon}}_{r_{\text{rec}}} = \begin{bmatrix} 3.2 & 0 & 0 \\ 0 & 2.048 & 0 \\ 0 & 0 & 2.56 \end{bmatrix}, \tag{2-42}$$

and (b) the isotropic medium $\bar{\bar{\epsilon}}_{r,\text{rec}} = 2.56\bar{\bar{I}}$. The anisotropic tensor is defined in Cartesian coordinates; it is necessary to perform the Cartesian-to-cylindrical coordinate transformation (described in Section 2.4 ) to make such representations conforming with (2-1)–(2-4) of our formulation. We considered several scenarios where the radius of the rod varied from $a = 1.8$ m to $a = 3.6$ m when the operating frequency was fixed in 50 MHz. These problems were considered before in [44] and [18] and are used here to validate our formulation for radially-open waveguides. The domain was discretized into 12 elements by using expansion functions of order 5 (with 484 DoF). Fig. 2.8 shows the normalized axial wavenumber of the fundamental mode as a function of the normalized rod diameter for both isotropic and anisotropic cases. Excellent agreement is observed versus the analytical results from [44] (the isotropic case) and from the Cartesian-coordinate-based SEM from [18] (the anisotropic case). It is worth noting that the presented cylindrical-coordinate-based SEM requires 76% less DoFs when compared with the approach in [18].

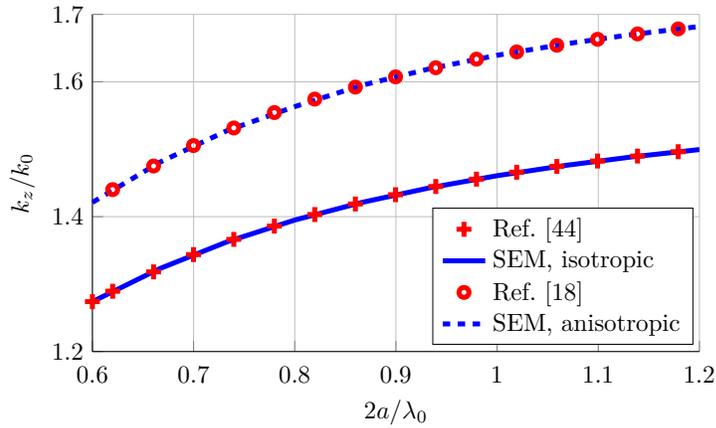

Figure 2.8: Normalized axial wavenumber for the hybrid $HE_{11}$ mode of a radially-unbounded circular waveguide.

### 2.6.3
### Anisotropic and Multilayered Coaxial Waveguides

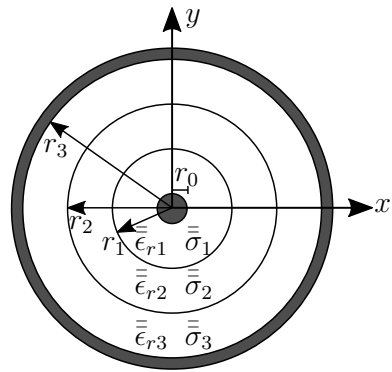

Figure 2.9: A coaxial waveguide is composed of three layers.

Next, we investigate the ability of our SEM formulation to solve complex media (lossy, inhomogeneous, and anisotropic) waveguide problems. We consider a coaxial waveguide with three layers defined by the radii $r_0 = 0.1$ m, $r_1 = 0.4$ m, $r_2 = 0.7$ m, $r_3 = 1.0$ m, and truncated internally at $r_0$ and externally at $r_3$ by PEC. The geometry of the problem is depicted in Fig. 2.9. To showcase the capabilities of the present SEM formulation, we will analyze three different scenarios where the media is non-reciprocal and described through non-symmetric and non-Hermitian tensors. The media considered are characterized by the following relative permittivity and electrical conductivity tensors represented in Cartesian coordinates:

$$\bar{\bar{\epsilon}}_{a,rec} = \begin{bmatrix} 3.5 & 2.5 & 0 \\ 2.3 & 2.9 & 0 \\ 0 & 0 & 3.1 \end{bmatrix}, \ \bar{\bar{\sigma}}_{a,rec} = 10^{-3} \begin{bmatrix} 7 & 4 & 0 \\ 3 & 5 & 0 \\ 0 & 0 & 6 \end{bmatrix} \quad (2\text{-}43a)$$

$$\bar{\bar{\epsilon}}_{b,rec} = \begin{bmatrix} 2.4 & 1.9 & 0 \\ 2.0 & 2.2 & 0 \\ 0 & 0 & 2.8 \end{bmatrix}, \ \bar{\bar{\sigma}}_{b,rec} = 10^{-3} \begin{bmatrix} 4 & 1 & 0 \\ 2 & 3 & 0 \\ 0 & 0 & 5 \end{bmatrix} \quad (2\text{-}43b)$$

$$\bar{\bar{\epsilon}}_{c,rec} = \begin{bmatrix} 2.1 & 1.4 & 0 \\ 1.3 & 2.2 & 0 \\ 0 & 0 & 1.9 \end{bmatrix}, \ \bar{\bar{\sigma}}_{c,rec} = 10^{-3} \begin{bmatrix} 3 & 1 & 0 \\ 2 & 2 & 0 \\ 0 & 0 & 4 \end{bmatrix} \quad (2\text{-}43c)$$

where the permeability is that of the vacuum. Since the above tensors are defined in Cartesian coordinates, it is necessary to perform the Cartesian-to-cylindrical coordinate transformation to make such representations conforming with (2-1)–(2-4) of our formulation. For case (a), we consider a homogeneous coaxial waveguide, where the parameters are defined by $\bar{\bar{\epsilon}}_{rj} = \bar{\bar{\epsilon}}_{a,rec}$ and $\bar{\bar{\sigma}}_j = \bar{\bar{\sigma}}_{a,rec}$ for all $j \in \{1, 2, 3\}$. For case (b), we consider a two-layer inhomogeneous coaxial waveguide with $\bar{\bar{\epsilon}}_{r1} = \bar{\bar{\epsilon}}_{a,rec}$, $\bar{\bar{\sigma}}_1 = \bar{\bar{\sigma}}_{a,rec}$, $\bar{\bar{\epsilon}}_{rj} = \bar{\bar{\epsilon}}_{b,rec}$ and $\bar{\bar{\sigma}}_j = \bar{\bar{\sigma}}_{b,rec}$ for $j \in \{2, 3\}$. For case (c), we consider a three-layer inhomogeneous coaxial waveguide, where the first layer is defined by $\bar{\bar{\epsilon}}_{r1} = \bar{\bar{\epsilon}}_{a,rec}$ and $\bar{\bar{\sigma}}_1 = \bar{\bar{\sigma}}_{a,rec}$, the second layer by $\bar{\bar{\epsilon}}_{r2} = \bar{\bar{\epsilon}}_{b,rec}$ and $\bar{\bar{\sigma}}_2 = \bar{\bar{\sigma}}_{b,rec}$, the outermost layer by $\bar{\bar{\epsilon}}_{r3} = \bar{\bar{\epsilon}}_{c,rec}$ and $\bar{\bar{\sigma}}_3 = \bar{\bar{\sigma}}_{c,rec}$. Assuming an operating frequency of $f = 100$ MHz, the results for the axial wavenumber $k_z$ from our SEM formulation are compared with FEM results [2] in the complex region $[0, 5]^2$. Excellent agreement is observed in Fig. 2.10. In this case, the computational domain is discretized into 12 elements with expansion functions of order 8, which implies 1348 DoF. In the FEM COMSOL model, we use a cubic-order discretization with the *extremely fine physics-controlled mesh* for obtaining accurate reference solutions. This results in 36858 DoF. It is worth noting that the presented formulation SEM requires about 96 % fewer DOFs when compared with this COMSOL solution.

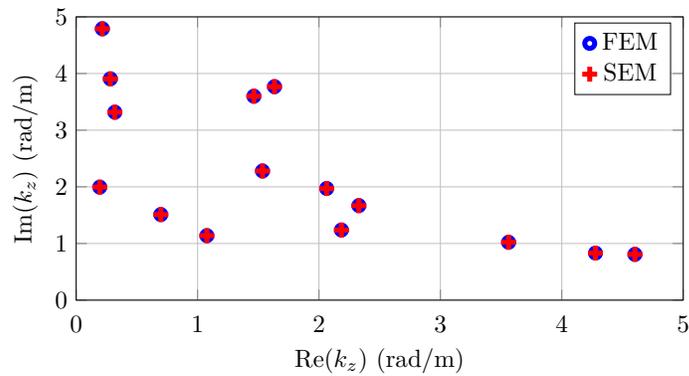

(a) One anisotropic layer

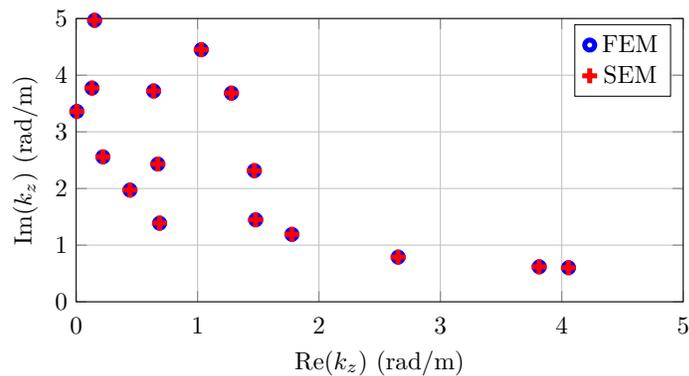

(b) Two anisotropic layers

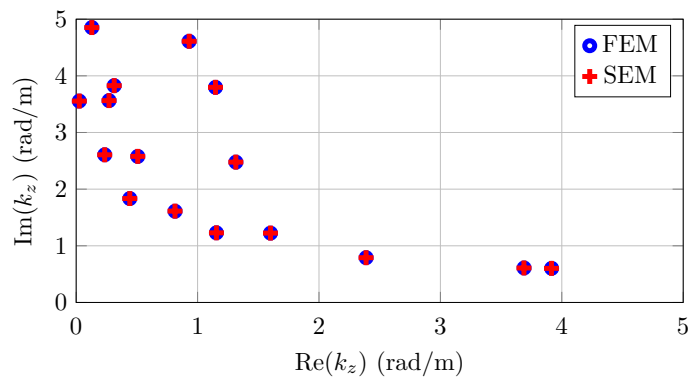

(c) Three anisotropic layers

Figure 2.10: Axial wavenumbers obtained over the complex region $[0, 5]^2$ of the multilayered coaxial waveguides filled with lossy and anisotropic media.

## 2.7
## Final Remarks

In this work, we investigated a novel SEM-based technique for solving the modal fields in cylindrical waveguides considering: (a) radially-bounded and radially-unbounded domains, (b) homogeneous and inhomogeneous filled waveguides, (c) lossless and lossy problems, and (d) isotropic and anisotropic media. The results obtained showed excellent agreement with analytical results from the literature and with FEM results. Our high-order expansion for the basis functions did not present numerical errors due to the Runge effect, and the observed convergence is exponential. The cylindrical-coordinate-based SEM represents a simplification for waveguides that present cylindrical-conforming boundaries. The present formulation requires fewer elements and DoFs when compared with the conventional Cartesian-coordinate-based SEM [18] and FEM [2].

# 3

# Higher-Order SEM for Rapid Analysis of Eccentric Coaxial Waveguides Filled with Lossy Anisotropic Media

## 3.1
## Introduction

The analysis of field propagation in non-concentric waveguides filled with lossy anisotropic media is essential to many problems, including sensing probes, microwave filters, geophysics exploration sensors, and the medical industry [9, 11, 21, 46]. Various approaches have been developed to compute wavenumbers and field patterns in these geometries. In [47], the cutoff wavenumber of an eccentric waveguide was calculated using translation addition theorems for cylindrical harmonic functions. In [48], a similar formulation was applied, combining cylindrical functions with sine and cosine laws to model the effect of small eccentricities. A related formulation was developed in [49] to calculate the electromagnetic scattering of the eccentric coating metamaterial. In [50], three methods were combined: conformal mapping, the method of intermediate problems, and the Rayleigh-Ritz method to produce rigorous results for eccentric and lunar eccentric waveguides. In [51], the cutoff wavenumbers for higher-order transverse magnetic (TM) and transverse electric (TE) modes were calculated using a derivation of addition's theorems from a previous work [52]. The contour integral method was then combined with Muller's root-finding to obtain the eigenvalues. In [53], a perturbation method was applied to derive analytical expressions for the cutoff wavenumber in circular and rectangular waveguides with small eccentric inner conductors.

Different approaches using numerical techniques have also been developed to tackle this type of problem. In [54, 55], cutoff wavenumbers were calculated using point-matching and finite-difference methods, respectively. In [56], the eigenvalue problem was solved via a meshless numerical method based on a combination of the fundamental and particular solutions for the Helmholtz equation. In [57], the Rayleigh-Ritz method was applied using polynomial approximation and superquadric functions to study waveguides with complex cross-sections, including eccentric waveguides. Related formulation using determinantal equations and Green's functions can be found in [58]. In the context of waveguides filled with complex media, as presented in [59], the cutoff wavenumbers of circular waveguides with homogeneous uniaxially and gyroelectric anisotropic material were determined using a volume integral equation and entire domain vectorial basis functions. In [60], the auxiliary-current vector norm-based method of auxiliary sources has been extended and applied to compute cutoff wavenumbers in waveguides of arbitrary shapes. In [61], the Helmholtz equation was solved in a bipolar coordinate system (BCS) by using the method of separation of variables. In the past, conformal mapping has proved robust for solving various wave-field problems. In [62], TE and TM modes were calculated in transmission lines with an eccentric circular inner conductor. In [63], the conformal mapping was used to transform an eccentric coaxial waveguide into an equivalent concentric problem, and the cutoff wavenumbers for TM and TE modes were found from the solution of the weighted Helmholtz equation. The concentric problem was solved by finite-difference in [64]. In [65], the eccentric domain was mapped into the rectangular domain, and the wavenumbers were obtained via the finite-element method (FEM). Similarly, in [12], the transformation optics (TO) was employed to map the eccentric problem into a concentric problem, and the cutoff wavenumbers and field patterns were obtained using perturbation techniques.

The spectral element method (SEM) has been successfully applied to various two-dimensional (2D) electromagnetic problems [18,19,22,23,25,66,67]. SEM-based methods are commonly used in conjunction with numerical mode matching (NMM) to treat three-dimensional problems. In [66,67], the NMM was applied to model cylindrically layered structures. Within the scope of modeling each axial layer, the authors employ the harmonic function in the azimuthal direction while resorting to a second-order one-dimensional (1D) SEM (for the second-order, the distribution nodes are uniform) for representing the field characteristics in the radial direction. In [18,19], a 2D-SEM formulation uses an expansion of fields in rectangular coordinates for modeling lossy anisotropic waveguides and problems with periodic boundary conditions. In [22,23], the 2D-SEM in cylindrical coordinates was applied to modeling multilayered concentric coaxial waveguides where a vector formulation (with coupled TM and TE modes) was required. The work [22] indicated that the 2D-SEM formulation based on cylindrical coordinates requires significantly fewer elements and degrees of freedom (DoF) compared to the conventional Cartesian approach [18] for cylindrical-conforming geometries. In this work, we expand our prior research in [25] by defining a normalized scalar Helmholtz equation in terms of decoupled TM and TE modes that enable the develop-ment of a high-speed technique for analyzing eccentric waveguides filled with uniaxially anisotropic media. With this in mind, the present study applies transformation optics to map the original eccentric waveguide problem into an equivalent concentric one, similar to the approach in [12]. Subsequently, the higher-order cylindrical SEM is employed to solve the mapped problem. This approach offers the advantage of simultaneously modeling a wide class of anisotropic media. The approach is not limited to small eccentricities, as the perturbation solution provided in [12]. The presented SEM formulation offers advancements compared to conventional 2D-SEM formulations for modeling eccentric problems, primarily due to the use of the TO, we can work with a coordinate system that conforms to the geometry. Furthermore, unlike the studies in [66, 67], our approach does not assume azimuthal symmetry. In contrast to the works presented in [22, 23], the present formulation can handle decoupled TM and TE modes, a feature not explored in the earlier studies, which confers computational efficiency to the method presented herein.

The main contribution of this work is to provide an accurate and numerically efficient technique that explores higher-order basis functions within a variational 2D-SEM formulation in cylindrical coordinates with exponential convergence characteristics that allow the modeling of eccentric coaxial waveguides filled with lossy anisotropic media. The remainder of this work is organized as follows. In Section 3.2, we present the problem statement based on the wave equations for TM and TE modes after employing the conformal TO. We introduce a novel variational formulation for axial electric and magnetic fields, providing a solution for a wide class of anisotropic medium scenarios. This formulation is particularly useful for addressing problems involving the evaluation of many different anisotropy configurations (such as for some inverse problems). Once a given scenario is solved, other scenarios can be promptly solved through a straightforward eigenvalue normalization, avoiding the need to solve the entire problem again from the start. In Section 3.3, we provide details of the SEM. We discuss the discretization of the problem domain into reference elements, the solution within a reference element, the selection of Lagrange basis functions, and the construction of the associated eigenvalue problem. In Section 3.4, we explore several examples to validate the present SEM against perturbational, analytical Graf's addition theorem (GAT), finite-element method (FEM), and finite-integration technique (FIT) solutions. We compare the eigenvalues and field representations in waveguides with small and large eccentricities. Numerical results demonstrate that our formulation is free from spurious modes, mitigates the Runge effect, and achieves high accuracy with a relatively small number of basis functions. Finally, in Section 3.5, we present concluding remarks.

## 3.2
## Variational Waveguide Problem

Consider an eccentric coaxial waveguide invariant along the axial direction, radially bounded by a perfect electric conductor (PEC) of radius $\tilde{r}_1$. The inner conductor is a circular PEC cylinder of radius $\tilde{r}_0$ with an offset $\tilde{d}$ from the $\tilde{z}$-axis. The eccentric geometry is referred to a cylindrical coordinate system denoted by $(\tilde{\rho}, \tilde{\phi}, \tilde{z})$, with the associated cross-section $\tilde{\Omega}$ depicted in Fig. 3.1(a). We assume a time-harmonic dependence of the form $e^{-i\omega t}$ and that the waveguide is filled with lossy uniaxially anisotropic media characterized by real-valued permeability tensor[1]

$$\tilde{\bar{\bar{\mu}}} = \mu_0 \tilde{\bar{\bar{\mu}}}_r \quad \text{with} \quad \tilde{\bar{\bar{\mu}}}_r = \text{diag}(\tilde{\mu}_{rs}, \tilde{\mu}_{rs}, \tilde{\mu}_{rz}), \tag{3-1}$$

and the complex-valued permittivity tensor

$$\tilde{\bar{\bar{\epsilon}}} = \epsilon_0 \tilde{\bar{\bar{\epsilon}}}_r + \frac{i}{\omega} \tilde{\bar{\bar{\sigma}}}, \tag{3-2}$$

with

$$\tilde{\bar{\bar{\epsilon}}}_r = \text{diag}(\tilde{\epsilon}_{rs}, \tilde{\epsilon}_{rs}, \tilde{\epsilon}_{rz}), \quad \tilde{\bar{\bar{\sigma}}} = \text{diag}(\tilde{\sigma}_s, \tilde{\sigma}_s, \tilde{\sigma}_z), \tag{3-3}$$

where $\tilde{\bar{\bar{\sigma}}}$ represents the conductive tensor, and $\omega = 2\pi f$ is angular frequency.

From conformal transformation optics principles [10, 36, 41], the original electromagnetic problem in the eccentric coordinates $(\tilde{\rho}, \tilde{\phi}, \tilde{z})$ can be transformed into a concentric problem with coordinates $(\rho, \phi, z)$ using

$$\mathbf{F} = \bar{\bar{J}}_{\text{TO}} \cdot \tilde{\mathbf{F}}, \quad \text{with} \ \ \mathbf{F} \in \{\mathbf{E}, \mathbf{H}\}, \tag{3-4}$$

$$\bar{\bar{p}} = |\bar{\bar{J}}_{\text{TO}}|^{-1} \bar{\bar{J}}_{\text{TO}} \cdot \tilde{\bar{\bar{p}}} \cdot \bar{\bar{J}}_{\text{TO}}^T, \quad \text{with} \ \ p \in \{\mu, \epsilon\}, \tag{3-5}$$

where $\bar{\bar{J}}_{\text{TO}}$ is the Jacobian of the transformation $(\tilde{\rho}, \tilde{\phi}, \tilde{z}) \rightarrow (\rho, \phi, z)$, and $|\bar{\bar{J}}_{\text{TO}}|$ is the respectively determinant. Assuming $z = \tilde{z}$ and $r_1 = \tilde{r}_1$, and uniaxially anisotropic media, we can simplify (3-5) and express the transformed constitutive tensor as

$$\bar{\bar{p}} = \begin{bmatrix} p_s & 0 & 0 \\ 0 & p_s & 0 \\ 0 & 0 & p_z \end{bmatrix} = \begin{bmatrix} \tilde{p}_s & 0 & 0 \\ 0 & \tilde{p}_s & 0 \\ 0 & 0 & |\bar{\bar{J}}_{\text{TO}}|^{-1}\tilde{p}_z \end{bmatrix}, \tag{3-6}$$

with $p \in \{\mu, \epsilon\}$, where

$$|\bar{\bar{J}}_{\text{TO}}|^{-1} = \frac{(1 - \tilde{x}_1/\tilde{x}_2)^2}{(1 - 2\rho\cos\phi/\tilde{x}_2 + \rho^2/\tilde{x}_2^2)^2}, \tag{3-7}$$

---

[1]Complex permeability can also be considered by allowing $\tilde{\mu}_{rs}$ and $\tilde{\mu}_{rz}$ to be complex-valued parameters.

and

$$\tilde{x}_{1,2} = \frac{-\tilde{c} \mp \sqrt{\tilde{c}^2 - 4\tilde{r}_1^2}}{2}, \quad \text{with} \quad \tilde{c} = \frac{\tilde{r}_0^2 - \tilde{r}_1^2 - \tilde{d}^2}{\tilde{d}}. \tag{3-8}$$

In the mapped domain, the axial constitutive parameters $\mu_z$ and $\epsilon_z$ exhibit variations in $\rho$ and $\phi$ directions. This is in contrast to the *constant* tensors assumed in (3-1) and (3-2).

To solve the concentric problem in the cross-section $\Omega$ depicted in Fig. 3.1(b), the electromagnetic fields can be decomposed into TM and TE modes [43]. We first solve for the axial field components $E_z$ and $H_z$, then compute the fields' transverse components in cylindrical coordinates. In concentric domain $\Omega$, the wave equation for $E_z$ and $H_z$ satisfies [68, Ch. 2]

$$\left( \nabla_s^2 + \frac{p_z(\rho, \phi)}{p_s} k_\rho^2 \right) F = 0, \tag{3-9}$$

where $\nabla_s^2$ is the transverse (to $z$) Laplacian operator, $k_\rho$ is the radial wavenumber defined by $k_\rho^2 = k_s^2 - k_z^2$ with $k_s^2 = \omega^2 \mu_s \epsilon_s$, and $F = E_z$ if $p_{s,z} = \epsilon_{s,z}$, and $F = H_z$ if $p_{s,z} = \mu_{s,z}$. The domain is truncated by a PEC; then the fields must satisfy $F = 0$ for $F = E_z$ and $\partial F/\partial \rho = 0$ for $F = H_z$ on $\rho \in \{r_0, r_1\}$.

We can use the Galerkin method to formulate a variational formulation for problem (3-9), which involves decomposing the fields into a finite set of basis functions. The resulting equation is given by

$$\langle \boldsymbol{\nabla}_s F, \boldsymbol{\nabla}_s W \rangle - k_\rho^2 \left\langle \frac{p_z(\rho, \phi)}{p_s} F, W \right\rangle = 0, \tag{3-10}$$

where $\langle \cdot, \cdot \rangle$ denotes an inner product in Hilbert space, $F$ is either $E_z$ or $H_z$, and $W$ is a test function in the space $H_0^1(\Omega)$ [35, 69, 70].

From (3-6), using $p_s = \tilde{p}_s$ and $p_z = \tilde{p}_z |\bar{\bar{J}}_{\text{TO}}|^{-1}$, we can normalize the media parameters and rewrite (3-10) in the following form

$$\langle \boldsymbol{\nabla}_s F, \boldsymbol{\nabla}_s W \rangle - \lambda \left\langle |\bar{\bar{J}}_{\text{TO}}|^{-1} F, W \right\rangle = 0, \tag{3-11}$$

where $\lambda = (\tilde{p}_z/\tilde{p}_s) k_\rho^2$. For uniaxially anisotropic media, it is important to notice that the axial field distribution in the $\Omega$-domain is solely determined by the problem's geometry. Therefore, once the values of $\tilde{r}_1$, $\tilde{r}_0$, and $\tilde{d}$ are defined, it is obtained an entire set of solutions that is independent of the parameters of the medium.

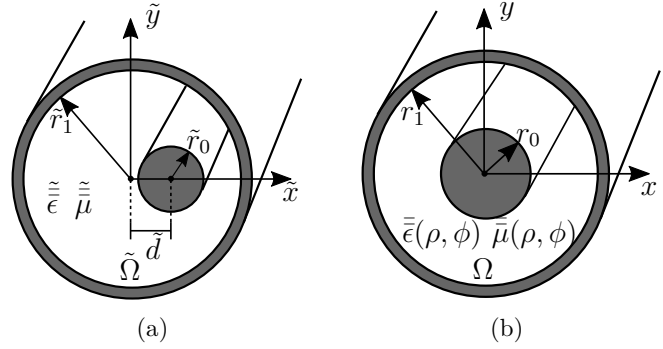

Figure 3.1: (a) Geometry of an eccentric coaxial waveguide invariant along the axial direction. (b) Geometry of the mapped concentric waveguide.

### 3.3
### Spectral Element Method

The spectral element method is an advanced implementation of the FEM, which also involves partitioning the computational domain into discrete elements. However, in SEM, the internal nodes of each element are located at the zeros of specific orthogonal polynomial families. This feature enables the avoidance of the Runge effect, which refers to large oscillations that occur in the polynomial approximation of regions near the edges of a surface in higher-order interpolations. As a result, SEM can achieve higher accuracy with fewer elements compared to a corresponding FEM implementation. This section presents an improvement formulation in cylindrical coordinates for the spectral element method modeling of coaxial waveguides.

### 3.3.1
### Discretized Domain

The cross-section $\Omega$ of the coaxial computational domain is first discretized into four elements $\Omega^e$, with $e \in \{1, 2, 3, 4\}$, as shown in Fig. 3.2. Each element is described by the polar coordinates $(\rho^e, \phi^e)$, which is then mapped into the reference coordinates $(\xi, \eta)$ via

$$\rho^e(\xi) = 0.5(\rho_0^e + \rho_1^e) + 0.5(\rho_1^e - \rho_0^e)\xi, \quad \forall \xi \in [-1, 1], \tag{3-12a}$$

$$\phi^e(\eta) = 0.5(\phi_0^e + \phi_1^e) + 0.5(\phi_1^e - \phi_0^e)\eta, \quad \forall \eta \in [-1, 1]. \tag{3-12b}$$

The superscript $e$ denotes the reference element parameters. This mapping allows us to relate the transverse points of the reference element covering $-1 \leq \xi \leq 1$ and $-1 \leq \eta \leq 1$ with the physical element covering $\rho_0^e \leq \rho^e \leq \rho_1^e$ and $\phi_0^e \leq \phi^e \leq \phi_1^e$. Fig. 3.3 illustrates the mapping between the reference element and a given physical element. The Jacobian matrix associated with

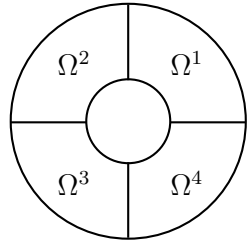

Figure 3.2: Partition of the computational domain into four elements.

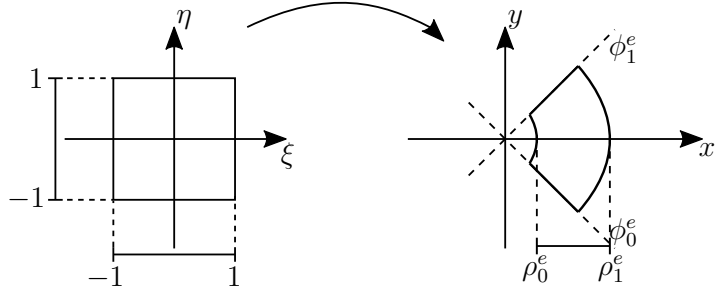

Figure 3.3: Mapping between the reference element and a given physical element.

this coordinate transformation is given by

$$\bar{\bar{J}}_s^e = \begin{bmatrix} \frac{\partial \rho^e}{\partial \xi} & \frac{\partial \rho^e}{\partial \eta} \\ \rho^e \frac{\partial \phi^e}{\partial \xi} & \rho^e \frac{\partial \phi^e}{\partial \eta} \end{bmatrix} \tag{3-13}$$

$$= \begin{bmatrix} 0.5(\rho_1^e - \rho_0^e) & 0 \\ 0 & 0.5\rho^e(\phi_1^e - \phi_0^e) \end{bmatrix}.$$

The above transformation matrix is diagonal, and as a result, the $\rho$- and $\eta$-directions are decoupled, as well as $\phi$ and $\xi$ directions [22, 23].

### 3.3.2
### Lagrange Interpolation

We employ Lagrange basis functions associated with Gauss-Lobatto-Legendre (GLL) sampling points to represent the field components in the reference coordinates. Specifically, for the coordinate $\xi$ with sampling points $\xi_k \in [-1, 1]$ for all $j = 1, 2, ..., N + 1$, we use the $N$th-order one-dimensional GLL basis function given by

$$\phi_j^N(\xi) = \prod_{k \neq j}^{N+1} \frac{\xi - \xi_k}{\xi_j - \xi_k} \quad \forall \xi \in [-1, 1], \tag{3-14}$$

where the $N + 1$ interpolation nodes $\xi_k$ are the roots of the $(N + 1)$-degree completed Lobatto polynomial $(1 - \xi^2)Lo_{N-1}(\xi) = 0$. A quick examination

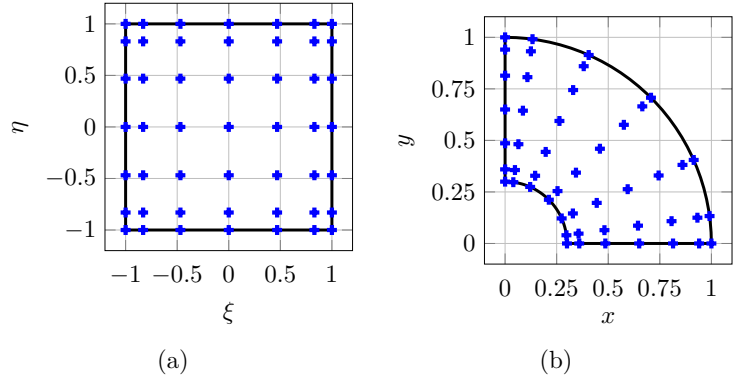

Figure 3.4: Node distribution in the (a) reference and (b) physical element, using expansion order $N = 6$.

reveals that these polynomials satisfy the cardinal interpolation property given by

$$\phi_j^N(\xi_k) = \begin{cases} 1 & \text{if } j = k \\ 0 & \text{if } j \neq k \end{cases}. \tag{3-15}$$

In other words, $\phi_j^N(\xi_k)$ acts as the Kronecker delta, i.e., $\phi_j^N(\xi_k) = \delta_{jk}$ [28].

We use tensor-product expansions to obtain the two-dimensional GLL basis function

$$\hat{\psi}_{ij}(\xi, \eta) = \phi_i^N(\xi)\phi_j^N(\eta) \ \ \forall (\xi, \eta) \in [-1, 1]^2. \tag{3-16}$$

To obtain the corresponding basis function $\psi_{ij}^e$ and its gradient in the physical element, we use the covariant mappings via [35]

$$\begin{cases} \psi_{ij}^e(\rho^e, \phi^e) = \hat{\psi}_{ij}(\xi, \eta) \\ \boldsymbol{\nabla}_s \psi_{ij}^e(\rho^e, \phi^e) = \bar{\bar{J}}_s^{e^{-1}} \hat{\boldsymbol{\nabla}}_s \hat{\psi}_{ij}(\xi, \eta) \end{cases}, \tag{3-17}$$

where $\boldsymbol{\nabla}_s = \left(\frac{\partial}{\partial \rho^e}, \frac{1}{\rho^e}\frac{\partial}{\partial \phi^e}\right)^T$ and $\hat{\boldsymbol{\nabla}}_s = \left(\frac{\partial}{\partial \xi}, \frac{\partial}{\partial \eta}\right)^T$ are the gradient operators in the physical and reference elements, respectively. The Jacobian matrix $\bar{\bar{J}}_s^e$ is defined by (3-13). The covariant mappings allow us to relate the basis functions and their gradients in the reference element to those in the physical element. As an illustrative example, the distribution of nodes in the reference and physical elements considering expansion function order $N = 6$ are depicted in Fig. 3.4. It is worth noting that the point distribution is denser towards the edges of the element, a characteristic that effectively mitigates the occurrence of Runge's phenomenon when employing higher-order expansions [28, 29]. The choice of interpolation basis functions results in $(N + 1)^2$ degrees of freedom on each element, leading to a total of $4N(N + 1)$ DoFs for the computational domain described in Fig. 3.2.

### 3.3.3
### Eigenvalue Problem

To construct the discrete form of problem (3-11), the axial field $F$ in the reference domain is expanded as

$$F(\xi, \eta) = \sum_{n=1}^{(N+1)^2} a_n \hat{\psi}_n(\xi, \eta). \qquad (3\text{-}18)$$

Note that to simplify the notation, we are replacing the double index $ij$ in (3-16) by a single index $n$. Following the Galerkin method, test functions are identical to the basis functions, namely, $W = \hat{\psi}_m$ for $m = 1, 2, ..., (N+1)^2$. Using (3-17) and (3-18), we can define elemental matrices that represent the inner products in (3-11). These matrices can be written in reference coordinates as

$$\mathbf{D}^e(m, n) = \int_{-1}^{1} \int_{-1}^{1} (\bar{\bar{J}}_s^{e-1} \boldsymbol{\hat{\nabla}}_s \hat{\psi}_n)^T \cdot (\bar{\bar{J}}_s^{e-1} \boldsymbol{\hat{\nabla}}_s \hat{\psi}_m) |\bar{\bar{J}}_s^e| d\xi d\eta, \qquad (3\text{-}19a)$$

$$\mathbf{M}^e(m, n) = \int_{-1}^{1} \int_{-1}^{1} |\bar{\bar{J}}_{\mathrm{TO}}|^{-1} \hat{\psi}_n \cdot \hat{\psi}_m |\bar{\bar{J}}_s^e| d\xi d\eta, \qquad (3\text{-}19b)$$

where $|\bar{\bar{J}}_{\mathrm{TO}}|^{-1}$ is given in (3-7), and $|\bar{\bar{J}}_s^e|$ denotes the determinant of the Jacobian matrix in (3-13).

The integrands in (3-19) are evaluated using the Gauss-Lobatto quadrature, and we used the same GLL nodes from the expansion basis functions in the quadrature (details are presented in Section 2.3.4). The global diffusion matrix $\mathbf{D}$ and global mass matrix $\mathbf{M}$ are constructed by assembling elemental matrices $\mathbf{D}^e$ and $\mathbf{M}^e$, respectively, where $e \in \{1, 2, 3, 4\}$. It is noteworthy that the matrix $\mathbf{D}^e$ remains azimuthally invariant, indicating that it remains the same for all elements $e$. Consequently, it only needs to be computed once. Conversely, the matrix $\mathbf{M}^e$ must be calculated four times. However, it is essential to mention that $\mathbf{M}^e$ is a diagonal matrix (see Section 3.3.4), which means that filling it incurs minimal computational cost.

Once the global diffusion and mass matrices are obtained, we can define the generalized eigenvalue problem

$$\mathbf{D}\bar{a} = \lambda \mathbf{M}\bar{a}, \qquad (3\text{-}20)$$

where $\lambda$ and $\bar{a}$ represent eigenvalues and eigenvectors, respectively. The eigenvalues are associated with the radial wavenumbers via

$$k_\rho = \sqrt{\frac{\tilde{p}_s}{\tilde{p}_z}} \lambda, \qquad (3\text{-}21)$$

while the eigenvectors are associated with the axial field coefficients in (3-18).

The axial wavenumbers can be obtained in the form

$$k_z = \sqrt{k_s^2 - k_\rho^2}. \tag{3-22}$$

The linear system (3-20) has $4N(N+1)$ DoFs for the case of $F = H_z$. In the $F = E_z$ case, the PEC boundary condition $E_z = 0$ must be satisfied for $\rho \in \{r_0, r_1\}$, which implies in $4N(N-1)$ DoFs.

### 3.3.4
### Elemental Mass Matrix

The use of a node distribution equal to that of the basis functions in the Lobatto quadrature introduces a significant numerical simplification: the elementary mass matrix defined in (3-19b) becomes diagonal. In order to show this, we evaluate the integrand in (3-19b) using (2-32), yielding

$$\mathbf{M}^e(m,n) = \sum_{p=1}^{N+1} \sum_{q=1}^{N+1} |\bar{\bar{J}}_{\text{TO}}(\rho^e(\xi_p), \phi^e(\eta_q))|^{-1} \hat{\psi}_n(\xi_p, \eta_q) \tag{3-23a}$$
$$\cdot \hat{\psi}_m(\xi_p, \eta_q) |\bar{\bar{J}}_s^e(\xi_p, \eta_q)| w_p w_q.$$

By using (3-16), the above becomes

$$\mathbf{M}^e(m,n) = \sum_{p=1}^{N+1} \sum_{q=1}^{N+1} |\bar{\bar{J}}_{\text{TO}}(\rho^e(\xi_p), \phi^e(\eta_q))|^{-1} \phi_i^N(\xi_p) \phi_j^N(\eta_q)$$
$$\cdot \phi_{i'}^N(\xi_p) \phi_{j'}^N(\eta_q) |\bar{\bar{J}}_s^e(\xi_p, \eta_q)| w_p w_q, \tag{3-23b}$$

where the index $n$ runs over the indices $(i,j)$ and $m$ over $(i',j')$. Specifically, $n = N(j-1) + i$ and $m = N(j'-1) + i'$. Next, using (3-15), we obtain

$$\mathbf{M}^e(m,n) = \sum_{p=1}^{N+1} \sum_{q=1}^{N+1} |\bar{\bar{J}}_{\text{TO}}(\rho^e(\xi_p), \phi^e(\eta_q))|^{-1} \delta_{ip} \delta_{jq} \tag{3-23c}$$
$$\cdot \delta_{i'p} \delta_{j'q} |\bar{\bar{J}}_s^e(\xi_p, \eta_q)| w_p w_q,$$

where $\delta_{rs}$ represents the Kronecker delta. Note that the sum has non-zero only terms for $(i,j) = (i',j') = (p,q)$. Thus, the elemental mass matrix can be defined as

$$\mathbf{M}^e(m,n) = \begin{cases} |\bar{\bar{J}}_{\text{TO}}(\rho^e(\xi_i), \phi^e(\eta_j))|^{-1} |\bar{\bar{J}}_s^e(\xi_i, \eta_j)| w_i w_j & \text{if } n = m \\ 0 & \text{if } n \neq m \ . \end{cases} \tag{3-23d}$$

This result significantly simplifies the numerical solution and reduces the computational cost.

## 3.4
## Numerical Results

This section shows numerical examples of solving eccentric coaxial waveguide problems using the above-presented SEM formulation. Firstly, we consider scenarios with small eccentricities and compare our SEM results against analytical solutions obtained from Graf's addition theorem (GAT) (detailed in Appendix A) and perturbative-based techniques described in [12], and the numerical solutions from the FEM and FIT solvers of CST [71]. Secondly, we investigate more complex waveguide structures where both the eccentricity offset and the inner radius vary. To validate our method, we compare the obtained results with FIT reference solutions from CST. Lastly, we explore waveguides filled with lossy anisotropic media, and the results are validated against FIT. The algorithms associated with the GAT and SEM methods were implemented using the MATLAB platform [42].

### 3.4.1
### Small Eccentricity

For our initial validation, we examine an eccentric waveguide characterized by $\tilde{r}_1 = 5$ mm, $\tilde{r}_0 = 0.05\tilde{r}_1$, and offset $\tilde{d} = 0.05\tilde{r}_1$, while assuming the medium is the vacuum. This problem was previously studied in [12,61] and is selected here for convergence analysis of our SEM-based technique. We want to verify the relative error $\varepsilon$ behavior as a function of the expansion functions orders $N$. Such error is defined in terms of the numerical SEM solution $f_{\text{SEM}}$ and analytical solution $f_{\text{GAT}}$ by Graf's addition theorem via

$$\varepsilon = \frac{|f_{\text{SEM}} - f_{\text{GAT}}|}{|f_{\text{GAT}}|}. \tag{3-24}$$

The relative errors of the first TE and TM modes are shown in Fig. 3.5. It can be observed that the error exponentially decreases as the polynomial expansion order of the SEM increases. The absence of numerical errors due to the Runge effect, typically associated with higher-order interpolations, is evident and aligns with theoretical predictions in [28,29].

We consider the same waveguide as above and proceed to examine the accuracy of the SEM solution for different eccentricity offsets $\tilde{d} \in \{0.05\tilde{r}_1, 0.10\tilde{r}_1, 0.15\tilde{r}_1, 0.20\tilde{r}_1\}$. This particular problem was previously presented in [12,25], and we include the GAT reference solution for comparison in this study. Fig. 3.6 shows the radial wavenumber of the first three TE and TM modes as a function of the normalized offset $\tilde{d}/\tilde{r}_1$ using the cavity-material perturbation method (CMPM) [72] and the regular perturbation method (RPM) from [12], FEM from [71], and GAT and SEM present in this work. Notice that CMPM does not yield accurate results for $\tilde{d} \geq 0.1\tilde{r}_1$. In particular, the $TM_{01}$ and $TM_{21}$ results from RPM present a noticeable error for $\tilde{d} = 0.2\tilde{r}_1$. The growth in the relative error as a function of the increase in the normalized displacement $\tilde{d}/\tilde{r}_1$ of the inner conductor is a limitation of the perturbation methods. The SEM, in contrast, showed excellent agreement in all simulated cases compared to the FEM and GAT solutions. This indicates that this approach can be extended to geometries that present a more significant offset of the inner conductor.

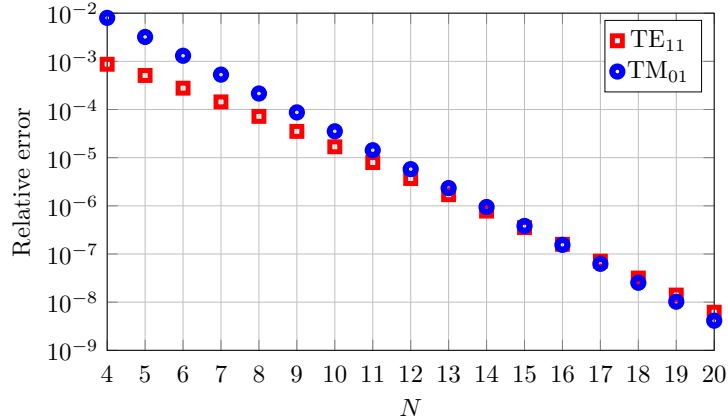

Figure 3.5: Relative error of the first TE and TM modes as a function of the expansion order $N$ used in our SEM algorithm.

Now, we investigate the ability of our SEM formulation to solve eccentric waveguides filled with anisotropic media. For this, we reproduce the cases considered in [12, 25], in which $\tilde{\epsilon}_s = \epsilon_0 \tilde{\epsilon}_{rs}$, $\tilde{\epsilon}_z = \epsilon_0 \tilde{\epsilon}_{rz}$, and $\tilde{\mu}_s = \tilde{\mu}_z = \mu_0$, where $\epsilon_0$ and $\mu_0$ are the vacuum values. The geometry has $\tilde{r}_1 = 5$ mm, $\tilde{r}_0 = 0.05\tilde{r}_1$, and offset $\tilde{d} = 0.2\tilde{r}_1$. Figs. 3.7 and 3.8 show results for the radial wavenumbers of the $TM_{01}$ mode as a function of $\tilde{\epsilon}_{rs}$ and $\tilde{\epsilon}_{rz}$. We can observe excellent agreement between dense-mesh FIT [71], GAT, and SEM. However, the CMPM and RPM present a visible error. A notable advantage of our SEM formulation is that we need to calculate the eigenvalues associated with the geometry of the problem only once. Whenever there are modifications to the parameters of the media, we can then recalculate the radial wavenumbers using (3-21).

The computational time cost required by the SEM (expansion function of order $N = 7$) is compared with that of the FIT, CMPM, RPM, and GAT in Table 3.1 for the anisotropic waveguide with $\tilde{\epsilon}_{rs} = 5$ and $\tilde{\epsilon}_{rz} = 1$ (Case I), with $\tilde{\epsilon}_{rz} = 5$ and $\tilde{\epsilon}_{rs} = 1$ (Case II). The present SEM and GAT demonstrate outstanding performance and significantly lower computational costs when compared to other techniques. The GAT and SEM algorithms were written in Matlab [42] running on a PC with 1.80 GHz Intel Core i7-8565U with four

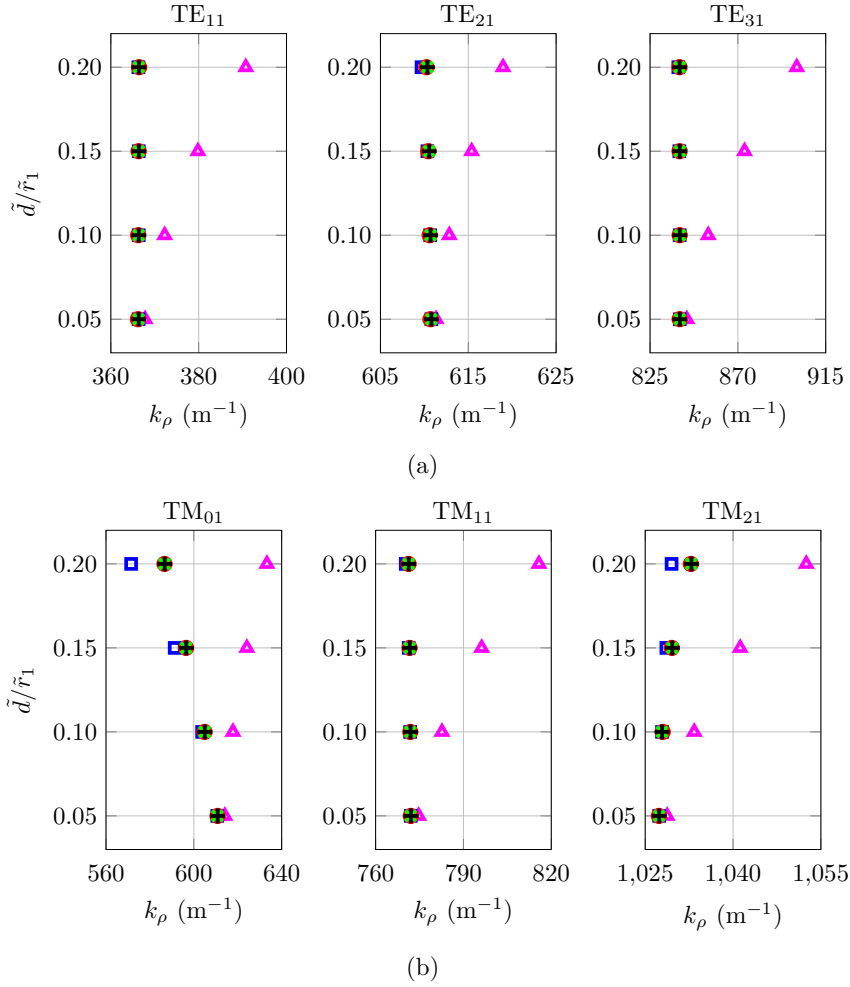

Figure 3.6: Radial wavenumbers for (a) TE and (b) TM modes as a function of the normalized eccentricity distance $\tilde{d}/\tilde{r}_1$ get by CMPM ($\triangle$), by RPM ($\square$), by FEM ($\circ$), by Graf's addition theorem ($\times$), and by SEM ($+$).

cores. The FIT, CMPM, and RMP times were obtained from [12].

Table 3.1: Computational cost in the problems with small eccentricities.

|  | FIT | CMPM | RPM | GAT | SEM |
|---|---|---|---|---|---|
| Case I | 15 min 44 s | 2.07 s | 16.72 s | 1.73 s | 0.64 s |
| Case II | 22 min 27.10 s | 2.43 s | 21.60 s | 2.33 s | 0.63 s |

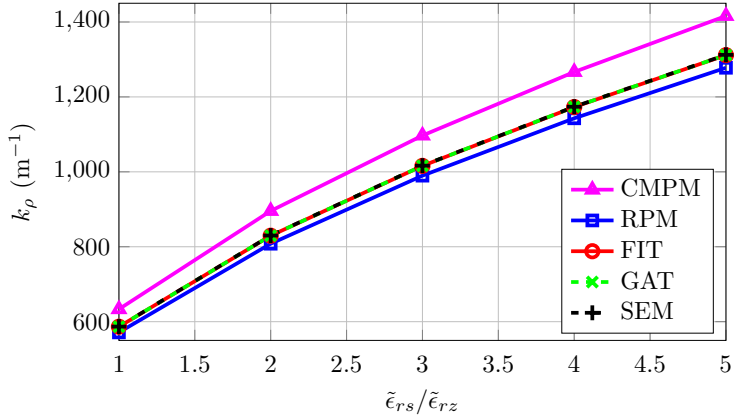

Figure 3.7: Radial wavenumbers for the $\text{TM}_{01}$ mode as a function of the anisotropic permittivity ratio $\tilde{\epsilon}_{rs}/\tilde{\epsilon}_{rz}$, with $\tilde{\epsilon}_{rz} = 1$.

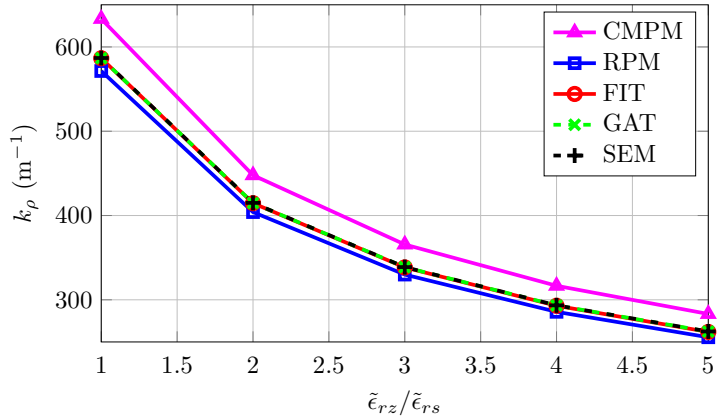

Figure 3.8: Radial wavenumbers for the $\text{TM}_{01}$ mode as a function of the anisotropic permittivity ratio $\tilde{\epsilon}_{rz}/\tilde{\epsilon}_{rs}$, with $\tilde{\epsilon}_{rs} = 1$.

### 3.4.2
### Large Eccentricity

In this section, we analyze the performance of SEM in modeling waveguides where the inner conductor has a large radius and a large offset. As a first example, we consider an eccentric coaxial waveguide characterized by $\tilde{r}_1 = 10$ mm, offset $\tilde{d} = 0.1\tilde{r}_0$, inner radius variation $\tilde{r}_0 \in \{0.1\tilde{r}_1, 0.2\tilde{r}_1, 0.3\tilde{r}_1, 0.4\tilde{r}_1, 0.5\tilde{r}_1, 0.6\tilde{r}_1, 0.7\tilde{r}_1, 0.8\tilde{r}_1\}$ and filled with lossless uniaxially anisotropic media characterized by $\tilde{\bar{\bar{\mu}}}_r = \text{diag}(2.2, 2.2, 2.7)$ and $\tilde{\bar{\bar{\epsilon}}}_r = \text{diag}(5.6, 5.6, 4.6)$. The radial wavenumbers for $\text{TE}_{11}$ and $\text{TE}_{21}$ modes are presented in Fig. 3.9 and exhibit significant agreement with FIT solutions. For this example, the expansion order $N = 7$ guarantees great accuracy of the results.

As a second example, we considered the same waveguide but with a fixed inner radius $\tilde{r}_0 = 0.1\tilde{r}_1$, and we varied the eccentricity offset $\tilde{d} \in$ $\{0.1\tilde{r}_1, 0.2\tilde{r}_1, 0.3\tilde{r}_1, 0.4\tilde{r}_1, 0.5\tilde{r}_1, 0.6\tilde{r}_1, 0.7\tilde{r}_1, 0.8\tilde{r}_1\}$. We use an expansion basis function order of $N = 10$ in this example. This larger order is required to ensure the convergence of the results, as the internal conductor has extreme displacement. Fig. 3.10 presents the results for the $\text{TM}_{01}$ and $\text{TM}_{11}$ radial wavenumbers. Fig. 3.11 show the normalized axial electric fields for different offset parameters for $\text{TM}_{01}$ mode. Additionally, Fig. 3.12 shows the normalized axial magnetic fields for $\text{TE}_{11}$ and $\text{TE}_{21}$ modes (with symmetric and antisymmetric configurations) for the case $\tilde{d} = 0.2\tilde{r}_1$. Our SEM technique demonstrates excellent agreement with the FIT solutions.

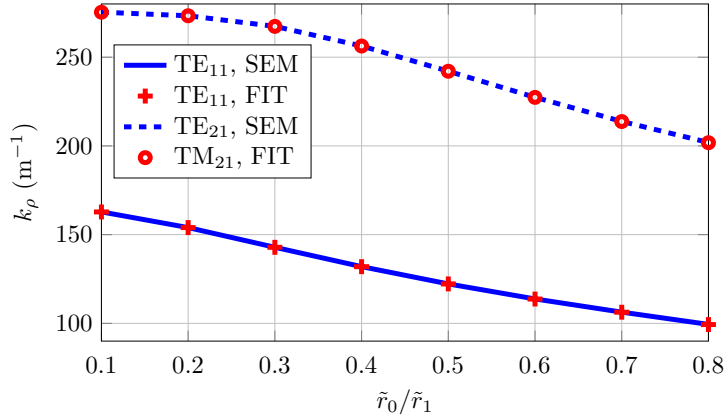

Figure 3.9: Radial wavenumbers for the $\text{TE}_{11}$ and $\text{TE}_{21}$ modes as a function of the inner ratio $\tilde{r}_0/\tilde{r}_1$ in an eccentric coaxial waveguide filled with anisotropic media.

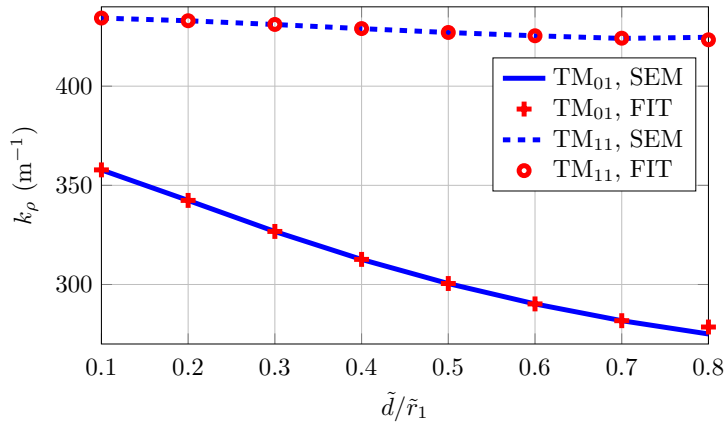

Figure 3.10: Radial wavenumbers for the $\text{TM}_{01}$ and $\text{TM}_{11}$ modes as a function of the offset $\tilde{d}/\tilde{r}_1$ in an eccentric coaxial waveguide filled with anisotropic media.

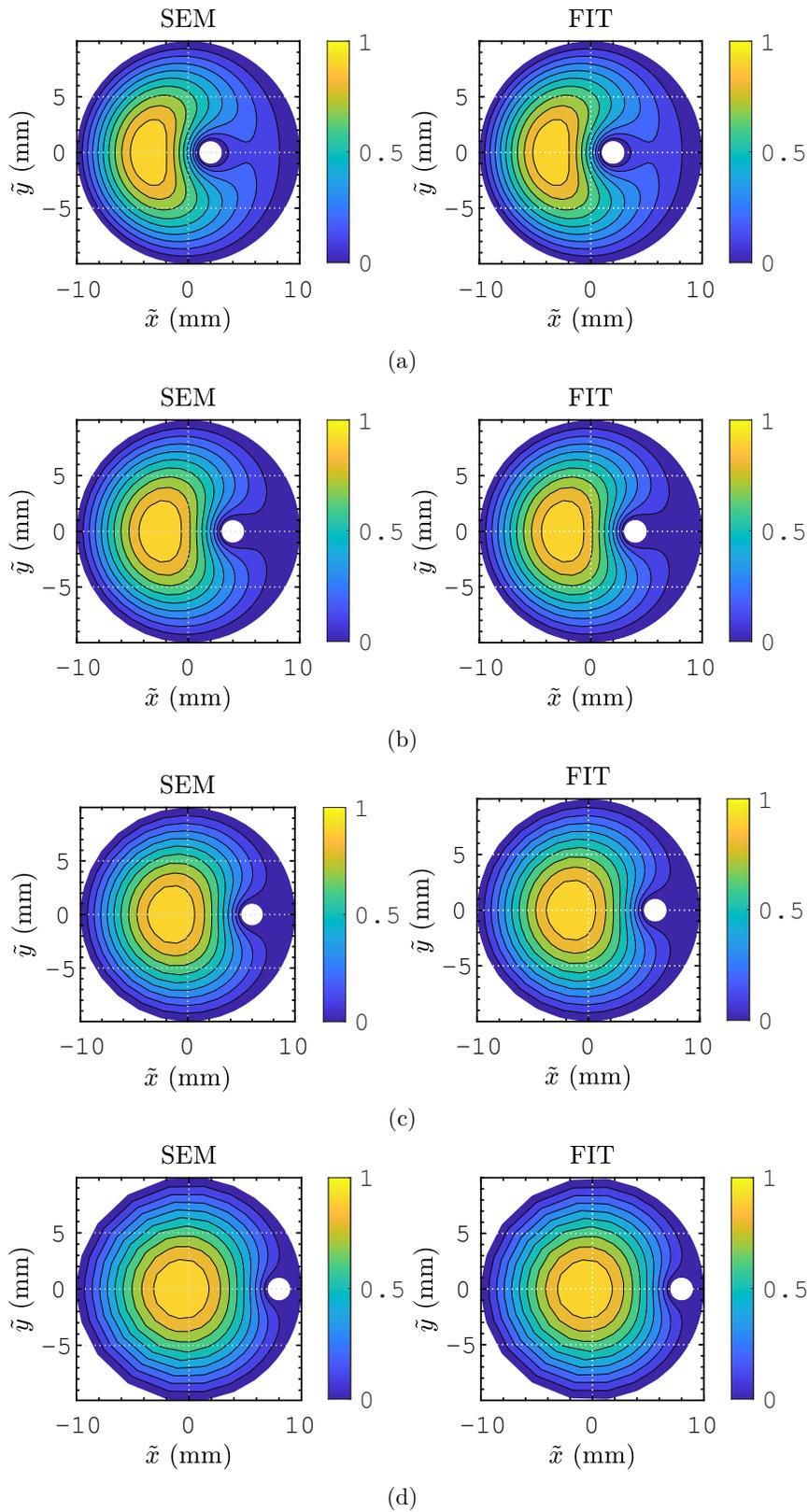

Figure 3.11: Normalized axial electric field patterns for the $TM_{01}$ mode calculated by using the SEM (left) and by the FIT (right) for different values of the eccentricities. (a) $\tilde{d} = 0.2\tilde{r}_1$. (b) $\tilde{d} = 0.4\tilde{r}_1$. (c) $\tilde{d} = 0.6\tilde{r}_1$. (d) $\tilde{d} = 0.8\tilde{r}_1$.

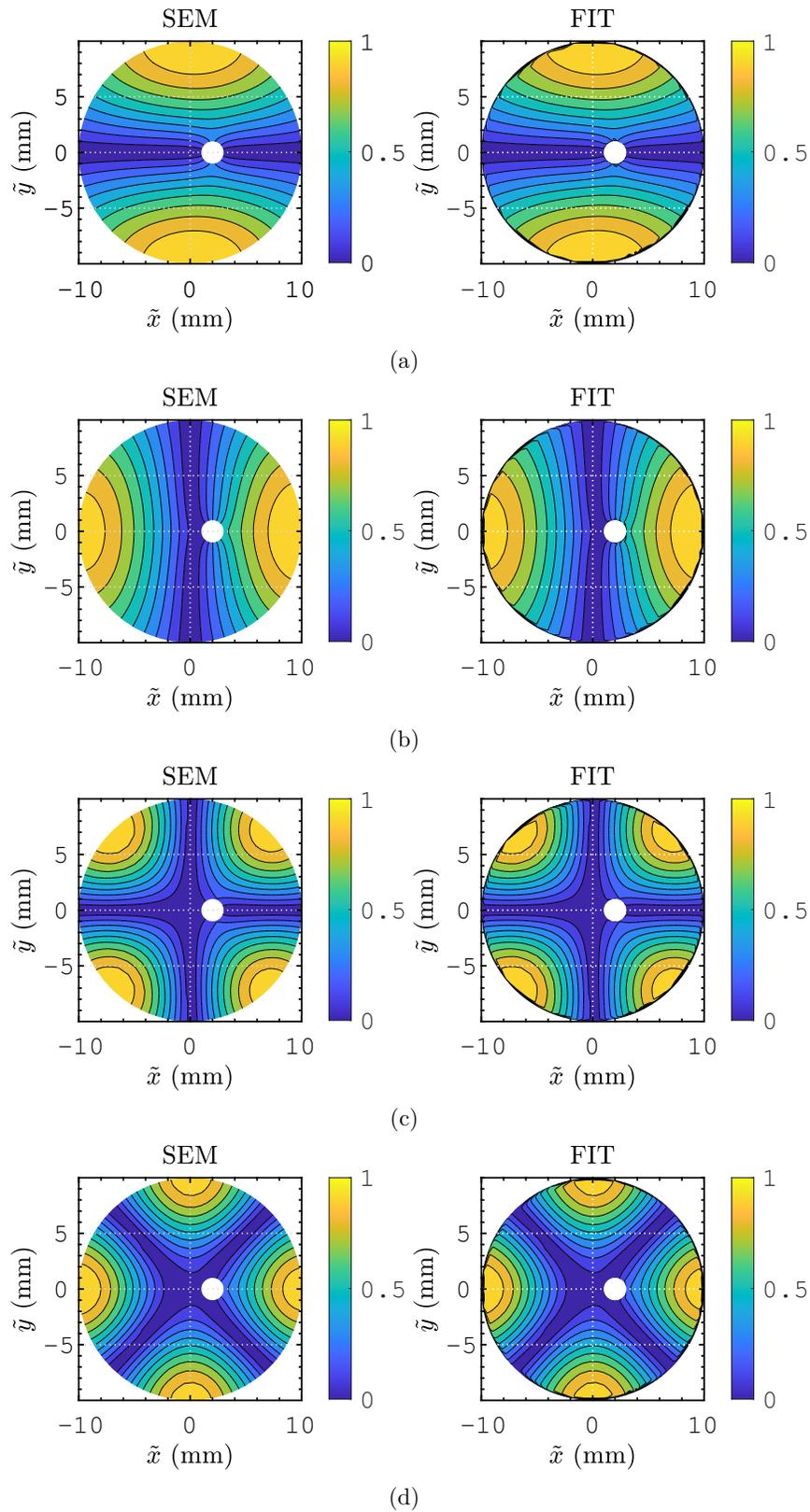

Figure 3.12: Normalized axial magnetic field patterns for the dominant TE modes were calculated using SEM (left) and FIT (right). (a) Antisymmetric $TE_{11}$ mode. (b) Symmetric $TE_{11}$ mode. (c) Antisymmetric $TE_{21}$ mode. (d) Symmetric $TE_{21}$ mode.

The computational cost of the SEM is compared to the FIT solution in Table 3.2 for the case with $\tilde{r}_0 = 0.8\tilde{r}_1$ shown in Fig. 3.9 (Case III). Additionally, a comparison is also made for the example with $\tilde{d} = 0.8\tilde{r}_1$ presented in Fig. 3.10 (Case IV). In the analysis of the computational time spent by the SEM, it should be noted that Case IV experienced a higher computational time compared to Case III. This was attributed to the necessity of employing an expansion order of $N = 10$ to ensure solution convergence. The SEM codes ran on the same computer used for Cases I and II in Section 3.4.1. All the FIT results were obtained using a dedicated HP Z800 Workstation with a dual quad-core 2.40-GHz Intel Xeon E5620 processor. Notice we have computed only the port modes by using the standard options with a discretization of 40 hexahedral cells per wavelength in the FIT solver of CST [71].

Table 3.2: Computational cost in problems with large eccentricities.

|          | FIT      | SEM     |
|----------|----------|---------|
| Case III | 10466 s  | 0.63 s  |
| Case IV  | 11606 s  | 4.15 s  |

To ensure a fair comparison, it is crucial to emphasize the number of modes represented in each numerical method. In the case of SEM, the number of modes calculated in each simulation is given by $8N^2$, which represents the sum of TE and TM modes for a given expansion order $N$. In Case III, an expansion order of $N = 7$ was employed, resulting in 392 modes. In Case IV, an expansion order of $N = 10$ was utilized, resulting in 800 modes. On the other hand, only 15 (TM+TE) modes were requested for FIT. However, due to the increased complexity of the geometry compared to Cases I and II, achieving convergence required additional computational time.

### 3.4.3
### Lossy Anisotropic Media

Finally, we investigate the ability of our SEM formulation to solve eccentric waveguides filled with lossy media. The media under investigation exhibit uniaxial anisotropy, as described in Section 3.4.2. Additionally, we introduce the conductive tensor $\bar{\bar{\sigma}} = \{0.38, 0.38, 0.34\}$ to characterize the media dissipation. The dimensions of the waveguide are defined as $\tilde{r}_1 = 10$ mm, $\tilde{r}_0 = 0.2\tilde{r}_1$, and $\tilde{d} = 0.3\tilde{r}_1$. We present the results for the real and imaginary components of the axial wavenumber $k_z$ as a function of the operating frequency $f$, ranging from 1 GHz to 10 GHz, in Fig. 3.13. For this example, we consider basis functions with order $N = 7$. The SEM results present excellent agreement with FIT results. Additionally, the normalized axial electric fields for $TM_{01}$ and $TM_{11}$, and operating frequency $f = 10$ GHz, are shown in Fig 3.14, with again agreement with the FIT solution.

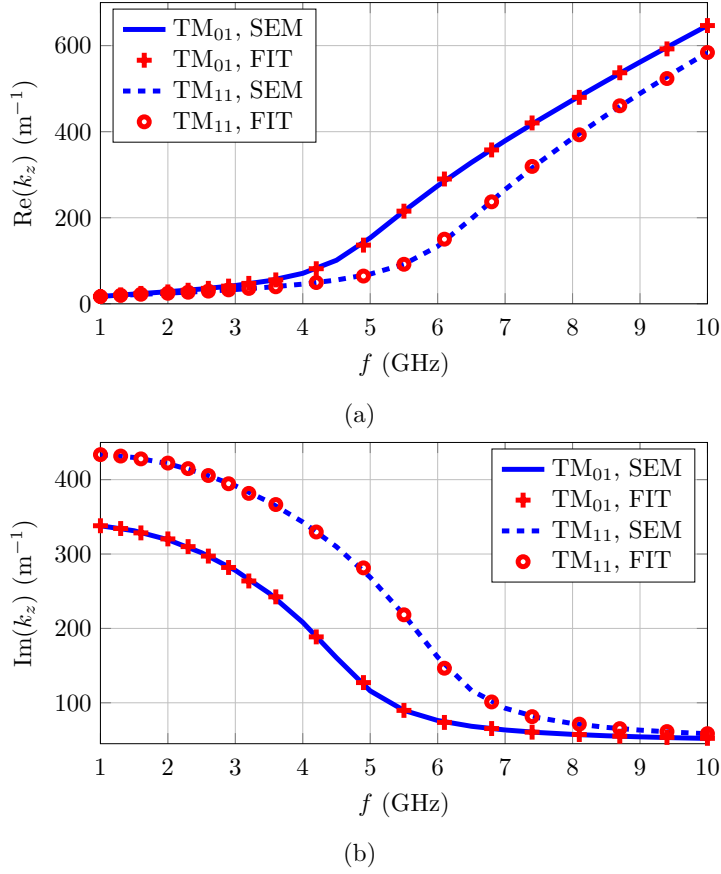

(a)

(b)

Figure 3.13: Axial wavenumbers for the $TM_{01}$ and $TM_{11}$ modes as a function of frequency in an eccentric coaxial waveguide filled with lossy anisotropic media.

As an additional validation, we consider the extreme case where $\tilde{r}_1 = 10$ mm, $\tilde{r}_0 = 0.25\tilde{r}_1$, and $\tilde{d} = 0.74\tilde{r}_1$, i.e., the waveguide has a large offset that places the inner conductor in close proximity to the outer PEC cylinder because $\tilde{r}_1 - \tilde{r}_0 - \tilde{d} = 0.01\tilde{r}_1$. The waveguide is filled with the same lossy anisotropic media and operates at the same frequency as before. To ensure convergence in this case, we consider expansion order $N = 13$. In Fig. 3.15, we compared the results for the axial wavenumbers with the FIT solution for the $TM_{01}$ and $TM_{11}$ modes. In addition, in Fig. 3.16, we compare the axial electrical field for the same modes, considering a fixed frequency of $f = 10$ GHz. Again, we can note excellent agreement between the two solutions.

It should be noted that according to (3-2), the variation in frequency affects the imaginary part of the complex permittivity. Consequently, the medium that fills the waveguide will change. Given that the geometric param-

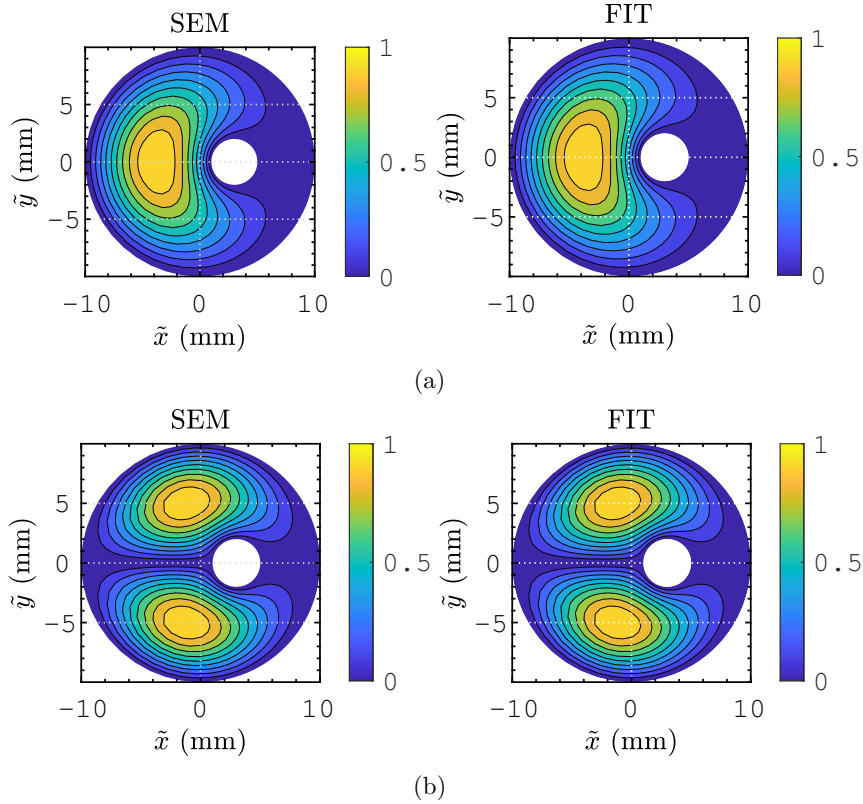

(a)

(b)

Figure 3.14: Normalized axial electric field patterns for the (a) $TM_{01}$ and
(b) $TM_{11}$ calculated using SEM (left) and using FIT (right).

eters of the problem remain unaltered, it suffices to calculate the eigenvalues
associated with this structure only once. Subsequently, the radial wavenumbers
$k_\rho$ are determined through (3-21), and the axial wavenumbers are calculated
by (3-22).

The Table 3.3 presents the total computational time required for cal-
culating the 19 frequency points in Fig. 3.13 using the SEM with expansion
function of order $N = 7$ (Case V), and Fig. 3.15 using expansion order $N = 13$
(Case VI). The time cost was compared with the FIT solution for the same
number of points. The simulations were performed on the same computers and
with the same configurations as those used in Cases III and IV. Once again, the
proposed technique shows an improvement of several orders of magnitude in
the computational performance compared to FIT, including when modeling
eccentric waveguides filled with lossy anisotropic media. As demonstrated
in the preceding analysis, our approach did not exhibit any restriction for
modeling extreme eccentricities and provided accurate solutions in a very short
simulation time.

The low computational cost obtained in the simulations is explained by the small size of the linear system resulting from the formulation using higher-order 2D Lagrange basis functions associated with GLL points. For TM modes are necessary $4N(N-1)$ DoFs and for TE modes are necessary $4N(N+1)$ DoFs. In particular, Case V requires only 168 DoFs, and Case VI requires 624 DoFs.

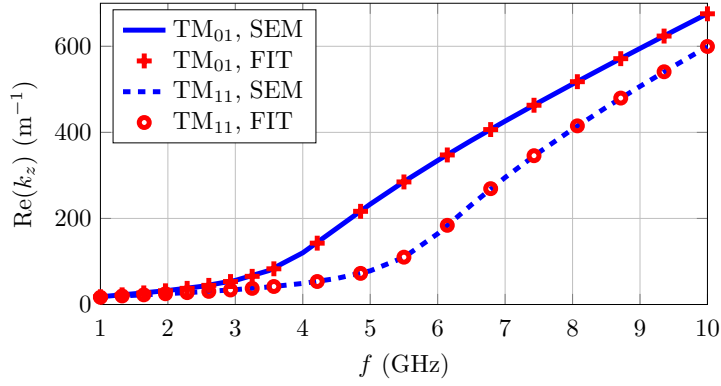

(a)

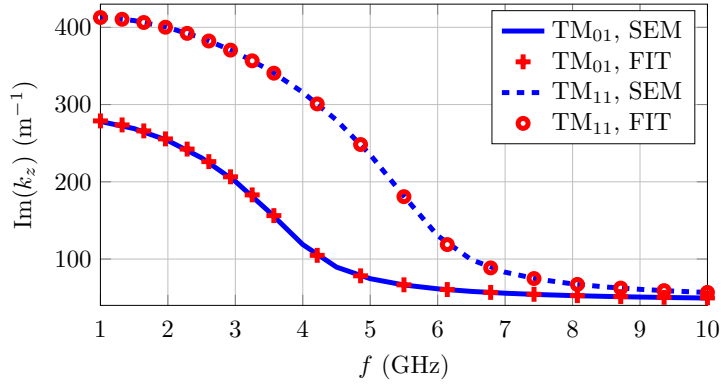

(b)

Figure 3.15: Axial wavenumbers for the $TM_{01}$ and $TM_{11}$ modes as a function of frequency in an eccentric coaxial waveguide filled with lossy anisotropic media with extreme offset.

Table 3.3: Computational cost in the problems with lossy media.

|         | FIT      | SEM     |
|---------|----------|---------|
| Case V  | 11778 s  | 0.66 s  |
| Case VI | 34141 s  | 17.16 s |

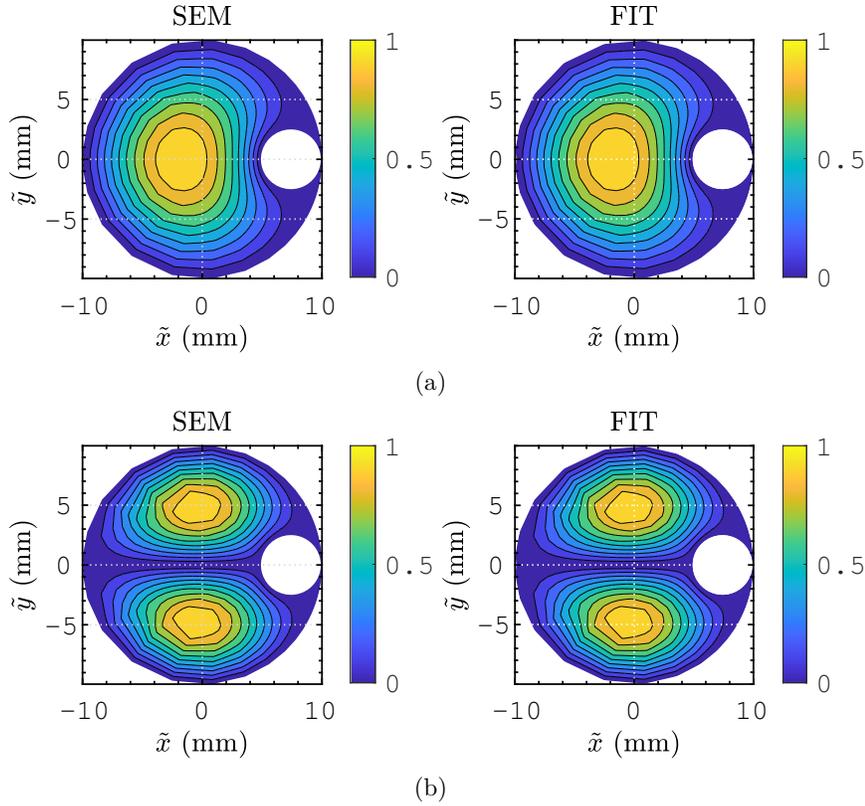

(a)

(b)

Figure 3.16: Normalized axial electric field patterns for the (a) TM$_{01}$ and (b) TM$_{11}$ calculated using SEM (left) and using FIT (right) for an eccentric coaxial waveguide with extreme offset.

## 3.5
## Final Remarks

This work investigated a novel two-dimensional SEM in cylindrical coordinates to rapidly analyze the modal fields in eccentric coaxial waveguides. Given the geometry of the problem, the proposed normalized variational formulation obtains the solution for an entire set of configurations independently of the medium's parameters. The method utilizes transformation optics to map the original eccentric waveguides into equivalent concentric problems. The proposed higher-order expansion for the two-dimensional basis functions is free from numerical errors caused by the Runge effect, and the observed convergence is exponential. Unlike perturbation-based techniques, our approach does not restrict the eccentricity offset and presents great accuracy and time-cost results for extreme geometric cases. The formulation accurately models a wide range of eccentric coaxial waveguides filled with lossy uniaxially anisotropic media. The results demonstrated excellent agreement with a dense-mesh reference FIT solution at a much reduced computational cost. It is worth mentioning that more complex waveguides can be modeled using SEM-based methods. Nonetheless, a vector field formulation is required when dealing with multilayered media and general permittivity and permeability tensors. With a few adaptations, the present SEM can provide a highly efficient tool for studying wave propagation in complex media, with potential applications in the modeling and designing of metamaterial devices and geophysical exploration sensors, for example.

# 4

# Higher-Order SEM to Model Eccentric Anisotropic Two-Layer Waveguides via Conformal Transformation Optics

## 4.1
## Introduction

The wave propagation analysis in non-concentric multilayered cylindrical waveguides is essential for many applications, such as metamaterial devices and geophysical exploration sensors. Recently, various methods have been employed to model circular waveguides partially filled with eccentric rods. In [73], an analytical method using a bipolar coordinate system was employed to obtain the cutoff frequencies for an isotropic medium. In [74], a perturbation method was developed to obtain axial wavenumbers and field patterns considering uniaxially anisotropic media. In [27], a similar problem was solved using transformation optics (TO) to map the original eccentric waveguide into an equivalent concentric problem. In this case, the spectral element method (SEM) in cylindrical coordinates was applied to obtain the solution of the mapped domain. The last approach is not restricted to small eccentricities like the perturbation method [74]. Furthermore, to model electromagnetic logging sensors in multieccentric cylindrically-layered earth formations, a pseudo-analytical method was developed using Graf's addition theorem for the cylindrical waves in [75] and solved using the TO and perturbation-based technique in [10].

The spectral element method has been applied successfully for a wide range of electromagnetic problems [18, 19, 22, 23, 25, 27]. The preliminary study [22] indicates that the SEM-based technique in cylindrical coordinates is more efficient, needing fewer elements and degrees of freedom (DoF), compared to the conventional Cartesian approach [18] for modeling cylindrical-conforming boundaries. Our formulation differs from that in [18, 19] because we represent the fields and media tensors in cylindrical coordinates. In this regard, the current work extends the prior TO studies from [10, 25, 27], which are limited to uniaxially anisotropic tensors, to include non-symmetric and non-Hermitian tensors. Subsequently, we employ the cylindrical SEM to model the mapped problem.

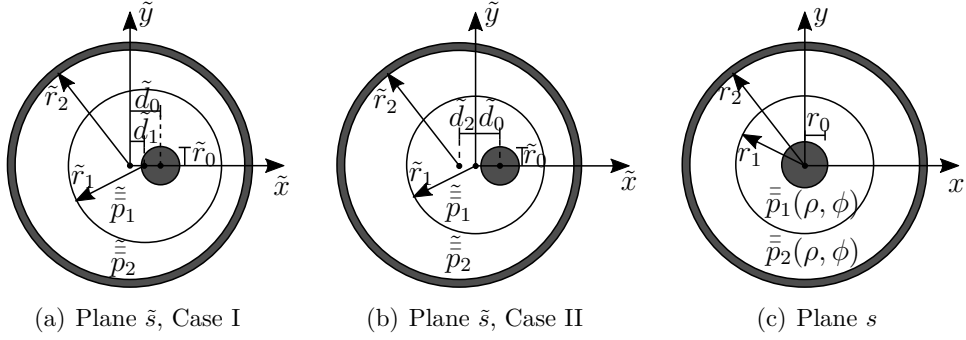

(a) Plane $\tilde{s}$, Case I        (b) Plane $\tilde{s}$, Case II        (c) Plane $s$

Figure 4.1: (a) Geometry of an eccentric two-layer coaxial waveguide considering offsets $\tilde{d}_2 = 0$, $\tilde{d}_0 \neq 0$, and $\tilde{d}_1 \neq 0$. (b) Geometry of an eccentric two-layer coaxial waveguide considering offsets $\tilde{d}_1 = 0$, $\tilde{d}_0 \neq 0$, and $\tilde{d}_2 \neq 0$. (c) Geometry of the mapped concentric two-layer coaxial waveguide.

## 4.2
## Formulation Overview

### 4.2.1
### Variational Waveguide Problem

Consider a non-concentric two-layer waveguide invariant along the axial direction, as depicted in Fig. 4.1(a). The eccentric geometry is referred to in a cylindrical coordinate system denoted by $(\tilde{\rho}, \tilde{\phi}, \tilde{z})$. The geometry is defined in terms of three circular interfaces associated with the radii $\tilde{r}_0$, $\tilde{r}_1$, and $\tilde{r}_2$. Additionally, such cylinders are offset with respect to the $\tilde{z}$-axis by the the distances $\tilde{d}_0$, $\tilde{d}_1$, and $\tilde{d}_2$. Positive offsets are defined as those that move the boundaries toward the $\tilde{x}$-axis. The waveguide is radially bounded by a perfect electric conductor (PEC) at the boundaries of $\tilde{r}_0$ and $\tilde{r}_2$. In Case I, we consider the scenario where the external circular PEC boundary is concentric to the $\tilde{z}$-axis, and the interface between the two layers is displaced $\tilde{d}_1$ from the $\tilde{z}$-axis. In Case II, we investigate the scenario where the interface is concentric, while the external PEC is offset from the $\tilde{z}$-axis by the distance $\tilde{d}_2$. Figs. 4.1(a) and 4.1(b) depict the corresponding cross-section views for Case I and Case II, respectively.

We assume a time-harmonic dependence of the form $e^{-i\omega t}$ and the media are characterized by the real-valued permeability tensor

$$\bar{\bar{\tilde{\mu}}}_j = \mu_0 \bar{\bar{\tilde{\mu}}}_{r,j}, \tag{4-1}$$

with

$$\bar{\bar{\tilde{\mu}}}_{r,j} = \begin{bmatrix} \tilde{\mu}_{\rho\rho,j} & \tilde{\mu}_{\rho\phi,j} & 0 \\ \tilde{\mu}_{\phi\rho,j} & \tilde{\mu}_{\phi\phi,j} & 0 \\ 0 & 0 & \tilde{\mu}_{zz,j} \end{bmatrix}, \tag{4-2}$$

and the complex-valued permittivity tensor

$$\bar{\bar{\tilde{\epsilon}}}_{c,j} = \bar{\bar{\tilde{\epsilon}}}_{r,j} + \frac{i}{\omega\epsilon_0}\bar{\bar{\tilde{\sigma}}}_j, \tag{4-3}$$

with

$$\bar{\bar{\tilde{\epsilon}}}_{r,j} = \begin{bmatrix} \tilde{\epsilon}_{\rho\rho,j} & \tilde{\epsilon}_{\rho\phi,j} & 0 \\ \tilde{\epsilon}_{\phi\rho,j} & \tilde{\epsilon}_{\phi\phi,j} & 0 \\ 0 & 0 & \tilde{\epsilon}_{zz,j} \end{bmatrix}, \bar{\bar{\tilde{\sigma}}}_j = \begin{bmatrix} \tilde{\sigma}_{\rho\rho,j} & \tilde{\sigma}_{\rho\phi,j} & 0 \\ \tilde{\sigma}_{\phi\rho,j} & \tilde{\sigma}_{\phi\phi,j} & 0 \\ 0 & 0 & \tilde{\sigma}_{zz,j} \end{bmatrix}, \tag{4-4}$$

expressed in cylindrical coordinates $(\tilde{\rho}, \tilde{\phi}, \tilde{z})$, where $j = 1, 2$ refer to the inner and outer layers, respectively.

From conformal transformation optics [10, 36, 41], we can relate the eccentric problems with an equivalent concentric problem (represented in Fig. 4.1(c)). The eccentric planes can be written by $\tilde{s} = \tilde{x} + i\tilde{y}$ and the concentric plane by $s = x + iy$. Let $z = \tilde{z}$, $r_2 = \tilde{r}_2$ for Case I, and $r_1 = \tilde{r}_1$ for Case II, we can map the two planes via

$$s = \tilde{x}_2 \frac{\tilde{s} - \tilde{x}_1}{\tilde{s} - \tilde{x}_2}, \tag{4-5}$$

where

$$\tilde{x}_{1,2} = \frac{-\tilde{c} \mp \sqrt{\tilde{c}^2 - 4\tilde{r}_j^2}}{2}, \text{ and } \tilde{c} = \frac{\tilde{r}_0^2 - \tilde{r}_j^2 - \tilde{d}_0^2}{\tilde{d}_0}, \tag{4-6}$$

with $j = 2$ for Case I, and $j = 1$ for Case II. The mapping in (4-5) satisfies Cauchy-Riemann conditions [36]

$$\frac{\partial x}{\partial \tilde{x}} = \frac{\partial y}{\partial \tilde{y}}, \quad \frac{\partial x}{\partial \tilde{y}} = -\frac{\partial y}{\partial \tilde{x}}. \tag{4-7}$$

The original electromagnetic problem in the eccentric coordinates $(\tilde{\rho}, \tilde{\phi}, \tilde{z})$ can be transformed into a concentric problem with coordinates $(\rho, \phi, z)$ using the vector and tensor relations

$$\mathbf{F} = \bar{\bar{J}}_{\text{TO}} \cdot \tilde{\mathbf{F}}, \text{ with } \mathbf{F} \in \{\mathbf{E}, \mathbf{H}\}, \tag{4-8}$$

$$\bar{\bar{p}} = |\bar{\bar{J}}_{\text{TO}}|^{-1}\bar{\bar{J}}_{\text{TO}} \cdot \bar{\bar{\tilde{p}}} \cdot \bar{\bar{J}}_{\text{TO}}^T, \text{ with } p \in \{\mu, \epsilon\}, \tag{4-9}$$

where $\bar{\bar{J}}_{\text{TO}}$ is the Jacobian of the transformation $(\tilde{\rho}, \tilde{\phi}, \tilde{z}) \rightarrow (\rho, \phi, z)$, and $|\bar{\bar{J}}_{\text{TO}}|$ is the respectively determinant. We can express $\bar{\bar{J}}_{\text{TO}}$ by using

$$\bar{\bar{J}}_{\text{TO}} = \bar{\bar{R}}(\phi) \cdot \bar{\bar{J}}_{xyz} \cdot \bar{\bar{R}}^T(\tilde{\phi}), \tag{4-10}$$

in which $\bar{\bar{J}}_{xyz}$ is the Jacobian associated with the transformation $(\tilde{x}, \tilde{y}, \tilde{z}) \rightarrow (x, y, z)$, and $\bar{\bar{R}}(\phi)$ is the transformation matrix for Cartesian-to-cylindrical components. Using the condition (4-7), we obtain

$$\bar{\bar{J}}_{xyz} = \begin{bmatrix} \frac{\partial x}{\partial \tilde{x}} & \frac{\partial x}{\partial \tilde{y}} & 0 \\ -\frac{\partial x}{\partial \tilde{y}} & \frac{\partial x}{\partial \tilde{x}} & 0 \\ 0 & 0 & 1 \end{bmatrix}, \ \ \bar{\bar{R}}(\phi) = \begin{bmatrix} \cos\phi & \sin\phi & 0 \\ -\sin\phi & \cos\phi & 0 \\ 0 & 0 & 1 \end{bmatrix}, \quad (4\text{-}11)$$

where

$$\frac{\partial x}{\partial \tilde{x}} = \left[ (\tilde{x} - \tilde{x}_2)^2 - \tilde{y}^2 \right] \frac{\tilde{x}_2(\tilde{x}_1 - \tilde{x}_2)}{\left[ (\tilde{x} - \tilde{x}_2)^2 + \tilde{y}^2 \right]^2}, \quad (4\text{-}12)$$

$$\frac{\partial x}{\partial \tilde{y}} = 2\tilde{y}(\tilde{x} - \tilde{x}_2) \frac{\tilde{x}_2(\tilde{x}_1 - \tilde{x}_2)}{\left[ (\tilde{x} - \tilde{x}_2)^2 + \tilde{y}^2 \right]^2}. \quad (4\text{-}13)$$

Based on the approach developed in [18, 23, 32, 35], we can define a variational waveguide problem for the concentric cylindrical cross-section $\Omega$ (illustrated in Fig. 4.1(c)). Now, we wish to determine all pairs $(k_z, (\mathbf{e}_s, e_z)) \in (\mathbb{C}, \mathbf{H}_0(\text{curl}, \Omega) \times H_0^1(\Omega))$ that satisfy

$$\begin{cases} \langle \mu_{zz}^{-1} \boldsymbol{\nabla}_s \times \mathbf{e}_s, \boldsymbol{\nabla}_s \times \mathbf{w} \rangle + ik_z \langle \bar{\bar{\mu}}_s^{-1} \boldsymbol{\nabla}_s e_z, \mathbf{w} \rangle - k_0^2 \langle \bar{\bar{\epsilon}}_{cs} \mathbf{e}_s, \mathbf{w} \rangle + k_z^2 \langle \bar{\bar{\mu}}_s^{-1} \mathbf{e}_s, \mathbf{w} \rangle = 0 \\ \langle \bar{\bar{\epsilon}}_{cs} \mathbf{e}_s, \boldsymbol{\nabla}_s w \rangle - ik_z \langle \epsilon_{cz} e_z, w \rangle = 0 \end{cases}$$
$$(4\text{-}14)$$

for all $(\mathbf{w}, w) \in \mathbf{H}_0(\text{curl}, \Omega) \times H_0^1(\Omega)$, where $\mathbf{e}_s = (e_\rho, e_\phi)^T$ and $e_z$ are the transversal and axial electric field components, $k_z$ is the axial wavenumber, and $\boldsymbol{\nabla}_s = \left( \frac{\partial}{\partial \rho}, \frac{1}{\rho} \frac{\partial}{\partial \phi} \right)^T$. The domain is truncated by a PEC, so the boundary conditions to be imposed on the boundary of the domain (at $\rho = r_0$ and $\rho = r_2$) are

$$\begin{aligned} \mathbf{e}_s \times \hat{n} &= \mathbf{0} \\ e_z &= 0 \end{aligned} \quad \text{on } \partial\Omega, \quad (4\text{-}15)$$

where $\hat{n}$ is the unit outward normal on $\partial\Omega$.

## 4.2.2
## Spectral Element Method

The cross-section $\Omega$ of the coaxial computation domain is first discretized into eight curvilinear elements $\Omega^e$, as shown in Fig. 4.2. Each element is described by the polar coordinates $(\rho^e, \phi^e)$, which is then mapped into the reference (Cartesian) coordinates $(\xi, \eta)$ via

$$\rho^e(\xi) = 0.5(\rho_0^e + \rho_1^e) + 0.5(\rho_1^e - \rho_0^e)\xi, \ \ \forall \, \xi \in [-1, 1], \quad (4\text{-}16a)$$

$$\phi^e(\eta) = 0.5(\phi_0^e + \phi_1^e) + 0.5(\phi_1^e - \phi_0^e)\eta, \ \ \forall \, \eta \in [-1, 1]. \quad (4\text{-}16b)$$

The superscript $e$ is used to indicate the parameters of the reference square element in $(\xi, \eta)$. Figs. 4.3 and 4.4 illustrate the considered mapping. The Jacobian matrix associated with this coordinate transformation is given by

$$\bar{\bar{J}}_s^e = \begin{bmatrix} \frac{\partial \rho^e}{\partial \xi} & \frac{\partial \rho^e}{\partial \eta} \\ \rho^e \frac{\partial \phi^e}{\partial \xi} & \rho^e \frac{\partial \phi^e}{\partial \eta} \end{bmatrix} = \begin{bmatrix} 0.5(\rho_1^e - \rho_0^e) & 0 \\ 0 & 0.5\rho^e(\phi_1^e - \phi_0^e) \end{bmatrix}. \qquad (4\text{-}17)$$

The above transformation matrix is diagonal, and as a result, the $\rho$- and $\eta$-directions are decoupled, as well as $\phi$ and $\xi$ directions [22, 23, 25].

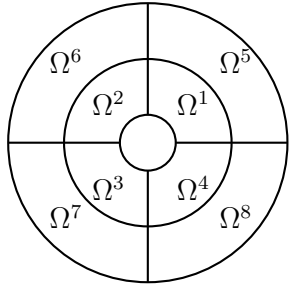

Figure 4.2: Discretization of a coaxial domain $\Omega$ considering eight curvilinear elements.

Our approach employs Lagrange interpolation associated with Gauss-Lobatto-Legendre (GLL) sampling points to represent field components in the reference coordinates. The $N$th-order one-dimensional GLL basis function can be written with

$$\phi_j^N(\xi) = \prod_{k \neq j}^{N+1} \frac{\xi - \xi_k}{\xi_j - \xi_k} \quad \forall \ \xi \in [-1, 1], \qquad (4\text{-}18)$$

where $j = 1, 2, ..., N + 1$, and $N + 1$ interpolation nodes $\xi_k$ correspond to the roots of the $(N + 1)$-degree completed Lobatto polynomial. The two-dimensional GLL basis function is derived through tensor-product expansions

$$\hat{\psi}_{z,ij}(\xi, \eta) = \phi_i^N(\xi)\phi_j^N(\eta) \quad \forall \ (\xi, \eta) \in [-1, 1]^2, \qquad (4\text{-}19)$$

and the mixed-order curl-conforming vector edge-based basis functions are written as

$$\hat{\Phi}_{\xi,ij}(\xi, \eta) = \phi_i^{N-1}(\xi)\phi_j^N(\eta)$$
$$\hat{\Phi}_{\eta,ij}(\xi, \eta) = \phi_i^N(\xi)\phi_j^{N-1}(\eta) \qquad (4\text{-}20)$$

for all $(\xi, \eta) \in [-1, 1]^2$. As an illustrative example, the GLL node distributions in the reference and physical elements considering expansion function order $N = 5$ are depicted in Fig. 4.3. Additionally, the curl-conforming node distributions with the same expansion order are represented in Fig. 4.4. Note that the node distributions are denser towards the edges of the element, a characteristic that effectively mitigates the occurrence of the Runge effect when employing higher-order expansions [28].

To formulate the discrete version of the problem presented in (4-14), we approximate the transversal electrical components in the reference domain $\hat{\mathbf{e}}_s = (\hat{e}_\xi, \hat{e}_\eta)$ using summations of vector edge-based basis functions

$$\hat{e}_\xi(\xi, \eta) = \sum_{j=1}^{N(N+1)} \hat{e}_{\xi,j} \hat{\Phi}_{\xi,j}(\xi, \eta), \tag{4-21a}$$

$$\hat{e}_\eta(\xi, \eta) = \sum_{j=1}^{N(N+1)} \hat{e}_{\eta,j} \hat{\Phi}_{\eta,j}(\xi, \eta), \tag{4-21b}$$

and the axial field $\hat{e}_z$ using summations of GLL basis functions

$$\hat{e}_z(\xi, \eta) = \sum_{j=1}^{(N+1)^2} \hat{e}_{z,j} \hat{\psi}_{z,j}(\xi, \eta). \tag{4-21c}$$

By employing the expansions (4-21) and applying the Galerkin method for test functions in the variational formulation (4-14), we obtain the following generalized eigenvalue problem

$$\left( -\mathbf{S}_{ss} - \mathbf{K}_{sz} \mathbf{M}_{zz}^{-1} \mathbf{K}_{zs} + k_0^2 \mathbf{M}_{ss1} \right) \hat{\mathbf{e}}_{s,m} = k_z^2 \mathbf{M}_{ss2} \hat{\mathbf{e}}_{s,m}, \tag{4-22}$$

where $k_z^2$ are the eigenvalues associated with the axial wavenumbers, and $\hat{\mathbf{e}}_{s,m}$ are the eigenvectors associated with the transversal electric field components. The matrices used in (4-22) have the same form as those found in [18, 23].

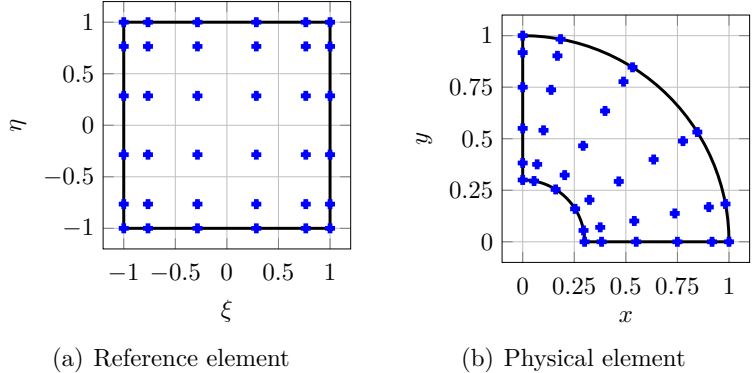

(a) Reference element    (b) Physical element

Figure 4.3: GLL node distributions using expansion order $N = 5$.

It is crucial to emphasize that the constitutive media parameters in the mapped concentric domain exhibit variations in the $\rho$ and $\phi$ directions. This is in contrast with the constant tensor expressions presented in (4-1)–(4-4).

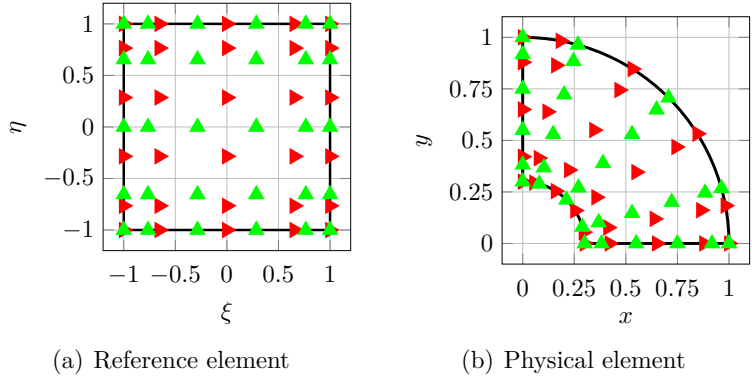

(a) Reference element          (b) Physical element

Figure 4.4: Curl-conforming node distributions using expansion order $N = 5$.

## 4.3
## Numerical Results

As our first validation, we consider the geometry for Case I, characterized by the radii $\tilde{r}_0 = 1$ mm, $\tilde{r}_1 = 5.36$ mm, $\tilde{r}_2 = 10$ mm, and offsets $\tilde{d}_0 = 2$ mm, $\tilde{d}_1 = 1.42$ mm, $\tilde{d}_2 = 0$ mm. The media are described by the relative permittivity and electrical conductivity diagonal tensors represented in cylindrical coordinates by

$$\bar{\bar{\epsilon}}_{r,1} = \begin{bmatrix} 5.7 & 0 & 0 \\ 0 & 2.3 & 0 \\ 0 & 0 & 4.1 \end{bmatrix}, \quad \bar{\bar{\sigma}}_1 = \begin{bmatrix} 0.38 & 0 & 0 \\ 0 & 0.25 & 0 \\ 0 & 0 & 0.30 \end{bmatrix}, \tag{4-23a}$$

$$\bar{\bar{\epsilon}}_{r,2} = \begin{bmatrix} 4.3 & 0 & 0 \\ 0 & 2.7 & 0 \\ 0 & 0 & 3.9 \end{bmatrix}, \quad \bar{\bar{\sigma}}_2 = \begin{bmatrix} 0.14 & 0 & 0 \\ 0 & 0.21 & 0 \\ 0 & 0 & 0.23 \end{bmatrix}. \tag{4-23b}$$

The permeability parameter is that of the vacuum. Fig. 4.5 shows the real and imaginary part of axial wavenumber $k_z$ for the first two hybrid modes as a function of the operating frequency $f$ varying from 0.5 GHz to 6.5 GHz. The results obtained with the proposed SEM show excellent agreement with those from finite-element method (FEM) from COMSOL Multiphysics [2].

With a second example, we investigate the accuracy of the SEM solution for Case II. The geometry parameters are defined by $\tilde{r}_0 = 1$ mm, $\tilde{r}_1 = 5$ mm, $\tilde{r}_2 = 9.96$ mm, $\tilde{d}_0 = 1$ mm, $\tilde{d}_1 = 0$ mm, and $\tilde{d}_2 = -2.68$ mm. We consider the waveguide filled with the same media as before. Fig. 4.6 shows the results for the first two hybrid modes as a function of frequency. Again, we obtained great agreement between the SEM and FEM solutions.

Now, we will explore the capabilities of the present SEM formulation to model waveguides that operate in low-frequency and are filled with more complex anisotropic materials, where the media are non-reciprocal and described by non-symmetric and non-Hermitian tensors. Consider the scenario of Case II defined by the radii $\tilde{r}_0 = 4$ in, $\tilde{r}_1 = 10$ in, $\tilde{r}_2 = 16.32$ in, and the eccentricity offsets $\tilde{d}_0 = 1$ in, $\tilde{d}_1 = 0$ in, $\tilde{d}_2 = -1.91$ in. The relative permittivity and electrical conductivity tensors are given by

$$\bar{\bar{\epsilon}}_{r,1} = \begin{bmatrix} 1.9 & 2.3 & 0 \\ 2.7 & 1.8 & 0 \\ 0 & 0 & 2.8 \end{bmatrix}, \quad \bar{\bar{\sigma}}_1 = 10^{-4} \begin{bmatrix} 5 & 4 & 0 \\ 3 & 2 & 0 \\ 0 & 0 & 6 \end{bmatrix}, \tag{4-24a}$$

$$\bar{\bar{\epsilon}}_{r,2} = \begin{bmatrix} 3.7 & 4.3 & 0 \\ 4.2 & 3.0 & 0 \\ 0 & 0 & 3.2 \end{bmatrix}, \quad \bar{\bar{\sigma}}_2 = 10^{-2} \begin{bmatrix} 4 & 5 & 0 \\ 2 & 3 & 0 \\ 0 & 0 & 5 \end{bmatrix}. \tag{4-24b}$$

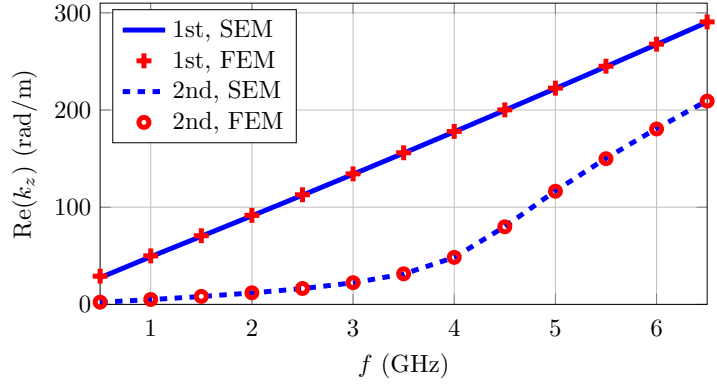

(a) Phase constant

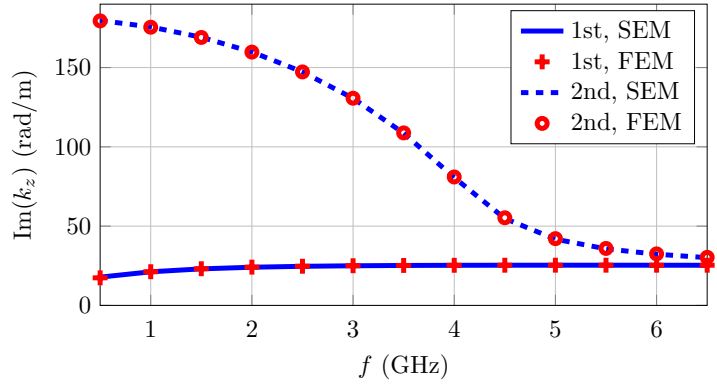

(b) Attenuation constant

Figure 4.5: Axial wavenumbers for the first two modes in an eccentric two-layer waveguide considering $\tilde{d}_2 = 0$.

The permeability is once again the same as that of the vacuum. Assuming operating frequency $f \in \{100, 150, 200\}$ MHz, the results for the axial wavenumbers $k_z$ from our formulation are compared with FEM results [2] in the complex region $[0, 10]^2$. Excellent agreement is observed in Fig. 4.7. For the case $f = 100$ MHz, the FEM solution presented one spurious mode, which can be seen in Fig. 4.7(a). Such eigenvalue is located near the least attenuated mode, and its non-physical characteristics were confirmed by inspecting the field distribution, which exhibits a highly irregular spatial distribution.

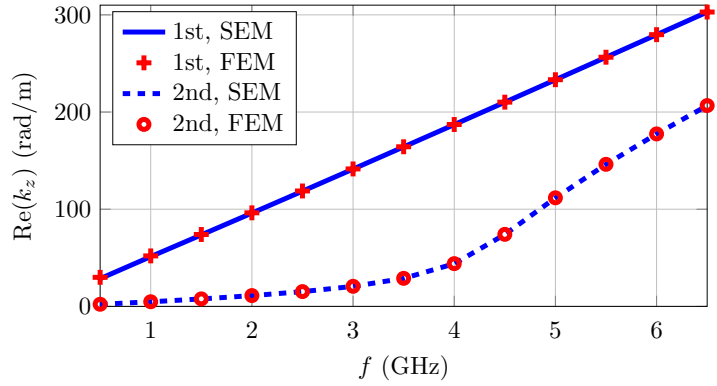

(a) Phase constant

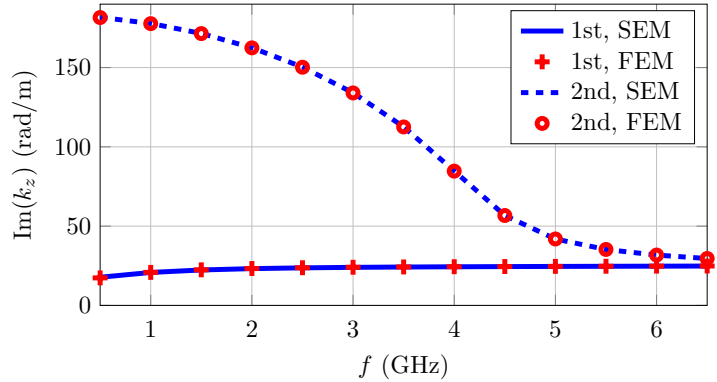

(b) Attenuation constant

Figure 4.6: Axial wavenumbers for the first two modes in an eccentric two-layer waveguide considering $\tilde{d}_1 = 0$.

For the three examples explored here, the SEM computation domains were each discretized into 8 elements with expansion functions of order 10, which implies 1444 DoF. In the COMSOL's FEM model, we employed a cubic-order discretization with the extremely fine physics-controlled mesh option to obtain accurate reference solutions. This led to 35982 DoF for the first example, 34608 DoF for the second, and 42606 DoF for the last example. The current SEM formulation requires roughly 4 % of the DoFs needed by the FEM for the first two examples and 3 % for the third.

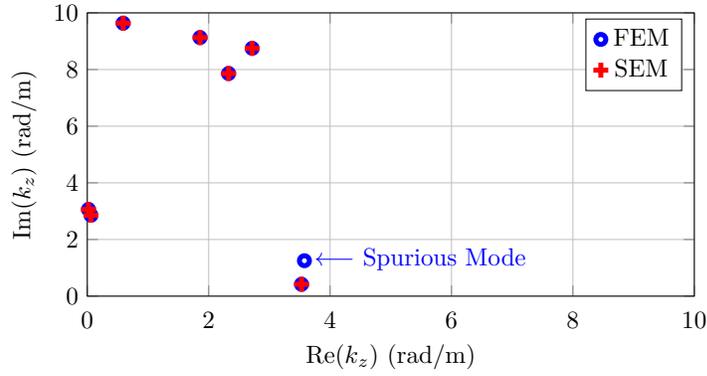

(a) $f = 100$ MHz

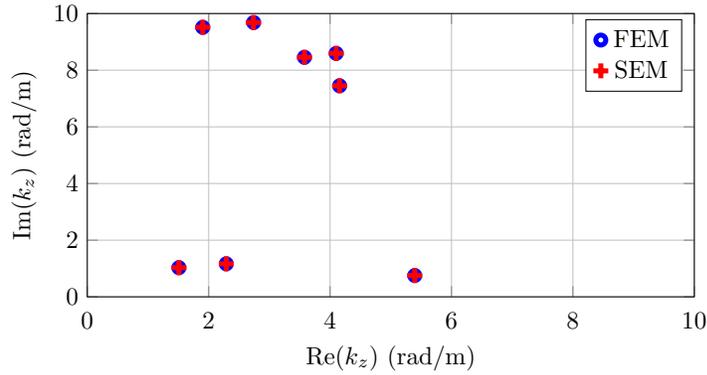

(b) $f = 150$ MHz

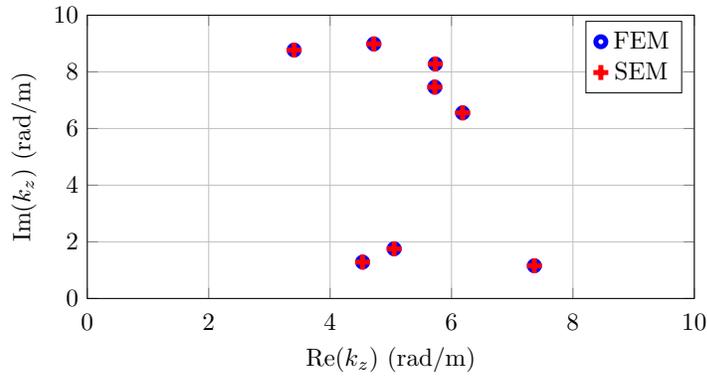

(c) $f = 200$ MHz

Figure 4.7: Axial wavenumbers obtained over the complex region $[0, 10]^2$ for non-concentric two-layer waveguides filled with lossy and anisotropic media.

## 4.4
## Final Remarks

This study presents an improved spectral element method tailored for analyzing eccentric two-layer waveguides with anisotropic non-reciprocal media. Through the application of conformal transformation optics, the complexity of the eccentric waveguide was transformed into a more tractable concentric counterpart. The cylindrical SEM was applied in the concentric domain and proved accurate to model waveguides filled with complex media (lossy, inhomogeneous, and anisotropic) for a wide frequency range. Our approach required far fewer DoF in comparison to the conventional FEM solution. In future works, we plan to incorporate a perfect matching layer into our current SEM model. This will enable the simulation of open radial boundaries and facilitate the modeling of logging while drilling sensors in eccentric boreholes with invasion zones for geophysical exploration.

# 5

# Electromagnetic Analysis of Logging Sensors for Oil Well in Geophysical Formations via SEM

## 5.1
## Introduction

Prospection of hydrocarbon reservoirs in deep subterranean formations is a task of fundamental importance to oil and gas exploration. Assessing these reservoirs often entails utilizing well-logging techniques, wherein sensors are deployed within drilled boreholes to gather data about the geological properties of the surrounding formations. These properties include porosity, which can be determined through the implementation of acoustic logging tools, and electrical resistivity (or conductivity), which can be assessed via electromagnetic logging tools. The latter, known as logging-while-drilling (LWD) sensors, typically consist of coil antennas encased within a metallic mandrel positioned within the borehole (see Fig. 5.1). Operating within a frequency range from 100 kHz to approximately 1 GHz, with a common frequency of around 2 MHz for commercially available sensors, these tools capture crucial information about the subsurface formations [76].

Typical oilfield well-logging scenarios include a series of fractured soil formations in which the electromagnetic constitutive parameters (permittivity and permeability) often exhibit uniformity within specific layers. However, these formations are anisotropic, potentially manifesting diverse electrical resistivities or conductivities across different orientations. The conductivity of the geological media correlates directly with the presence of water or oil within the formation. By discerning the resistivity of these strata, it becomes feasible to ascertain the saturation of hydrocarbons within the formation and its relationship with soil porosity.

In conventional drilling operations, the drill collar imparts rotational motion and downward force to the drill bit, which pulverizes and fragments rock formations. To facilitate this process, the borehole is flooded with drilling mud to cool and lubricate the bit. This drilling mud may possess varying conductive properties, depending on whether it is water-based or oil-based, and can infiltrate the surrounding formation, particularly in regions of high poros-ity adjacent to the borehole. Consequently, an invasion zone is established, characterized by intermediate electrical properties between the borehole mud and the native formation materials [77–79].

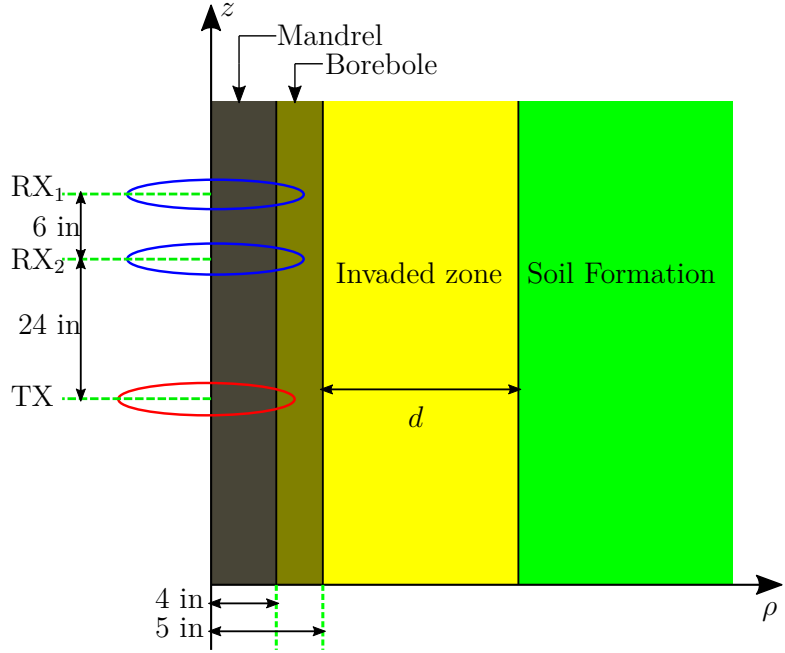

Figure 5.1: The geometry of a well-logging tool in a geophysical formation with invaded zone.

In this work, we developed a computational tool for modeling electromagnetic well-logging sensors via the fields obtained by the recently introduced higher-order SEM formulation discussed in [22–24]. Our approach differs from [21, 66] for using the SEM in cylindrical coordinates, a higher-order expansion for the two-dimensional basis function, and horizontal coil antennas instead of dipoles.

## 5.2
## Variational Waveguide Problem

Consider a cylindrical waveguide structure invariant in the axial direction $z$ and excited by a generic source at the position $z = z_T$. This source generates the fields $\mathbf{E}^+$ and $\mathbf{H}^+$ traveling in the positive $z$-direction and the fields $\mathbf{E}^-$ and $\mathbf{H}^-$ traveling in the negative $z$-direction. The electromagnetic fields can be expressed in terms of the waveguide modes as

$$\mathbf{E}^\pm = \sum_{p=1}^\infty A_p^\pm [\mathbf{e}_{s,p}(\rho,\phi) \pm \hat{z} e_{z,p}(\rho,\phi)] e^{\pm i k_{z,p}(z-z_T)} \tag{5-1a}$$

$$\mathbf{H}^\pm = \sum_{p=1}^\infty A_p^\pm [\pm \mathbf{h}_{s,p}(\rho,\phi) + \hat{z} h_{z,p}(\rho,\phi)] e^{\pm i k_{z,p}(z-z_T)}, \tag{5-1b}$$

for $z \gtrless z_T$, where $\mathbf{e}_{s,p}$, $\mathbf{h}_{s,p}$ and $e_{z,p}$, $h_{z,p}$ are transverse and axial fields, respectively, associated with the axial wavenumber $k_{z,p} \in \mathbb{C}$, and $A_p^{\pm}$ are the complex amplitudes determined by the source. In our discrete (computational) domain, these fields and wavenumbers will be translated into eigenvectors and the associated eigenvalues.

On the grounds of the approach developed in [18,23,32,35] and assuming invariance along $z$ and source-free region, a variational formulation can be obtained by rewriting Maxwell's equations in complex anisotropic media as a system of equations for the electric field components transverse and parallel to $z$, $\mathbf{e}_s$ and $e_z$ respectively, defined on the cylindrical waveguide cross-section $\Omega \subset \mathbb{R}^2$. This essentially reduces the 3D waveguide problem to a 2D problem. We wish to determine all pairs $(k_z, (\mathbf{e}_s, e_z)) \in (\mathbb{C}, \mathbf{H}_0(\mathrm{curl}, \Omega) \times H_0^1(\Omega))$ that satisfy

$$\begin{cases} \langle \mu_{zz}^{-1} \boldsymbol{\nabla}_s \times \mathbf{e}_s, \boldsymbol{\nabla}_s \times \mathbf{w} \rangle + ik_z \langle \bar{\bar{\mu}}_s^{-1} \boldsymbol{\nabla}_s e_z, \mathbf{w} \rangle \\ \qquad - k_0^2 \langle \bar{\bar{\epsilon}}_{cs} \mathbf{e}_s, \mathbf{w} \rangle + k_z^2 \langle \bar{\bar{\mu}}_s^{-1} \mathbf{e}_s, \mathbf{w} \rangle = 0 \\ \qquad \langle \bar{\bar{\epsilon}}_{cs} \mathbf{e}_s, \boldsymbol{\nabla}_s w \rangle - ik_z \langle \epsilon_{cz} e_z, w \rangle = 0 \end{cases} \tag{5-2}$$

for all $(\mathbf{w}, w) \in \mathbf{H}_0(\mathrm{curl}, \Omega) \times H_0^1(\Omega)$, where $\mathbf{e}_s = (e_\rho, e_\phi)^T$ and $\boldsymbol{\nabla}_s = \left( \frac{\partial}{\partial \rho}, \frac{1}{\rho} \frac{\partial}{\partial \phi} \right)^T$. The domain is truncated by a perfect electric conductor (PEC).

The variational waveguide problem for electrical field components (5-2) can be solved by the recently introduced SEM technique discussed in [22–24]. Considering anisotropic media represented by diagonal tensors $p = \mathrm{diag}\{p_\rho, p_\phi, p_z\}$, we can get the magnetic field components in cylindrical coordinates via

$$h_\rho = \frac{1}{ik_z} \left[ \frac{\partial h_z(\rho, \phi)}{\partial \rho} - i\omega\epsilon_0 \epsilon_{c\phi} e_\phi(\rho, \phi) \right] \tag{5-3}$$

$$h_\phi = \frac{1}{ik_z} \left[ \frac{1}{\rho} \frac{\partial h_z(\rho, \phi)}{\partial \phi} + i\omega\epsilon_0 \epsilon_{c\rho} e_\rho(\rho, \phi) \right] \tag{5-4}$$

$$h_z = \frac{1}{i\omega\mu_0 \mu_{zz}} \left[ \frac{1}{\rho} \frac{\partial (\rho e_\phi(\rho, \phi))}{\partial \rho} - \frac{1}{\rho} \frac{\partial e_\rho(\rho, \phi)}{\partial \phi} \right], \tag{5-5}$$

where the constitutive media parameters can be inhomogeneous. The derivative in (5-3)-(5-5) can be solved via finite-difference [80].

## 5.3
## Electromagnetic Well-Logging Sensors

### 5.3.1
### Mode Excitation from Sources

We can determine the unknown amplitudes $A_p^{\pm}$ in (5-1) by using the Lorentz reciprocity theorem [81, Ch. 1]. Considering the two sets of sources $\mathbf{J}_1$, $\mathbf{M}_1$ and $\mathbf{J}_2$, $\mathbf{M}_2$, which generate the fields $\mathbf{E}_1$, $\mathbf{H}_1$ and $\mathbf{E}_2$, $\mathbf{H}_2$, respectively, in the cylindrical volume $V$ enclosed by the closed surface $S$, we can derive

$$\oint_S \left(\mathbf{E}_1 \times \mathbf{H}_2 - \mathbf{E}_2 \times \mathbf{H}_1\right) \cdot d\mathbf{s} = \int_V \left(\mathbf{E}_2 \cdot \mathbf{J}_1 - \mathbf{E}_1 \cdot \mathbf{J}_2 + \mathbf{H}_1 \cdot \mathbf{M}_2 - \mathbf{H}_2 \cdot \mathbf{M}_1\right) dv, \tag{5-6}$$

where $d\mathbf{s}$ is the differential outward normal on the surface $S$, and $dv$ is the cylindrical differential volume in $V$.

For a radially bounded waveguide terminated by a perfectly electric conductor, the portion of the surface integral over the waveguide walls vanishes because the tangential electric fields are zero, i.e., $\mathbf{E} \times \mathbf{H} \cdot \hat{\rho} = \mathbf{H} \cdot (\hat{\rho} \times \mathbf{E}) = 0$. Consequently, this reduces the integration (5-6) to be just over the top and bottom guide cross-section $S_0^{\pm}$ as

$$\int_{S_0^{\pm}} \left(\mathbf{E}_1 \times \mathbf{H}_2 - \mathbf{E}_2 \times \mathbf{H}_1\right) \cdot \hat{z}\, ds =$$
$$\int_V \left(\mathbf{E}_2 \cdot \mathbf{J}_1 - \mathbf{E}_1 \cdot \mathbf{J}_2 + \mathbf{H}_1 \cdot \mathbf{M}_2 - \mathbf{H}_2 \cdot \mathbf{M}_1\right) dv. \tag{5-7}$$

Furthermore, for reciprocal media [82, Ch. 15], the modes are orthogonal over the guide cross-section in the sense of

$$\int_{S_0^{\pm}} \mathbf{E}_p^{\pm} \times \mathbf{H}_{p'}^{\pm} \cdot \hat{z}\, ds = \int_{S_0^{\pm}} \mathbf{E}_p^{\pm} \times \mathbf{H}_{p'}^{\mp} \cdot \hat{z}\, ds = 0, \quad \text{for all} \quad p \neq p'. \tag{5-8}$$

We can find the complex amplitudes $A_p^{\pm}$ by selecting $\mathbf{E}_1$ and $\mathbf{H}_1$ as the fields radiated by the sources. In other words, $\mathbf{E}_1 = \mathbf{E}^{\pm}$ and $\mathbf{H}_1 = \mathbf{H}^{\pm}$, depending on propagation direction. Firstly, let $\mathbf{E}_2$ and $\mathbf{H}_2$ be the waveguide mode traveling in the negative $z$-direction, i.e.,

$$\mathbf{E}_2 = \mathbf{E}_{p'}^{-}, \tag{5-9}$$
$$\mathbf{H}_2 = \mathbf{H}_{p'}^{-}.$$

Substituting (5-9) into (5-7) and using $\mathbf{J}_1 = \mathbf{J}$, $\mathbf{M}_1 = \mathbf{M}$, and $\mathbf{J}_2 = \mathbf{0}$, $\mathbf{M}_2 = \mathbf{0}$, we get

$$\int_{S_0^+} A_p^+ \left(\mathbf{E}_p^+ \times \mathbf{H}_{p'}^- - \mathbf{E}_{p'}^- \times \mathbf{H}_p^+\right) \cdot \hat{z}\, ds$$
$$- \int_{S_0^-} A_p^- \left(\mathbf{E}_p^- \times \mathbf{H}_{p'}^- - \mathbf{E}_{p'}^- \times \mathbf{H}_p^-\right) \cdot \hat{z}\, ds = \int_V \left(\mathbf{E}_{p'}^- \cdot \mathbf{J} - \mathbf{H}_{p'}^- \cdot \mathbf{M}\right) dv. \tag{5-10}$$

Under the orthogonal condition (5-8) the integrands in (5-10) are only different from zero for $p = p'$, and using (5-1) we can rewrite (5-10) as

$$\int_{S_0^+} A_p^+ \left[ e_{\rho,p} h_{\phi,p} - e_{\phi,p} h_{\rho,p} - (e_{\rho,p} h_{\phi,p} - e_{\phi,p} h_{\rho,p}) \right] ds = \int_V \left( \mathbf{E}_p^- \cdot \mathbf{J} - \mathbf{H}_p^- \cdot \mathbf{M} \right) dv. \tag{5-11}$$

Then, we can obtain $A_p^+$ for all $p \in \mathbb{N}$ via

$$A_p^+ = \frac{S_p^+}{N_p^+}, \tag{5-12}$$

where

$$S_p^+ = \int_V \left[ \left( e_{\rho,p}\hat{\rho} + e_{\phi,p}\hat{\phi} - e_{z,p}\hat{z} \right) \cdot \mathbf{J} \right.$$
$$\left. + \left( h_{\rho,p}\hat{\rho} + h_{\phi,p}\hat{\phi} - h_{z,p}\hat{z} \right) \cdot \mathbf{M} \right] e^{-ik_{z,p}z} \, dv, \tag{5-13}$$

and

$$N_p^+ = -2 \int_{S_0^+} \left( e_{\rho,p} h_{\phi,p} - e_{\phi,p} h_{\rho,p} \right) ds. \tag{5-14}$$

Similarly, to get $A_p^-$, we can rewrite (5-7) considering $\mathbf{E}_2 = \mathbf{E}_{p'}^+$ and $\mathbf{H}_2 = \mathbf{H}_{p'}^+$ as

$$\int_{S_0^+} A_p^+ \left( \mathbf{E}_p^+ \times \mathbf{H}_{p'}^+ - \mathbf{E}_{p'}^+ \times \mathbf{H}_p^+ \right) \cdot \hat{z} \, ds$$
$$- \int_{S_0^-} A_p^- \left( \mathbf{E}_p^- \times \mathbf{H}_{p'}^+ - \mathbf{E}_{p'}^+ \times \mathbf{H}_p^- \right) \cdot \hat{z} \, ds = \int_V \left( \mathbf{E}_{p'}^+ \cdot \mathbf{J} - \mathbf{H}_{p'}^+ \cdot \mathbf{M} \right) dv. \tag{5-15}$$

Using the orthogonal condition (5-8) and the fields (5-1), we can obtain $A_p^-$ for all $p \in \mathbb{N}$ via

$$A_p^- = \frac{S_p^-}{N_p^-}, \tag{5-16}$$

where

$$S_p^- = \int_V \left[ \left( e_{\rho,p}\hat{\rho} + e_{\phi,p}\hat{\phi} + e_{z,p}\hat{z} \right) \cdot \mathbf{J} \right.$$
$$\left. - \left( h_{\rho,p}\hat{\rho} + h_{\phi,p}\hat{\phi} + h_{z,p}\hat{z} \right) \cdot \mathbf{M} \right] e^{+ik_{z,p}z} \, dv, \tag{5-17}$$

and

$$N_p^- = -2 \int_{S_0^-} \left( e_{\rho,p} h_{\phi,p} - e_{\phi,p} h_{\rho,p} \right) ds. \tag{5-18}$$

In summary, the complex amplitudes associated with the $p$th waveguide modes in (5-1) are given by

$$A_p^\pm = \frac{S_p^\pm}{N_p},\tag{5-19}$$

for all $p \in \mathbb{N}$, where the source excitation is defined by

$$S_p^\pm = \int_V \left[(e_{\rho,p}\hat{\rho} + e_{\phi,p}\hat{\phi} \mp e_{z,p}\hat{z}) \cdot \mathbf{J} \pm (h_{\rho,p}\hat{\rho} + h_{\phi,p}\hat{\phi} \mp h_{z,p}\hat{z}) \cdot \mathbf{M})\right] e^{\mp ik_{z,p}z}dv,\tag{5-20}$$

and normalization term is $N_p = N_p^+ = N_p^-$ for $S_0^+ = S_0^-$.

### 5.3.2
### Electric Coil Antenna Source

The operating wavelengths in geophysical applications are typically much larger than the wire antenna thickness. Hence, we can model the antennas as filament-like current source densities. We consider the excitation of a horizontal electric coil antenna placed at $z = z_T$ and defined by

$$\mathbf{J} = I_T \delta(\rho - \rho_T)\delta(z - z_T)\hat{\phi},\tag{5-21}$$

where $I_T$ is the antenna constant current. Considering $\mathbf{M} = \mathbf{0}$, the propagation in the positive $z$-direction and replacing (5-21) in (5-20), we get

$$S_p^+ = e^{-ik_{z,p}z_T} \int_0^{2\pi} e_{\rho,p}(\rho_T,\phi)\, \rho_T d\phi,\tag{5-22}$$

and the normalization term is given by

$$N_p = -2 \int_{S_0} \left(e_{\rho,p}(\rho,\phi)h_{\phi,p}(\rho,\phi) - e_{\phi,p}(\rho,\phi)h_{\rho,p}(\rho,\phi)\right) \rho d\rho d\phi.\tag{5-23}$$

### 5.3.3
### Received Voltage

The well-logging tools employ horizontal and tilted coil antennas as receiving sensors, and the induced voltage (electromotive force) at receptors RXs is the main output parameter of interest for geophysical prospecting. In this work, horizontal coil antennas are considered for the receiving sensors. The voltage induced can be obtained by integrating the total electric field along the coil path in $z = z_{RX}$. The received voltage can be written as

$$V_{RX} = -\sum_{p=1}^{\infty} \left[ A_p^+ e^{+ik_{z,p}(z_{RX}-z_T)} \int_0^{2\pi} e_{\phi,p}(\rho_{RX}, \phi) \; \rho_{RX} d\phi \right], \qquad (5\text{-}24)$$

where $\rho_{RX}$ is the radial position for the coil antenna.

## 5.4
## Numerical Results

To illustrate the application of the technique developed here, we consider a logging-well tool consisting of one transmitter and two receiver horizontal coil antennas wrapped around a perfect electric conducting mandrel inside a vertical borehole represented in Fig. 5.1. In all validation cases, the axial positions of transmitter TX and receivers RX$_2$ and RX$_1$ are $z_T$, $z_T + 24$ in, and $z_T + 30$ in, respectively. The concentric mandrel, coil antennas, and borehole radii are 4 in, 4.5 in, and 5 in, respectively. The antenna is excited with $I_T = 1$ A and operates at 2 MHz. Given the low operational frequency, we may approximate the remaining constitutive parameters to those of a vacuum medium. In the SEM formulation, we consider the computational domain partitioned into 20 elements and expansion order for the basis function $N = 8$. The computational domain was truncated at 45 in by PEC. In geophysical prospecting applications, the output parameters of interest are the amplitude ratio (AR) voltage and phase difference (PD) voltage between the two receiver antennas.

As a first validation, we consider the case where the borehole is filled with mud $\sigma_{\text{mud}} = 10$ S/m, or oil $\sigma_{\text{oil}} = 0.0005$ S/m. The surrounding soil formation varies from $\sigma_{\text{soil}} = 1$ S/m to $\sigma_{\text{soil}} = 6$ S/m. Figs. 5.2 and 5.3 present results for AR and PD, demonstrating excellent agreement between the present SEM and the axial mode-matching (AMM) solution from [15].

As a second scenario, we consider the case with the oil-invaded zone of $d = 10$ in, represented in Fig. 5.1. Initially, we consider a homogeneous invaded zone with $\sigma_{\text{zone}} = 0.5$ S/m; the borehole is filled with the same oil conductive as before. Fig. 5.4 shows good agreement between the SEM results and the AMM from [15]. Subsequently, with an example case, we consider a similar problem, but now the oil-invaded zone is inhomogeneous and described by the continuous function

$$\sigma_{\text{zone}} = A\rho + B, \qquad (5\text{-}25)$$

where

$$A = \frac{\sigma_2 - \sigma_1}{\rho_2 - \rho_1}, \text{ and } B = \sigma_2 - A\rho_2, \qquad (5\text{-}26)$$

and the other parameters are $\sigma_1 = 0.0005$ S/m, $\sigma_2 = 0.5$ S/m, $\rho_1 = 5$ in,

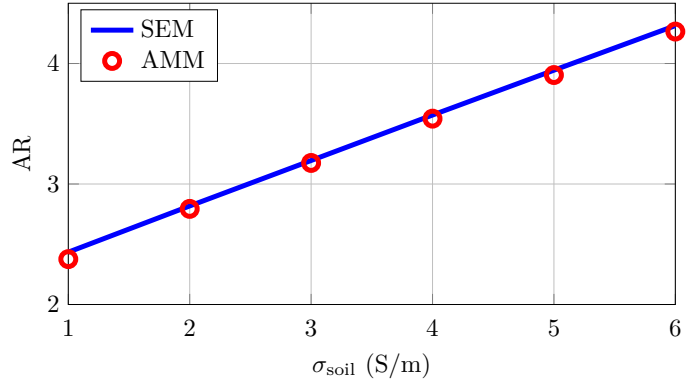

(a) Amplitude ratio

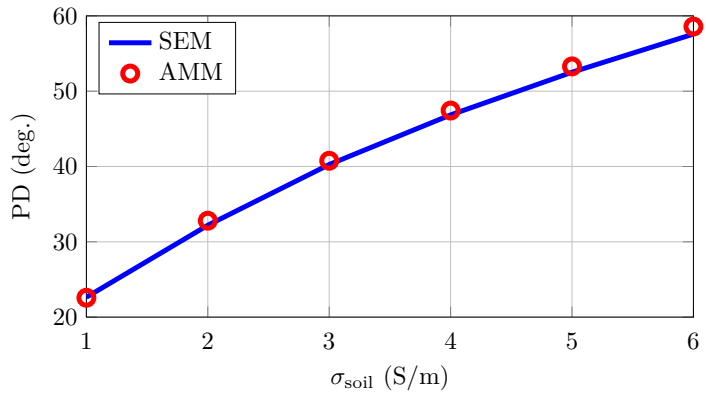

(b) Phase difference

Figure 5.2: Voltage AR and PD for the logging-well tool where the borehole is filled with mud as a function of conductive soil formation.

$\rho_2 = 15$ in. The results for this case are shown in Fig. 5.4, but with no validation for comparison available, as the AMM is unable to handle inhomogeneous radial layers.

As a last validation case, we consider the example previously presented in [83]. Consider the configuration depicted in Fig. 5.1, where the borehole is filled with oil $\sigma_{\text{oil}} = 0.0005$ S/m, the soil formation is invariant with $\sigma_{\text{soil}} = 5$ S/m, and the dimension of oil-invaded zone varies from $d = 0$ to $d = 10$ in and is filled with $\sigma_{\text{zone}} = 0.01$ S/m. Fig. 5.5 shows good agreement between the SEM solution and the axial mode-matching (AMM) and radial mode-matching (RMM). Relative errors observed in comparison with reference methods can be attributed to either an incomplete emulation of the radial-unbounded domain or to the mesh employed for representing the components of the electric and magnetic fields within the mode source formulation.

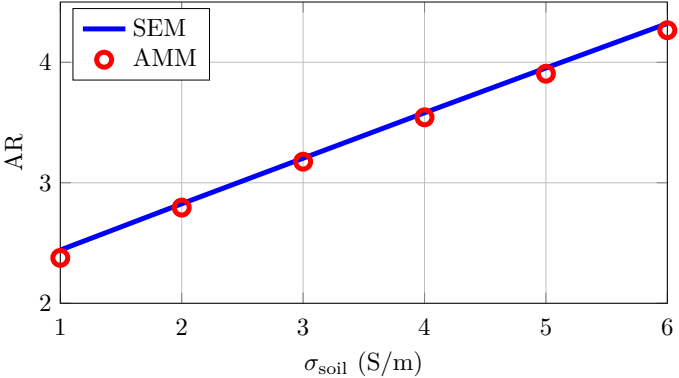

(a) Amplitude ratio

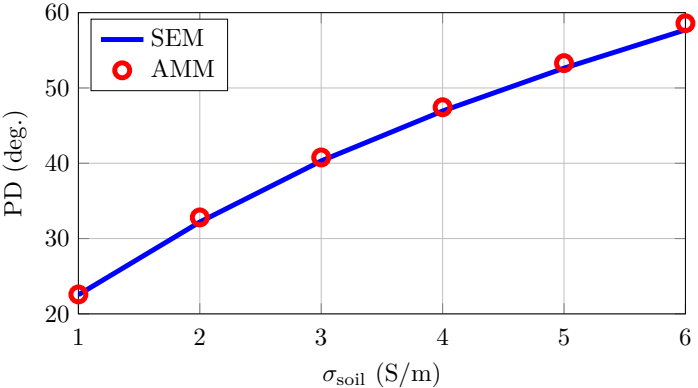

(b) Phase difference

Figure 5.3: Voltage AR and PD for the logging-well tool where the borehole is filled with oil as a function of conductive soil formation.

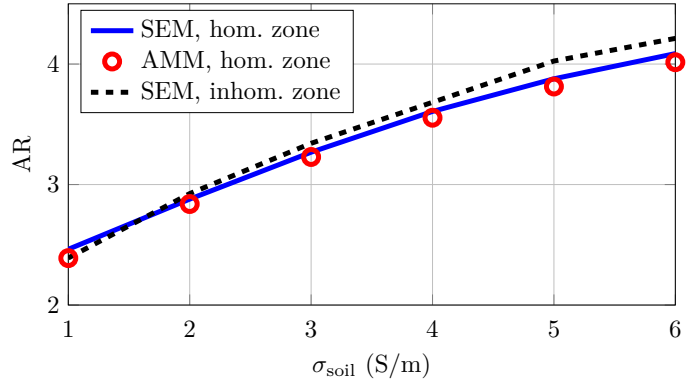

(a) Amplitude ratio

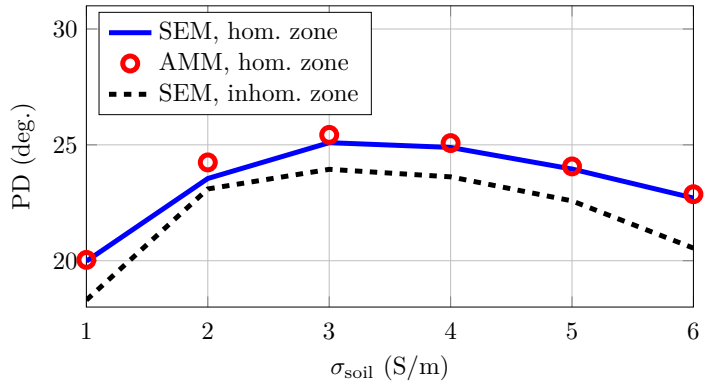

(b) Phase difference

Figure 5.4: Voltage AR and PD for the logging-well tool where the borehole is filled with oil, and there are homogeneous and inhomogeneous oil-invaded zones as a function of conductive soil formation.

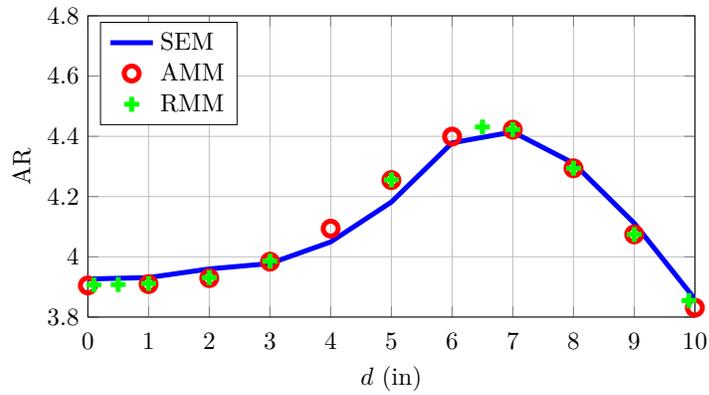

(a) Amplitude ratio

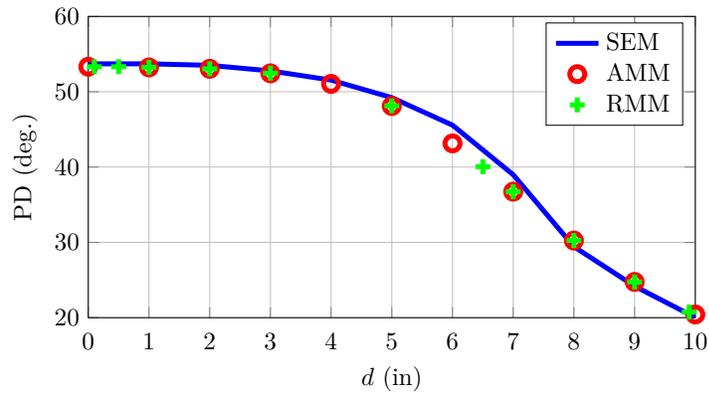

(b) Phase difference

Figure 5.5: Voltage AR and PD for the logging-well tool with an oil-invaded zone as a function of the invasion thickness $d$.

## 5.5
## Final Remarks

We have presented a new approach for modeling well-logging sensors via the fields obtained by higher-order SEM in cylindrical coordinates. Alternative SEM formulations found in the literature [21, 66] are less representative than the one presented here, as those works often consider simpler antenna models, such as magnetic dipoles. Preliminary results presented in this chapter align with those obtained through the mode-matching method. Our technique offers the advantage of modeling inhomogeneous oil invasions, which is not feasible with standard mode-matching solutions.

Continuing this research, we suggest integrating models for transmitting and receiving tilted-coil antennas into the SEM framework outlined here. Such an enhancement will facilitate the precise characterization of advanced logging sensors, paving the way for a more nuanced understanding of inhomogeneous and anisotropic geophysical formations. This endeavor extends the applicability of our SEM formulation and promises to enrich the toolkit available for exploring complex subsurface formations.

# 6
# Conclusion

## 6.1
## Conclusions and Future Works

The electromagnetic analysis of cylindrical-layered structures filled with complex anisotropic media typically involves employing computational software utilizing techniques like finite differences or finite elements. Although these methods are relatively easy to implement, they require extremely fine domain discretization and a significant number of degrees of freedom (DoF). In contrast, the suggested higher-order spectral element method (SEM) boasts remarkable accuracy while requiring fewer elements and DoFs, rendering it an attractive alternative. These characteristics stem from the exponential convergence of spectral methods and the simplification of field expansions in cylindrical coordinates, which is particularly beneficial for waveguides featuring cylindrical-conforming boundaries.

In Chapter 2, we investigated an SEM-based technique for solving the modal fields in cylindrical waveguides considering radially-bounded and radially-unbounded domains, homogeneous and inhomogeneous filled waveguides, lossless and lossy problems, and isotropic and anisotropic media. The results obtained showed excellent agreement with analytical results from the literature and FEM reference results. The cylindrical-coordinate-based SEM requires fewer elements and DoFs when compared with the conventional Cartesian-coordinate-based SEM and FEM.

In Chapter 3, the work investigated the SEM formulation to rapidly analyze the modal fields in eccentric coaxial waveguides. Given the geometry of the problem, the proposed normalized variational formulation obtains the solution for an entire set of configurations independently of the medium's parameters. Unlike perturbation-based techniques, our approach does not restrict the eccentricity offset and presents great accuracy and time-cost results for extreme geometric cases. The formulation accurately models a wide range of eccentric coaxial waveguides filled with lossy uniaxially anisotropic media. The results demonstrated excellent agreement with a dense-mesh reference FIT solution at a much reduced computational cost.

In Chapter 4, the study presented an improved spectral element method for analyzing eccentric two-layer waveguides with anisotropic non-reciprocal media. The complexity of the eccentric waveguide was transformed into a more tractable concentric counterpart through the application of conformal transformation optics. The cylindrical SEM was applied in the concentric domain and proved accurate to model waveguides filled with non-reciprocal media for a wide frequency range. Our approach required far fewer DoF in comparison to the conventional FEM solution.

In Chapter 5, we have presented a new approach for modeling well-logging sensors via the fields obtained by higher-order SEM in cylindrical coordinates. Preliminary results agree with the mode-matching method for the amplitude ratio and phase difference voltage. Our technique can be applied to model inhomogeneous oil invasions, which is not possible with the standard mode-matching formulation.

In future research, we plan to explore the convergence analysis of solutions obtained from variational waveguide problems for electric and magnetic fields as a function of expansion order. In geophysical applications, we aim to incorporate tilted coil antennas to evaluate anisotropic formations and consider eccentric oil invasion using transformation optics. Finally, we'd like to utilize the fields obtained via SEM as a basis for a mode-matching method, allowing us to model geophysical formations with axial layers more efficiently.

## 6.2
## Published Works

The author has contributed the following published works to the literature:

– Raul O. Ribeiro, Guilherme S. Rosa, José R. Bergmann, and Fernando L. Teixeira, "A Higher-Order Spectral Element Method to Model Eccentric Anisotropic Two-Layer Waveguides via Conformal Transformation Optics", *2024 18th European Conference on Antennas and Propagation (EuCAP)*, to appear.

– Raul O. Ribeiro, Johnes R. Gonçalves, Guilherme S. Rosa, José R. Bergmann, and Fernando L. Teixeira, "Higher Order Spectral Element Method for Rapid Analysis of Eccentric Coaxial Waveguides Filled with Lossy Anisotropic Media", *IEEE Transactions on Microwave Theory and Techniques*, Nov 22, 2023, doi: 10.1109/TMTT.2023.3332249.

– Raul O. Ribeiro, Guilherme S. Rosa, José R. Bergmann, and Fernando L. Teixeira, "Spectral Element Method for Modeling Anisotropic Circular Waveguides Loaded with Eccentric Rods via Conformal Transformation Optics", *2023 IEEE International Symp. on Antennas and Propag. and USNC-URSI Radio Sci. Meeting*, Sep 7, 2023, doi: 10.1109/USNC-URSI52151.2023.10237726.

– Raul O. Ribeiro, Johnes R. Gonçalves, Fernando L. Teixeira, José R. Bergmann and Guilherme S. Rosa, "Spectral Element Method for Modeling Eccentric Coaxial Waveguides Filled with Anisotropic Media via Conformal Transformation Optics", *2023 IEEE/MTT-S International Microwave Symposium - IMS 2023*, Jul 28, 2023, doi: 10.1109/IMS37964.2023.10187907.

– Raul O. Ribeiro, Guilherme S. Rosa, José R. Bergmann, and Fernando L. Teixeira, "A Novel High-Order Spectral Element Method for the Analysis of Cylindrical Waveguides Filled with Complex Anisotropic Media", *2023 17th European Conference on Antennas and Propagation (EuCAP)*, May 31, 2023, doi: 10.23919/EuCAP57121.2023.10133096.

– Raul O. Ribeiro, Guilherme S. Rosa, and José R. Bergmann, "An improved formulation of the spectral element method in cylindrical coordinates for the analysis of anisotropic-filled waveguide devices", *2022 IEEE International Symposium on Antennas and Propagation and USNC-URSI Radio Science Meeting (APS/URSI)*, Sep 21, 2022, doi: 10.1109/AP-S/USNC-URSI47032.2022.9886414.

## 6.3
**Grants and Awards**

The author received the following grants and awards:

– 2023 Region 9 Student Travel Grant for the IEEE International Symposium on Antennas and Propagation and USNC-URSI Radio Science Meeting, Portland, OR, USA.

– Finalist of the 2023 EuCAP (European Conference on Antennas and Propagation) Best Student Paper, Florence, Italy.

– 2022 IEEE Antennas and Propagation Society Fellowship (APSF).

– 2022 Doctoral Sandwich Abroad Fellowship by CNPq of Brazil's Ministry of Science, Technology, and Innovation.

– 2022 Region 9 Student Travel Grant for the IEEE International Symposium on Antennas and Propagation, Denver, CO, USA.

– Honorable Mention, 2022 IEEE AP-S Student Paper Competition, Denver, CO, USA.

# A

# Eccentric Coaxial Waveguide Solution via Graf's Addition Theorem

The axial fields that solve the problem in Fig. 3.1(a) can be expressed in terms of symmetrical (concerning the $\tilde{x}$ axis and identified via $s = 0$) and antisymmetric (identified via $s = 1$) configurations. Accordingly, the field solutions can be written as a truncated Bessel-Fourier series [43, Ch. 5] in the eccentric coordinate system via

$$\tilde{F}(\tilde{\rho}, \tilde{\phi}) = \sum_{j=1}^{N_G+1-s} c_j^{\text{TX}} \, B_{m,m}^{\text{TX}}(k_\rho^{\text{TX}} \tilde{\rho}) \, \Phi(m\tilde{\phi}), \tag{A-1}$$

where TX $\rightarrow$ TM for $\tilde{F} = \tilde{E}_z$, TX $\rightarrow$ TE for $\tilde{F} = \tilde{H}_z$, with $m = j - 1 + s$,

$$B_{m,n}^{\text{TM}}(k_\rho^{\text{TM}} \tilde{\rho}) = \left[ Y_m(k_\rho^{\text{TM}} \tilde{r}_1) \, J_n(k_\rho^{\text{TM}} \tilde{\rho}) - J_m(k_\rho^{\text{TM}} \tilde{r}_1) \, Y_n(k_\rho^{\text{TM}} \tilde{\rho}) \right], \tag{A-2}$$

$$B_{m,n}^{\text{TE}}(k_\rho^{\text{TE}} \tilde{\rho}) = \left[ Y_m'(k_\rho^{\text{TE}} \tilde{r}_1) \, J_n(k_\rho^{\text{TE}} \tilde{\rho}) - J_m'(k_\rho^{\text{TE}} \tilde{r}_1) \, Y_n(k_\rho^{\text{TE}} \tilde{\rho}) \right], \tag{A-3}$$

and $\Phi(m\tilde{\phi}) = \cos(m\tilde{\phi} - s\pi/2)$. In the above, Bessel and Newman functions of order $m$ are denoted by $J_m(\cdot)$ and $Y_m(\cdot)$, respectively, and the prime $'$ represents the derivative of the cylindrical functions with respect to their arguments. Also, we have $k_\rho^{\text{TM}} = (\tilde{\epsilon}_z/\tilde{\epsilon}_s)^{1/2} \, k_\rho$ and $k_\rho^{\text{TE}} = (\tilde{\mu}_z/\tilde{\mu}_s)^{1/2} \, k_\rho$. The axial wavenumber is again obtained via (3-22).

The field expansion in (A-1) already satisfy the boundary conditions at $\tilde{\rho} = \tilde{r}_1$. On the inner PEC at $\tilde{\rho} = \tilde{r}_0$, the remaining boundary conditions can be enforced via a coordinate transformation from a polar coordinate system with the origin at $(\tilde{x} = \tilde{d}, \tilde{y} = 0)$ to the origin $(\tilde{x} = 0, \tilde{y} = 0)$ and then by using the Graf's addition theorem [84, p. 363], [52]. As result, by recasting all the unknown coefficients $c_j^{\text{TX}}$ into the column-vector $\bar{c}^{\text{TX}}$, we can obtain the homogeneous system of equations

$$\mathbf{B}^{\text{TX}} \, \bar{c}^{\text{TX}} = \bar{0}, \tag{A-4}$$

where the associated matrices $\mathbf{B}^{\mathrm{TX}}$ for TM and TE modes are given by

$$\mathbf{B}^{\mathrm{TM}}(m+1-s, n+1-s) = B_{m,n}^{\mathrm{TM}}(k_\rho^{\mathrm{TM}}\tilde{r}_0)\, D_{m,n}(k_\rho^{\mathrm{TM}}\tilde{d}), \qquad \text{(A-5)}$$

$$\mathbf{B}^{\mathrm{TE}}(m+1-s, n+1-s) = B_{m,n}^{\mathrm{TE}}(k_\rho^{\mathrm{TE}}\tilde{r}_0)\, D_{m,n}(k_\rho^{\mathrm{TE}}\tilde{d}), \qquad \text{(A-6)}$$

with $m = s, s+1, \ldots, N_G$, $n = s, s+1, \ldots, N_G$, and

$$D_{m,n}(k_\rho^{\mathrm{TX}}\tilde{d}) \;=\; \left[ J_{n-m}(k_\rho^{\mathrm{TX}}\tilde{d}) + (-1)^{(m+s)}\, J_{n+m}(k_\rho^{\mathrm{TX}}\tilde{d}) \right]\, f_J(k_\rho^{\mathrm{TX}}\tilde{d}). \quad \text{(A-7)}$$

In the above, we have introduced the scaling factor $f_J(x) = (\pi x/2)^{1/2}$.

The non-trivial solution of (A-4) is given by $\det(\mathbf{B}^{\mathrm{TX}}) = 0$, and the radial wavenumbers $k_\rho^{\mathrm{TX}}$ are obtained by tracking the zeros of this transcendental equation. We have used Muller's method for finding the zeros using the algorithm described in [85, p. 466] and [86, App. E]. After that, for each found $k_\rho^{\mathrm{TX}}$, we then use the singular value decomposition of matrix $\mathbf{B}^{\mathrm{TX}}$ to compute the associated eigenvector $\bar{c}^{\mathrm{TX}}$ and obtain the corresponding axial field via (A-1).

The above formulas are a generalized version of that in [52] adequate to handle lossy uniaxially anisotropic media. It is important to emphasize that numerical issues on the zero-finding process will occur if the original formulas in [52] are used because of the presence of poles in such characteristic equations. The presented forms in (A-5) and (A-6) prevent such problems. In addition, the presence of the scaling function $f_J(x)$ alleviates other numerical stability issues related to the employment of a large number of cylindrical harmonics, e. g., $N_G$ large than 10 in (A-1). When $\tilde{d}/\tilde{r}_1$ is small, the convergence of the fundamental modes obtained using the GAT is archived using a relatively small $N_G$. Still, for large $\tilde{d}/\tilde{r}_1$, the procedure described here may present numerical issues due to the need of high orders cylindrical functions that cause under- and over-flow on the evaluation of $\det(\mathbf{B}^{\mathrm{TX}})$.

We have used $N_G = 4$ for all the GAT computations whose results are presented in Section 3.4.1. Note that the number of harmonics $N_G$ in the series (A-1) is not directly related to the expansion order $N$ of the SEM solution.

# B
# SEM for Modeling Anisotropic Circular Waveguides Loaded with Eccentric Rods

## B.1
## Introduction

The analysis of field propagation involving complex media with non-concentric radial layers is essential to a large class of problems including metamaterial devices and geophysical exploration. Several approaches have been developed to compute wavenumbers and field patterns in these structures [10,46,75,87]. Recently, an approximated method in a bipolar coordinate system was used to obtain the cutoff frequencies of circular waveguides containing an offset in the center conductor [61] or with an eccentrically dielectric-loaded [73]. Similar problems were solved in [12] and [74] using transformation optics principles to map the non-concentric problem into an equivalent concentric problem and perturbation techniques to compute the axial wavenumbers.

The spectral element method (SEM) has been applied successfully to various electromagnetic problems [18, 19, 22]. Preliminary results [22] suggest that a SEM formulation based on cylindrical coordinates requires considerably fewer resources than the conventional Cartesian approach for cylindrical-conforming boundaries. In this work, we propose a new SEM formulation based on transformation optics to map a circular waveguide problem loaded with an eccentric dielectric rod into an equivalent concentric problem. This problem was tackled by perturbation techniques in [74], which are limited to small eccentricity offsets. Here, the SEM in cylindrical coordinates enables the analysis of structures with larger offsets.

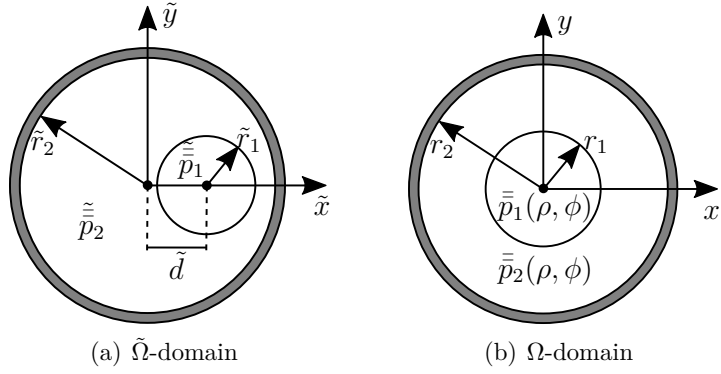

Figure B.1: Cross-section of a circular waveguide loaded with an (a) eccentric
rod, and (b) a mapped concentric rod.

## B.2
## Formulation Overview

Consider a circular waveguide invariant along the axial direction, radially
bounded by a perfect electric conductor (PEC) of radius $\tilde{r}_2$, and loaded with an
inner circular dielectric cylinder of radius $\tilde{r}_1$ with an offset $\tilde{d}$ from $\tilde{z}$-axis. The
eccentric structure is referred to as a cylindrical coordinate system denoted by
$(\tilde{\rho}, \tilde{\phi}, \tilde{z})$. The associated cross-section $\tilde{\Omega}$ is depicted in Fig. B.1(a). We assume
the time-harmonic dependence $e^{-i\omega t}$, and that the waveguide is filled with
uniaxially anisotropic mediums characterized by real-valued permeability

$$\bar{\bar{\tilde{\mu}}}_j = \mu_0 \bar{\bar{\tilde{\mu}}}_{r,j} \ \text{ with } \ \bar{\bar{\tilde{\mu}}}_{r,j} = \mathrm{diag}(\tilde{\mu}_{rs,j}, \tilde{\mu}_{rs,j}, \tilde{\mu}_{rz,j}), \tag{B-1}$$

and the complex-valued permittivity

$$\bar{\bar{\tilde{\epsilon}}}_j = \epsilon_0 \bar{\bar{\tilde{\epsilon}}}_{r,j} + \frac{i}{\omega} \bar{\bar{\tilde{\sigma}}}_j \tag{B-2}$$

with

$$\bar{\bar{\tilde{\epsilon}}}_{r,j} = \mathrm{diag}(\tilde{\epsilon}_{rs,j}, \tilde{\epsilon}_{rs,j}, \tilde{\epsilon}_{rz,j}), \ \bar{\bar{\tilde{\sigma}}}_j = \mathrm{diag}(\tilde{\sigma}_{s,j}, \tilde{\sigma}_{s,j}, \tilde{\sigma}_{z,j}), \tag{B-3}$$

where $j = 1, 2$ refer to the inner (dielectric rod) and outer regions, respectively.

From conformal transformation optics [10,36], the original problem in the
eccentric coordinates $(\tilde{\rho}, \tilde{\phi}, \tilde{z})$ can be transformed into a concentric problem
with the coordinates $(\rho, \phi, z)$ using

$$\mathbf{F} = \bar{\bar{J}} \cdot \tilde{\mathbf{F}}, \ \text{ with } \ \mathbf{F} \in \{\mathbf{E}, \mathbf{H}\}, \tag{B-4}$$

$$\bar{\bar{p}}_j = |\bar{\bar{J}}|^{-1} \bar{\bar{J}} \cdot \bar{\bar{\tilde{p}}}_j \cdot \bar{\bar{J}}^T, \ \text{ with } \ p_j \in \{\mu_j, \epsilon_j\}, \tag{B-5}$$

where $\bar{\bar{J}}$ is the Jacobian of the transformation $(\tilde{\rho}, \tilde{\phi}, \tilde{z}) \to (\rho, \phi, z)$. Assuming
$z = \tilde{z}$ and $r_2 = \tilde{r}_2$, we can express the transformed constitutive tensors as [10]

$$\bar{\bar{p}}_j = \mathrm{diag}(p_{s,j}, p_{s,j}, p_{z,j}(\rho, \phi)) = \mathrm{diag}(\tilde{p}_{s,j}, \tilde{p}_{s,j}, |\bar{\bar{J}}|^{-1}\tilde{p}_{z,j}), \tag{B-6}$$

with $p_j \in \{\mu_j, \epsilon_j\}$, $j \in \{1, 2\}$, where

$$|\bar{\bar{J}}|^{-1} = \frac{(1 - \tilde{x}_1/\tilde{x}_2)^2}{(1 - 2\rho\cos\phi/\tilde{x}_2 + \rho^2/\tilde{x}_2^2)^2}, \tag{B-7}$$

$$\tilde{x}_{1,2} = \frac{-c \mp \sqrt{c^2 - 4\tilde{r}_2^2}}{2}, \quad \text{and} \quad r_1 = \tilde{r}_1 \left( \frac{\tilde{x}_2}{\tilde{x}_2 - \tilde{d}} \right). \tag{B-8}$$

Note that in the mapped domain the axial parameters $\mu_{z,j}$ and $\epsilon_{z,j}$ vary in the
$\rho$ and $\phi$ directions.

The concentric problem in the $\Omega$-domain (depicted in Fig. B.1(b) can
be solved using the variational formulation in [18] associated with the recently
introduced SEM technique discussed in [22].

## B.3
## Results and Discussion

As validation, we consider an anisotropic circular waveguide with $\tilde{r}_2 =$
10 mm, $\bar{\bar{\tilde{\epsilon}}}_{r,2} = \text{diag}\{2.7, 2.7, 2.0\}$, $\bar{\bar{\tilde{\sigma}}}_2 = \text{diag}\{0.12, 0.12, 0.15\}$, loaded with an
eccentric dielectric rod having $\tilde{r}_1 = 0.3\tilde{r}_2$, $\tilde{d} = 0.3\tilde{r}_2$, $\bar{\bar{\tilde{\epsilon}}}_{r,1} = \text{diag}\{5.6, 5.6, 4.6\}$,
$\bar{\bar{\tilde{\sigma}}}_1 = \text{diag}\{0.38, 0.38, 0.34\}$. Other non-declared constitutive parameters are
assumed equal to the vacuum ones. A similar geometry was considered in [73,
74]; however, [73] assumes isotropic media, and [74] assumes lossless media.
Fig. B.2 shows the real and imaginary part of axial wavenumber $k_z$ for the
$\text{TM}_{01}$ and $\text{HE}_{11}$ modes as a function of the operating frequency $f$ varying
from 1 GHz to 10 GHz. In view of (B-2), for $f = 1$ GHz the conduction and
displacement current have the same order of magnitude, while for $f = 10$ GHz
the displacement current is ten times greater. The results obtained with
the proposed SEM show excellent agreement with those from the finite-
element method (FEM) from COMSOL Multiphysics [2]. The computational
domain was discretized into 12 elements with expansion functions of order 6,
which implies 744 DoF. In the FEM COMSOL model, we use a cubic-order
discretization with the extremely fine physics-controlled mesh option to obtain
accurate reference solutions. This results in 19627 DoF. The number of DoFs
required by the present SEM formulation is 3.8 % of that required by the FEM
solution.

In conclusion, this work presents an efficient and accurate technique for
analyzing complex circular waveguides filled with eccentric dielectric rods. Our
formulation is not restricted to small eccentricities like perturbative-based
techniques and requires significantly fewer DoF compared to the reference
FEM solution. With a few adaptations, this formulation can provide a very
efficient tool for the study of wave propagation in complex media, with potential applications in the modeling of metamaterial devices and geophysical exploration sensors.

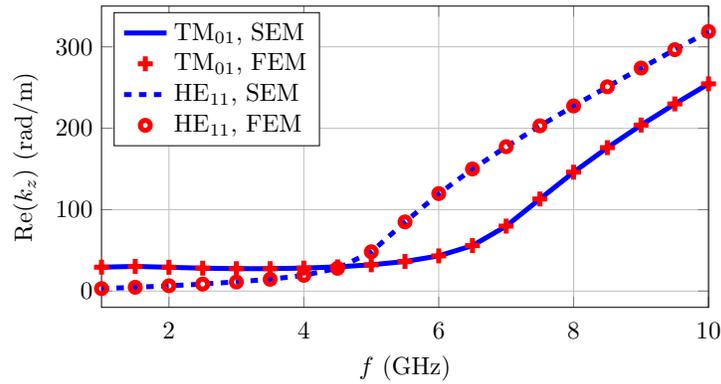

(a) Phase constant

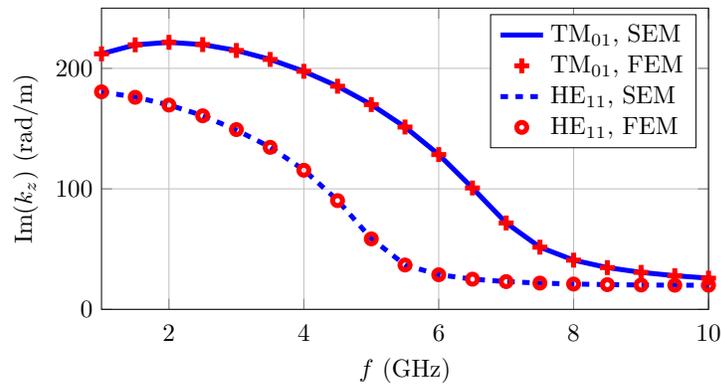

(b) Attenuation constant

Figure B.2: Axial wavenumbers for the two first modes in a circular waveguide loaded with an eccentric rod.